\documentclass{article}
\usepackage[english]{babel}
\usepackage{a4wide}

\RequirePackage{amsmath,amsfonts,amssymb,amsthm,mathtools,bm,siunitx}
\RequirePackage{graphicx,xcolor}
\graphicspath{{./fig/}}
\usepackage{ulem}

\newcommand{\sib}[1]{[\si{#1}]}

\theoremstyle{plain}
\newtheorem{problem}{Problem}
\newtheorem{assumption}{Assumption}

\usepackage[hidelinks]{hyperref}

\title{A mixed-dimensional model for direct current simulations in presence of a thin high-resistivity liner}
\author{
Alessio Fumagalli$^1$ \and
Lorenzo Panzeri$^1$ \and
Luca Formaggia$^1$ \and
Anna Scotti$^1$ \and
Diego Arosio$^2$}

\date{$^1$ Politecnico di Milano\\%
$^2$ Universit\`a di Modena e Reggio Emilia}

\begin{document}
\maketitle

\begin{abstract}
    In this work we present a mixed-dimensional mathematical model to obtain the electric potential and current density in direct current simulations when a thin liner is included in the modelled domain. The liner is used in landfill management to prevent leakage of leachate from the waste body into the underground and is made of a highly-impermeable high-resistivity plastic material.
    The electrodes and the liner have diameters and thickness respectively that
    are much smaller than their other dimensions, thus their numerical
    simulation might be too costly in an equi-dimensional setting. Our
    approach is to approximate them as objects of lower dimension and derive the corresponding equations. The obtained mixed-dimensional model is validated against laboratory experiments of increasing complexity showing the reliability of the proposed mathematical model.
\end{abstract}

\section{Introduction}

The direct current (DC) resistivity method is a geophysical technique used for
non-invasive investigations in a wide range of applications and is best suited
to detect subsurface saturation and pollution, as underground fluids and
dissolved ions strongly control subsurface resistivities \cite{Binley2005}.

Our target application is the geoelectrical monitoring of municipal solid waste
landfills (MSWLFs). In more detail, we want to explore the capability of the DC
method to characterize the high- or low-density polyethylene (HDPE or LDPE)
liner that is generally placed underneath the landfill to prevent leachate
leakage in the underground. The method is promising to discriminate between the
plastic liner, that is practically an electrical insulator, from leachate that
can have resistivity as low as 1 \sib{\ohm\meter}. The scientific literature reports several
case studies where geoelectrical methods have been applied to landfills with
different design, size, and confining materials, consisting of diverse types of
waste, and located either in urbanized or remote areas
\cite{Bernstone2000,MEJU2000115,Leroux2010,VAUDELET2011738,DeCarlo2013,Tsourlos2014,DIMAIO2018629}.
In many cases, the aim of the
studies was to delineate the extent of the buried waste body and DC measurements
have been conducted together with Induced Polarization ones
(\cite{Dahlin2010,DEDONNO2017302}). Some numerical and down-scaled laboratory
studies have been devoted to detecting and locating ruptures and holes in the
plastic liner with set ups with different levels of complexity
(\cite{Frangos1997,Binley1997,Binley2003,Ling2019,Aguzzoli2020}).
However, the geoelectrical method aimed to analyse the integrity of the
geomembrane at the field scale still lacks thorough validation and issues
related to the survey design, the investigation depth and to the spatial
resolution need to be addressed (\cite{aguzzoli2021inversion}).

In the simplest DC configuration, two electrodes inject current $I$ in the
investigated medium and two electrodes measure the potential difference $\Delta
V$.
The ratio $\Delta V / I$  multiplied by a factor $K$ related to the geometrical arrangement
of the electrodes yields the apparent electrical resistivity $\rho_a$, which is a
function of the resistivities of the investigated media and changes according to
$K$.

Up to date, geoelectrical equipment easily handles tens to hundreds of
electrodes for 2D, 3D and 4D (i.e., 3D-time lapse) surveys (\cite{Loke2019})
allowing for the tomographic reconstruction of the subsurface resistivities
(electrical resistivity tomography - ERT). ERT constitutes a non-linear and
ill-posed inverse problem in which the non-linearity is due to the equation
linking the subsurface resistivity values to the measured apparent resistivity
data, whereas the ill-posedness is mainly related to incomplete data coverage,
considerable number of unknown parameters to be estimated and noise
contamination in the data (\cite{Menke2018}).
Whether the inversion of geoelectrical
data is performed with a deterministic or probabilistic approach, a robust and
effective numerical code is needed to solve the forward problem. Though both
commercial (e.g., \cite{Loke2022}) and open-source forward geoelectrical codes
(\cite{RUCKER2017106,BLANCHY2020104423})
have been developed, in this work we
propose a mathematical model taking into account the peculiar aspects of the
geoelectrical monitoring of MSWLFs related to the different sizes of the
elements to be considered. The liner is a few millimetres thick and has lateral
extent of tens of thousands square meters; similarly, the ratio of the diameter
to the length of the electrodes is typically one to forty centimetres (
\cite{Ruecker2011}). It is clear that representing these items as three-dimensional
objects requires a computational grid that might be unnecessary refined, as
discussed in \cite{Berre2020a} for a similar context, with obvious computational cost.
Accordingly, we introduce new conceptual models to lighten the computational
cost of the forward problem while maintaining the accuracy. We consider a model
reduction technique to approximate the electrodes as one-dimensional objects of
the same length as their three-dimensional representation, and the liner as
two-dimensional object of the same extent. Example of model reduction to
one-dimensional objects can be found, for instance, in \cite{Peaceman1978,Peaceman1983}
and more
recently in \cite{Cerroni2019,Gjerde2019,Gjerde2020,Kuchta2021,Berrone2022}
while for two-dimensional
objects in
\cite{Martin2005,DAngelo2011,Fumagalli2016a,Scialo2017,Nordbotten2018}.
A new set of equations will be
derived for the electrodes and the liner along
with appropriate interface conditions for their coupling with the surrounding
domain. With the reduced model, a coarser grid can be used, or bigger problems
can be solved, which will reduce ill-posedness and computational cost.

The outcomes of the model-based code are then compared with results obtained
through down-scaled experiments in the laboratory and, whenever possible, with
theoretical results. In particular, we consider several settings by changing the
depth of the geomembrane and the presence or absence of holes.

To solve the model numerically two cell-centred finite volume schemes, multi-point flux approximation (MPFA) \cite{Aavatsmark2007,Aavatsmark2002} and two-point flux
approximation (TPFA) \cite{Aavatsmark2007a} methods, are exploited and the results are
compared with experimental data to identify the most suitable approach in terms
of computational costs and accuracy.

This article is organized as follows. In Section \ref{sec:model}, first we introduce the
mathematical model for a three-dimensional domain, and then the
mixed-dimensional model for the electrodes and the liner. Section \ref{sec:approximation} is devoted
to the spatial discretization of the equations by discussing appropriate
strategies to make their solution more efficient. In Section \ref{sec:experiments} we consider two
different settings to validate the proposed model: in the first one the liner is
flat and changes its depth, while in the second one the liner has a box shape
with the possible presence of a hole. Finally, we draw the conclusions in
Section \ref{sec:conclusions}.

\section{The mathematical model}\label{sec:model}

We
indicate
the domain where the geoelectrical equations will be applied with
$\Upsilon\subset \mathbb{R}^3$, with boundary $\partial \Upsilon$ and outward
unit normal $\bm{n}_\partial$. We consider the
following main variables: $\bm{E}:\Upsilon \rightarrow \mathbb{R}^3$, the electric
field in \sib{\volt\per\meter}, $\bm{J}:\Upsilon \rightarrow \mathbb{R}^3$,
the current density field in \sib{\ampere\per\square\meter}, and
$\varphi: \Upsilon \rightarrow \mathbb{R}$, the electric potential
in \sib{\volt}.

We present now the governing equations.
Ohm-Kirchhoff's law states that $\bm{J} = \sigma \bm{E}$,
with $\sigma:\Upsilon \rightarrow \mathbb{R}_{>0}$ the electric
conductivity of the medium in
\sib{\siemens\per\metre}. The latter is the inverse of the electric resistivity
$\rho$, measured in \sib{\ohm\metre}.
Moreover, as a consequence of the static Maxwell-Faraday equation, $\nabla \times
\bm{E} =\bm{0}$, and by the Helmholtz decomposition we get that the electric field
$\bm{E}$ can be related to the
gradient of the electric potential $\varphi$ as $\bm{E} = -\nabla \varphi$.
We can thus establish a constitutive relation between the current density field $\bm{J}$ and
the electric potential $\varphi$ as follow
\begin{gather}\label{eq:current_potential}
    \bm{J} + \sigma \nabla \varphi = \bm{0}.
\end{gather}
Gauss' law, or charge conservation equation for $\bm{J}$, can be expressed as
\begin{gather}\label{eq:continuity}
    \nabla \cdot \bm{J} = q,
\end{gather}
with $q:\Upsilon \rightarrow \mathbb{R}$ being the source of volumetric charge
density in \sib{\ampere\per\cubic\meter}.
For simplicity, we
consider an isolated body so at the boundary $\partial \Upsilon$ we set
$\bm{J}\cdot\bm{n}_\partial = 0$.
By combining \eqref{eq:current_potential} and \eqref{eq:continuity} the
problem can be written as:
\begin{problem}[The DC problem]\label{eq:the_system}
    The direct current problem is: find $(\bm{J}, \varphi)$ in $\Upsilon$ such that
    \begin{gather*}
        \begin{aligned}
            &\begin{aligned}
                &\bm{J} + \sigma \nabla \varphi = \bm{0}\\
                &\nabla \cdot \bm{J} = q
            \end{aligned}
            && \text{in } \Upsilon\\
            &\bm{J} \cdot \bm{n}_\partial = 0
            && \text{on } \partial \Upsilon
        \end{aligned}
    \end{gather*}
\end{problem}
In Subsection \ref{subsec:electrode}, we present a suitable model for the electrodes
represented as objects of dimension 1 immersed in $\Upsilon$ and in Subsection
\ref{subsec:liner} the bi-dimensional model for the liner.

\subsection{A model for the electrodes}\label{subsec:electrode}

Depending on the chosen electrode configuration, multiple electrodes are immersed in
$\Upsilon$
and can be used, in pairs, to inject electrical current and measure
the electric potential differences.

We want to introduce an effective and efficient mathematical model able to
describe their influence in Problem \ref{eq:the_system}.  For simplicity, we
consider a single electrode, being immediate the extension to multiple (separated)
electrodes, and set $q=0$.  Even if an electrode is a three-dimensional object,
typically
a cylinder of stainless steel, its radius is very small compared to its length
and domain size. Typically, the ratio between its radius and length is
in the order of 1/30 $\sim$ 1/50. Our idea is thus to approximate the electrode by an object of dimension
1, making the gridding process much easier and with fewer cells, essential
when dealing with real geometries.

\begin{figure}[tb]
    \centering
    \resizebox{0.25\textwidth}{!}{\fontsize{0.75cm}{2cm}\selectfont
\begingroup%
  \makeatletter%
  \providecommand\color[2][]{%
    \errmessage{(Inkscape) Color is used for the text in Inkscape, but the package 'color.sty' is not loaded}%
    \renewcommand\color[2][]{}%
  }%
  \providecommand\transparent[1]{%
    \errmessage{(Inkscape) Transparency is used (non-zero) for the text in Inkscape, but the package 'transparent.sty' is not loaded}%
    \renewcommand\transparent[1]{}%
  }%
  \providecommand\rotatebox[2]{#2}%
  \newcommand*\fsize{\dimexpr\f@size pt\relax}%
  \newcommand*\lineheight[1]{\fontsize{\fsize}{#1\fsize}\selectfont}%
  \ifx\svgwidth\undefined%
    \setlength{\unitlength}{272.250359bp}%
    \ifx\svgscale\undefined%
      \relax%
    \else%
      \setlength{\unitlength}{\unitlength * \real{\svgscale}}%
    \fi%
  \else%
    \setlength{\unitlength}{\svgwidth}%
  \fi%
  \global\let\svgwidth\undefined%
  \global\let\svgscale\undefined%
  \makeatother%
  \begin{picture}(1,0.83471096)%
    \lineheight{1}%
    \setlength\tabcolsep{0pt}%
    \put(0,0){\includegraphics[width=\unitlength,page=1]{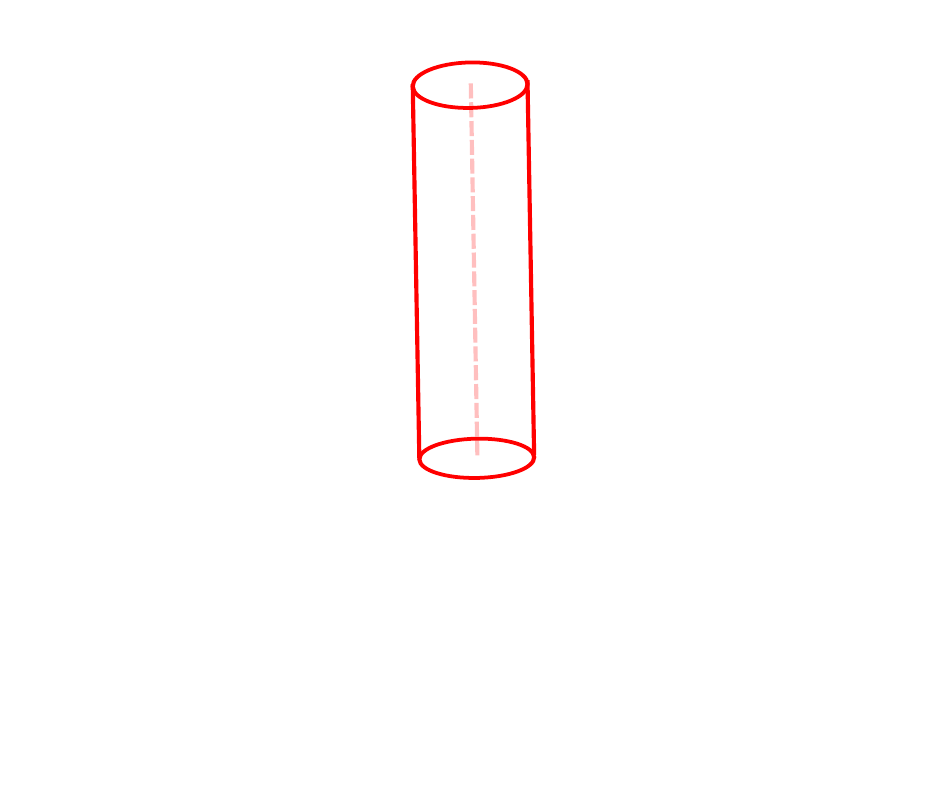}}%
    \put(0.03935454,0.56473752){\color[rgb]{0,0,0}\makebox(0,0)[lt]{\lineheight{1.25}\smash{\begin{tabular}[t]{l}$\Upsilon$\end{tabular}}}}%
    \put(0,0){\includegraphics[width=\unitlength,page=2]{electrode_full.pdf}}%
    \put(0.56800123,0.5329003){\color[rgb]{1,0,0}\makebox(0,0)[lt]{\lineheight{1.25}\smash{\begin{tabular}[t]{l}$\Gamma$\end{tabular}}}}%
    \put(0.51093711,0.28988601){\color[rgb]{0,0,0}\makebox(0,0)[lt]{\lineheight{1.25}\smash{\begin{tabular}[t]{l}$r$\end{tabular}}}}%
  \end{picture}%
\endgroup%
}%
    \hspace{0.1\textwidth}%
    \resizebox{0.25\textwidth}{!}{\fontsize{0.75cm}{2cm}\selectfont
\begingroup%
  \makeatletter%
  \providecommand\color[2][]{%
    \errmessage{(Inkscape) Color is used for the text in Inkscape, but the package 'color.sty' is not loaded}%
    \renewcommand\color[2][]{}%
  }%
  \providecommand\transparent[1]{%
    \errmessage{(Inkscape) Transparency is used (non-zero) for the text in Inkscape, but the package 'transparent.sty' is not loaded}%
    \renewcommand\transparent[1]{}%
  }%
  \providecommand\rotatebox[2]{#2}%
  \newcommand*\fsize{\dimexpr\f@size pt\relax}%
  \newcommand*\lineheight[1]{\fontsize{\fsize}{#1\fsize}\selectfont}%
  \ifx\svgwidth\undefined%
    \setlength{\unitlength}{272.250359bp}%
    \ifx\svgscale\undefined%
      \relax%
    \else%
      \setlength{\unitlength}{\unitlength * \real{\svgscale}}%
    \fi%
  \else%
    \setlength{\unitlength}{\svgwidth}%
  \fi%
  \global\let\svgwidth\undefined%
  \global\let\svgscale\undefined%
  \makeatother%
  \begin{picture}(1,0.83471096)%
    \lineheight{1}%
    \setlength\tabcolsep{0pt}%
    \put(0,0){\includegraphics[width=\unitlength,page=1]{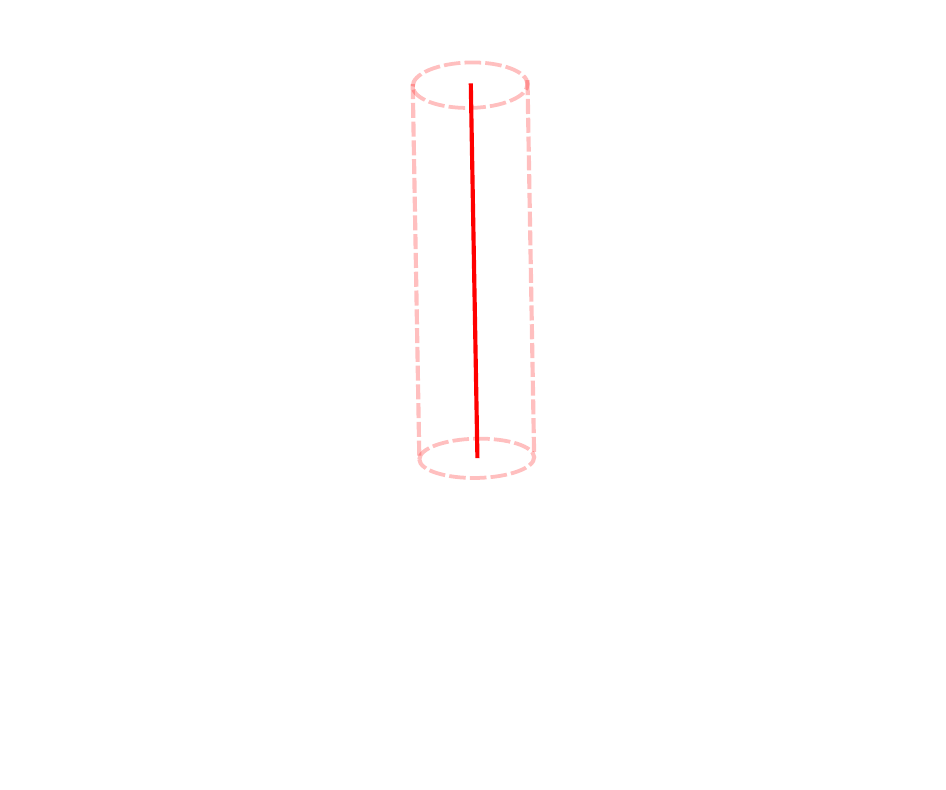}}%
    \put(0.03935454,0.56473752){\color[rgb]{0,0,0}\makebox(0,0)[lt]{\lineheight{1.25}\smash{\begin{tabular}[t]{l}$\Upsilon$\end{tabular}}}}%
    \put(0,0){\includegraphics[width=\unitlength,page=2]{electrode_reduced.pdf}}%
    \put(0.50739519,0.5329003){\color[rgb]{1,0,0}\makebox(0,0)[lt]{\lineheight{1.25}\smash{\begin{tabular}[t]{l}$\gamma$\end{tabular}}}}%
  \end{picture}%
\endgroup%
}
    \caption{An electrode $\Gamma$ represented as three-dimensional object
    (right)
    and the reduced electrode $\gamma$ as one-dimensional object (left).}
    \label{fig:electrode}
\end{figure}

Let us consider Figure \ref{fig:electrode} on the left where a three-dimensional electrode
$\Gamma$ is immersed in $\Upsilon$, with $ \mathring{\Gamma} \cap
\mathring{\Upsilon} =\emptyset$. We assume that one end of the electrode is
in contact with the top of $\partial \Upsilon$.
We suppose $\Gamma$ to be a cylinder with
radius $r$ and longitudinal axis $\gamma$ of length $l$, described as
\begin{gather*}
    \Gamma = D \times \gamma
    \quad \text{with} \quad
    D = \{ (\xi \cos \theta, \xi \sin
    \theta) \text{ with } \xi \in [0, r) \text{ and } \theta \in [0, 2\pi) \}.
\end{gather*}
We indicate with $\bm{v}_\Gamma$ the unit vector aligned with line
$\gamma$. The equations of
Problem \ref{eq:the_system} apply to $\Gamma$ with coefficients that are different
than those of the surrounding domain $\Upsilon$. We
indicate variables in $\Gamma$ with a subscript.
To couple the two systems, suitable conditions are needed
on the common boundary between $\Upsilon$ and $\Gamma$, we call
$\partial_i \Gamma = \overline{\Gamma} \cap \overline{\Upsilon}$
the internal part of $\partial \Gamma$. Since
the electrode is a cylinder, we further divide
$\partial_i \Gamma$ into two disjoint portions $\partial_l \Gamma$ and
$\partial_b \Gamma$, the former being the lateral surface and the latter the bottom
one. Finally, $\partial_t \Gamma$ is the top part of the boundary of
$\Gamma$, which is not in contact with the interior of $\Upsilon$.
Conservation of the current density through $\partial_i \Gamma$ and continuity
of the electric potential are required,
\begin{gather*}
    \begin{aligned}
        &\bm{J} \cdot \bm{n}_\Gamma = \bm{J}_\Gamma \cdot \bm{n}_\Gamma\\
        &\varphi = \varphi_\Gamma
    \end{aligned}
    \quad \text{on } \partial_i \Gamma
\end{gather*}
where $\bm{n}_\Gamma$ is the unit normal of $\partial_i \Gamma$ pointing from
$\Upsilon$ toward $\Gamma$.

To avoid excessive refinement of the computational grid to
capture the three-dimensional nature of the electrode, we employ
a model reduction strategy by approximating $\Gamma$ with its centre line
$\gamma$ while the surrounding domain, with an abuse of notation, will be
still indicated with $\Upsilon$, see Figure \ref{fig:electrode} on the right.
We will derive a new set of equations for $\gamma$ and coupling conditions
between $\gamma$ and $\Upsilon$, by following the
approach presented in
\cite{Peaceman1978,Peaceman1983,DAngelo2008,Gjerde2019,Gjerde2020,Kuchta2021}
derived for wells in porous media and here adapted
to the electrodes. We assume the following:
\begin{assumption}[Electrode radius]
    We assume that the radius $r$ of the electrode is much smaller than any
    other characteristic length in the domain.
\end{assumption}
\begin{assumption}[Lateral exchange]\label{ass:lateral_exchange}
    We assume that the current density is exchanged between the electrode
    and the surrounding media through the lateral boundary of the electrode $\partial_l
    \Gamma$, the
    exchange through $\partial_b \Gamma$ is assumed to be negligible.
\end{assumption}
The line current density exchanged between the electrode, now represented as $\gamma$, and the
surrounding media is denoted as
$j_\gamma$
in \sib{\ampere\per\metre} and approximated by the relation
\begin{gather}\label{eq:varsigma_electrde}
    j_\gamma + \sigma_\gamma (\varphi_\gamma - \varphi) = 0
    \quad \text{on } \gamma,
\end{gather}
where $\sigma_\gamma$ in \sib{\siemens\per\metre} is the conductivity of the interface between the electrode
and the surrounding media, whose meaning will be discussed in the following.
The reduced variables for the electrode $\varphi_\gamma$ in \sib{\volt}
and $\bm{J}_\gamma$ now in
\sib{\ampere} are defined as
\begin{gather*}
    \varphi_\gamma(z) = \frac{1}{\pi r^2}\int_D
    \varphi_\Gamma(\xi, \theta, z) d \xi d \theta
    \quad \text{and} \quad
    \bm{J}_\gamma(z) = \int_D (\bm{v}_\Gamma \otimes \bm{v}_\Gamma)
    \bm{J}_\Gamma (\xi, \theta, z) d \xi d \theta.
\end{gather*}
By following the model reduction procedure considered in the aforementioned
works, we can introduce a new model for the electrode represented as an object
$\gamma$ of dimension $n-2$
\begin{gather}\label{eq:electrode_dc}
    \begin{aligned}
        &\bm{J}_\gamma + \pi r^2 \sigma_\Gamma \nabla \varphi_\gamma = \bm{0}\\
        &\nabla \cdot \bm{J}_\gamma - j_\gamma = 0
    \end{aligned}
    \quad \text{in } \gamma,
\end{gather}
where the gradient and divergence are now computed along $\gamma$.
The second equation is a conservation equation which accounts for the
current along $\gamma$ and its exchange with
$\Upsilon$, represented as $j_\gamma$.
On the base of
assumption \ref{ass:lateral_exchange}, we set the so-called tip condition on the
bottom boundary of $\gamma$
\begin{gather}\label{eq:tip_condition}
    \bm{J}_\gamma \cdot \bm{v}_\Gamma = 0
    \quad \text{on } \partial_b \gamma.
\end{gather}
As mentioned, the electrodes inject a prescribed current density
$\overline{J}_\Gamma$ in \sib{\ampere\per\square\meter}, or
measure potential difference, thus $\overline{J}_\Gamma = 0$ in this latter
case. Accordingly, we employ one of the following conditions on
$\partial_t \gamma$:
\begin{gather*}
    \begin{aligned}
        &\bm{J}_\gamma \cdot \bm{n}_\partial = i_\gamma\\
        &\bm{J}_\gamma \cdot \bm{n}_\partial = 0
    \end{aligned}
    \quad \text{on } \partial_t \gamma,
\end{gather*}
for current and potential electrodes, respectively, with $i_\gamma = \pi r^2
\overline{J}_\Gamma$ being the current in \sib{\ampere}.

In $\Upsilon$, the electrode is seen as an immersed line acting as
source/sink of electric current. Introducing  the Dirac delta $\delta_\gamma$ distributed along $\gamma$ with measure
\sib{\per\square\metre}, the new set of equations for $\Upsilon$ reads
\begin{gather}\label{eq:cont_electr}
    \begin{aligned}
        &\bm{J} + \sigma \nabla \varphi = \bm{0}\\
        &\nabla \cdot \bm{J} + j_\gamma \delta_\gamma =0
    \end{aligned}
    \quad \text{in } \Upsilon
\end{gather}
Solutions of the previous equations exhibit a logarithmic
singularity in $\varphi$ with strong derivatives in the vicinity of $\gamma$.
Assuming a radial current density exiting from the electrode, following the
approach given \cite{Peaceman1978,Peaceman1983}, it is possible to write the coefficient in
\eqref{eq:varsigma_electrde} as
\begin{gather*}
    \sigma_\gamma = \frac{2\pi \sigma}{\ln(r_e/r) + S},
\end{gather*}
which is a semi-discrete parameter that is used to mimic the effect of the
singularity around $\gamma$, when a discretization in space is applied, without
the need of over-refining the grid.
The parameter $r_e$ in \sib{\metre} is the radius at which the electric potential in the medium
$\varphi$ is equal to the averaged
grid cell potential.
$r_e$ is normally taken equal to $r_e\approx 0.2 h$ with $h$, in
\sib{\metre}, the cell diameter.
$S$ in
\sib{\cdot} is the so-called skin-factor, which gives an additional electric
potential drop due to the actual features of the electrode and of the
surrounding zone.
$S$ can be used to model rust or other imperfections of the electrode, and also
to model the imperfect contact between the electrode and the ground.
 We are finally ready to introduce the new mixed-dimensional problem.
\begin{problem}[The mixed-dimensional electrode DC problem]\label{pb:electrode}
    Find $(\bm{J}, \varphi)$ in $\Upsilon$ and $(j_\gamma, \bm{J}_\gamma, \varphi_\gamma)$
    in $\gamma$ such that
    \begin{gather*}
        \begin{aligned}
            &\bm{J} + \sigma \nabla \varphi = \bm{0}\\
            &\nabla \cdot \bm{J} + j_\gamma \delta_\gamma =  0
        \end{aligned}
        \quad \text{in } \Upsilon,
        \qquad
        \begin{aligned}
            &\bm{J}_\gamma + \pi r^2 \sigma_\Gamma \nabla \varphi_\gamma = \bm{0}\\
            &\nabla \cdot \bm{J}_\gamma - j_\gamma = 0
        \end{aligned}
        \quad \text{in } \gamma,
    \end{gather*}
    coupled with the following interface condition
    \begin{gather*}
        j_\gamma + \sigma_\gamma (\varphi_\gamma - \varphi) = 0,
        \quad \text{on } \gamma
    \end{gather*}
    and with boundary conditions given by
    \begin{gather*}
        \begin{aligned}
            &\bm{J} \cdot \bm{n}_\partial = 0 && \text{on } \partial \Upsilon,\\
            &\bm{J}_\gamma \cdot \bm{n}_\partial = \pi r^2 \overline{J}_\Gamma
            && \text{on } \partial_t \gamma,\\
            &\bm{J}_\gamma \cdot \bm{v}_\Gamma = 0 &&  \text{on } \partial_b
            \gamma.
        \end{aligned}
    \end{gather*}
\end{problem}

\subsection{A model for the liner}\label{subsec:liner}

A liner is placed below the waste body
to avoid the infiltration of leachate in the underground, that
can be
dangerous for the groundwater.
For simplicity, we assume that the liner consists of a planar portion of a highly
resistive geomembrane whose thickness is much smaller than other lengths
in the domain. For real applications, the ratio between its thickness and lateral extension
is in the order of $\sim 10^{-5}$. To simplify gridding and avoid
excessive refinement in representing the liner as a three-dimensional object we
present a two-dimensional model of the liner with suitable coupling conditions
with the surrounding media.

In the case of a liner composed by several planar parts,
e.g. having a box shape, we will apply the new model to each part and use suitable
conditions to connect them. This aspect will be detailed later.

\begin{figure}[tb]
    \centering
    \resizebox{0.25\textwidth}{!}{\fontsize{0.75cm}{2cm}\selectfont
\begingroup%
  \makeatletter%
  \providecommand\color[2][]{%
    \errmessage{(Inkscape) Color is used for the text in Inkscape, but the package 'color.sty' is not loaded}%
    \renewcommand\color[2][]{}%
  }%
  \providecommand\transparent[1]{%
    \errmessage{(Inkscape) Transparency is used (non-zero) for the text in Inkscape, but the package 'transparent.sty' is not loaded}%
    \renewcommand\transparent[1]{}%
  }%
  \providecommand\rotatebox[2]{#2}%
  \newcommand*\fsize{\dimexpr\f@size pt\relax}%
  \newcommand*\lineheight[1]{\fontsize{\fsize}{#1\fsize}\selectfont}%
  \ifx\svgwidth\undefined%
    \setlength{\unitlength}{272.250359bp}%
    \ifx\svgscale\undefined%
      \relax%
    \else%
      \setlength{\unitlength}{\unitlength * \real{\svgscale}}%
    \fi%
  \else%
    \setlength{\unitlength}{\svgwidth}%
  \fi%
  \global\let\svgwidth\undefined%
  \global\let\svgscale\undefined%
  \makeatother%
  \begin{picture}(1,0.83471096)%
    \lineheight{1}%
    \setlength\tabcolsep{0pt}%
    \put(0.03935454,0.56473752){\color[rgb]{0,0,0}\makebox(0,0)[lt]{\lineheight{1.25}\smash{\begin{tabular}[t]{l}$\Upsilon$\end{tabular}}}}%
    \put(0,0){\includegraphics[width=\unitlength,page=1]{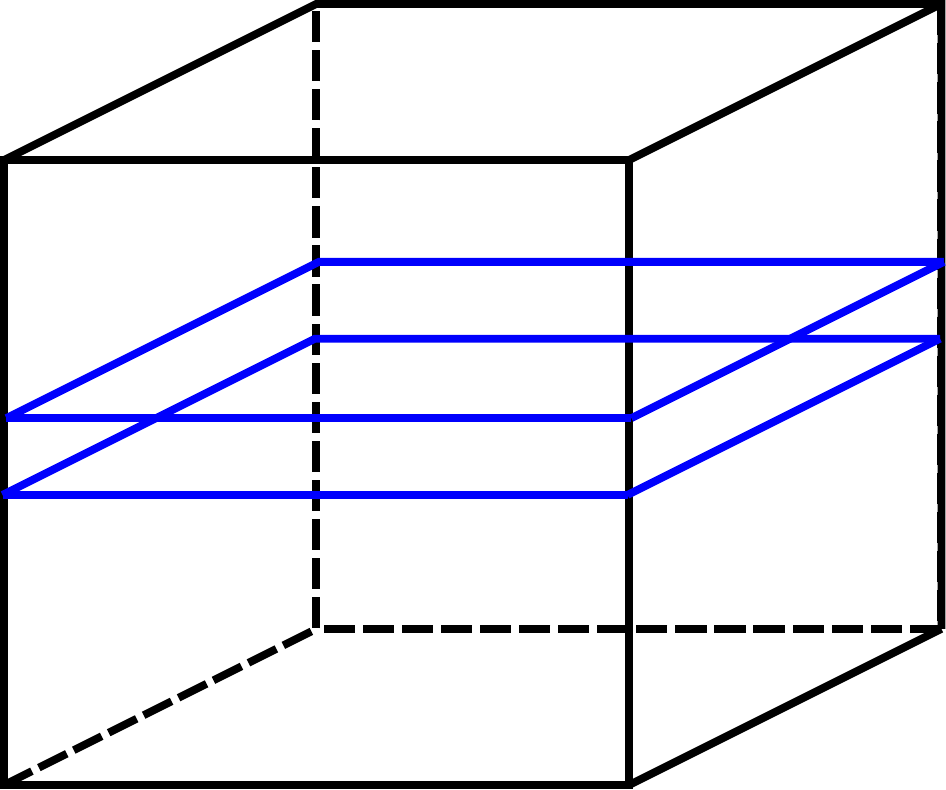}}%
    \put(0.8154503,0.31783799){\color[rgb]{0,0,1}\makebox(0,0)[lt]{\lineheight{1.25}\smash{\begin{tabular}[t]{l}$\Lambda$\end{tabular}}}}%
    \put(0,0){\includegraphics[width=\unitlength,page=2]{liner_full.pdf}}%
    \put(1.0284,0.50510625){\color[rgb]{0,0,0}\makebox(0,0)[lt]{\lineheight{1.25}\smash{\begin{tabular}[t]{l}$\varepsilon$\end{tabular}}}}%
  \end{picture}%
\endgroup%
}%
    \hspace{0.1\textwidth}%
    \resizebox{0.25\textwidth}{!}{\fontsize{0.75cm}{2cm}\selectfont
\begingroup%
  \makeatletter%
  \providecommand\color[2][]{%
    \errmessage{(Inkscape) Color is used for the text in Inkscape, but the package 'color.sty' is not loaded}%
    \renewcommand\color[2][]{}%
  }%
  \providecommand\transparent[1]{%
    \errmessage{(Inkscape) Transparency is used (non-zero) for the text in Inkscape, but the package 'transparent.sty' is not loaded}%
    \renewcommand\transparent[1]{}%
  }%
  \providecommand\rotatebox[2]{#2}%
  \newcommand*\fsize{\dimexpr\f@size pt\relax}%
  \newcommand*\lineheight[1]{\fontsize{\fsize}{#1\fsize}\selectfont}%
  \ifx\svgwidth\undefined%
    \setlength{\unitlength}{272.250359bp}%
    \ifx\svgscale\undefined%
      \relax%
    \else%
      \setlength{\unitlength}{\unitlength * \real{\svgscale}}%
    \fi%
  \else%
    \setlength{\unitlength}{\svgwidth}%
  \fi%
  \global\let\svgwidth\undefined%
  \global\let\svgscale\undefined%
  \makeatother%
  \begin{picture}(1,0.83471096)%
    \lineheight{1}%
    \setlength\tabcolsep{0pt}%
    \put(0.03935454,0.56473752){\color[rgb]{0,0,0}\makebox(0,0)[lt]{\lineheight{1.25}\smash{\begin{tabular}[t]{l}$\Upsilon$\end{tabular}}}}%
    \put(0,0){\includegraphics[width=\unitlength,page=1]{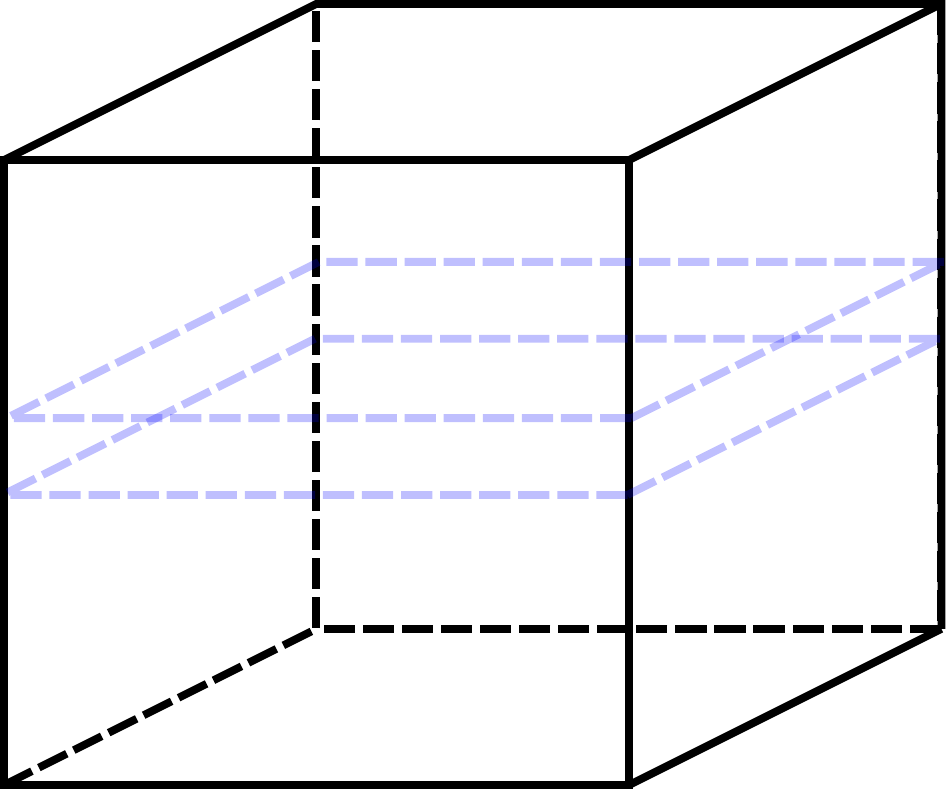}}%
    \put(0.8154503,0.31783799){\color[rgb]{0,0,1}\makebox(0,0)[lt]{\lineheight{1.25}\smash{\begin{tabular}[t]{l}$\lambda$\end{tabular}}}}%
    \put(0,0){\includegraphics[width=\unitlength,page=2]{liner_reduced.pdf}}%
  \end{picture}%
\endgroup%
}
    \caption{The liner $\Lambda$ represented as three-dimensional
    object (right) and the reduced liner $\lambda$ as two dimensional
    object (left).}
    \label{fig:liner}
\end{figure}

Let us consider Figure \ref{fig:liner} on the left where a
liner $\Lambda$  is immersed in $\Upsilon$.
It can be described as
a domain with
thickness $\varepsilon$ in $\sib{\metre}$ and central (portion of a) plane $\lambda$,
\begin{gather*}
    \Lambda = \lambda \times \left(-\frac{\varepsilon}{2},
    \frac{\varepsilon}{2}\right).
\end{gather*}
We indicate with
$\bm{n}_\Lambda$ the unit normal to the top boundary of $\Lambda$ pointing from
$\Upsilon$ toward $\Lambda$. Equations of
Problem \ref{eq:the_system}
apply also to $\Lambda$ with different conductivity than the surrounding $\Upsilon$,
therefore, as done for the electrode, we propose now a suitable reduced model
where the liner is approximated with its central plane $\lambda$, by
following the same approach as in
\cite{Martin2005,DAngelo2011,Nordbotten2018,Fumagalli2019a}.
First, we divide the boundary $\partial \Lambda$ of
$\Lambda$ into two parts: the lateral surface of thickness $\varepsilon$ is
named $\partial_l \Lambda$, and the bottom and top surfaces are identified as
$\partial_s \Lambda$.
If the liner ends inside $\Upsilon$, we can identify the portion of its
lateral boundary that is in contact with the internal part of $\Upsilon$.
We call this internal part of
the boundary as $\partial_i \Lambda \subset \partial_l \Lambda$.
Let us state the following assumptions.
\begin{assumption}[Liner thickness]
    We assume that the thickness $\varepsilon$ of the liner is much smaller than any
    other sizes in the domain.
\end{assumption}
\begin{assumption}[Lateral exchange]
    We assume that current density is exchanged between the liner
    and the surrounding media through $\partial_s \Lambda$, the
    exchange though $\partial_i \Lambda$ is assumed to be negligible.
\end{assumption}
We substitute the equi-dimensional representation of the liner $\Lambda$ with
its lower-dimensional counterpart $\lambda$,
so that the mesh size is not constraned by its thickness
$\varepsilon$. We present now the mathematical model for $\lambda$. Note that
the liner is in contact with $\Upsilon$ with its top
and bottom surfaces.

For each side of $\lambda$, the current density exchanged between the liner and the surrounding media, following
\eqref{eq:varsigma_electrde}, is denoted as
$j_\lambda$ in \sib{\ampere\per\square\metre}
and approximated by
\begin{gather}\label{eq:varsigma}
    \varepsilon j_\lambda + \sigma_\lambda (\varphi_\lambda - \varphi) = 0
    \quad \text{on } \lambda.
\end{gather}
Being $\lambda$ of codimension 1, it does not introduce
singular solutions and thus we do not have the need to tune $\sigma_\lambda$
as we did for $\sigma_\gamma$.
Indeed, the  value of $\sigma_\lambda$
in \sib{\siemens\per\metre} simply
represents
the conductivity of the liner along its normal.
The reduced variables are $\varphi_\lambda$ in \sib{\volt}
and $\bm{J}_\lambda$ now in
\sib{\ampere\per\metre} and defined as
\begin{gather*}
    \varphi_\lambda(x, y) = \frac{1}{\varepsilon}\int_{-\frac{\varepsilon}{2}}^{\frac{\varepsilon}{2}}
    \varphi_\Lambda(x, y, z) d z
    \quad \text{and} \quad
    \bm{J}_\gamma(x, y) = \int_{-\frac{\varepsilon}{2}}^{\frac{\varepsilon}{2}}
    (I - \bm{n}_\Lambda \otimes \bm{n}_\Lambda)
    \bm{J}_\Lambda (x, y, z) d z,
\end{gather*}
with $I$ the identity matrix.
We can now introduce a model for the liner
represented as an object $\lambda$ of dimension $n-1$
\begin{gather*}
    \begin{aligned}
        &\bm{J}_\lambda + \varepsilon \sigma_\Lambda \nabla \varphi_\lambda = \bm{0}\\
        &\nabla \cdot \bm{J}_\lambda - j_\lambda = 0
    \end{aligned}
    \quad \text{in } \lambda.
\end{gather*}
The gradient and divergence are now computed over the plane $\lambda$.
To couple this model with Problem \ref{eq:the_system} in $\Upsilon$, we simply
require current conservation so that $\bm{J}\cdot \bm{n}_\Lambda =
j_\lambda$ on each side of $\lambda$.
In the experiments we will consider the boundary of $\lambda$ to be
in contact with $\partial \Upsilon$, and thus inheriting the same boundary conditions, or
immersed into $\Upsilon$ in which case we impose null current exhange, as in \eqref{eq:tip_condition}.

When the liner geometry is more complex but can still be split into multiple
planar polygons, we can follow the same strategy for each one of them. Since the liner is composed by a
homogeneous material, at the interface between two planes we
simply require that the current density is conserved and that the electric
potential is continuous.

The problem is thus formalized as follows.
\begin{problem}[The mixed-dimensional liner DC problem]\label{eq:liner}
    Find $(\bm{J}, \varphi)$ in $\Upsilon$ and $(j_\lambda, \bm{J}_\lambda,
    \varphi_\lambda)$
    in $\lambda$ such that
    \begin{gather*}
        \begin{aligned}
            &\bm{J} + \sigma \nabla \varphi = \bm{0}\\
            &\nabla \cdot \bm{J}=  0
        \end{aligned}
        \quad \text{in } \Upsilon
        \qquad
        \begin{aligned}
            &\bm{J}_\lambda + \varepsilon \sigma_\Lambda \nabla \varphi_\lambda = \bm{0}\\
            &\nabla \cdot \bm{J}_\lambda - j_\lambda = 0
        \end{aligned}
        \quad \text{in } \lambda
    \end{gather*}
    coupled with the following interface condition on both sides of $\lambda$
    \begin{gather*}
        \varepsilon j_\lambda + \sigma_\lambda ( \varphi_\lambda - \varphi) = 0
        \quad \text{on } \lambda
    \end{gather*}
    with boundary conditions given by
    \begin{gather*}
        \begin{aligned}
            &\bm{J} \cdot \bm{n}_\partial = 0 && \text{on } \partial \Upsilon\\
            &\bm{J}_\lambda \cdot \bm{n}_\partial = 0 &&  \text{on } \partial
            \lambda
        \end{aligned}
    \end{gather*}
\end{problem}

\subsection{The complete mixed-dimensional model}

\begin{figure}[tb]
    \centering
    \resizebox{0.25\textwidth}{!}{\fontsize{0.75cm}{2cm}\selectfont
\begingroup%
  \makeatletter%
  \providecommand\color[2][]{%
    \errmessage{(Inkscape) Color is used for the text in Inkscape, but the package 'color.sty' is not loaded}%
    \renewcommand\color[2][]{}%
  }%
  \providecommand\transparent[1]{%
    \errmessage{(Inkscape) Transparency is used (non-zero) for the text in Inkscape, but the package 'transparent.sty' is not loaded}%
    \renewcommand\transparent[1]{}%
  }%
  \providecommand\rotatebox[2]{#2}%
  \newcommand*\fsize{\dimexpr\f@size pt\relax}%
  \newcommand*\lineheight[1]{\fontsize{\fsize}{#1\fsize}\selectfont}%
  \ifx\svgwidth\undefined%
    \setlength{\unitlength}{272.250359bp}%
    \ifx\svgscale\undefined%
      \relax%
    \else%
      \setlength{\unitlength}{\unitlength * \real{\svgscale}}%
    \fi%
  \else%
    \setlength{\unitlength}{\svgwidth}%
  \fi%
  \global\let\svgwidth\undefined%
  \global\let\svgscale\undefined%
  \makeatother%
  \begin{picture}(1,0.83471096)%
    \lineheight{1}%
    \setlength\tabcolsep{0pt}%
    \put(0,0){\includegraphics[width=\unitlength,page=1]{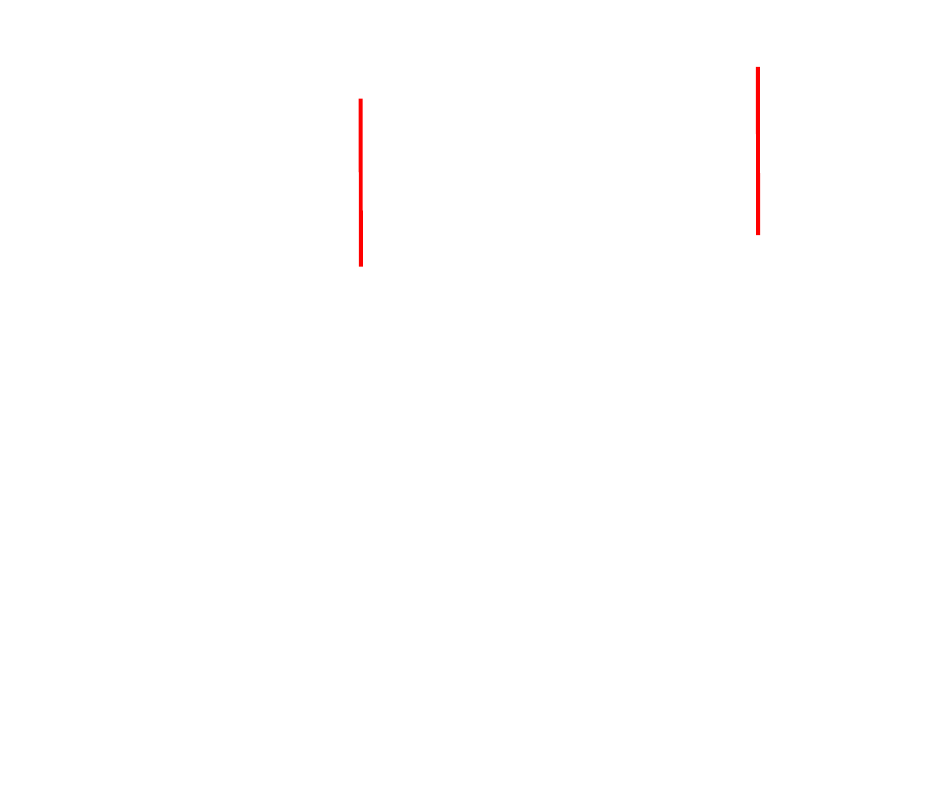}}%
    \put(0.03935454,0.56473752){\color[rgb]{0,0,0}\makebox(0,0)[lt]{\lineheight{1.25}\smash{\begin{tabular}[t]{l}$\Upsilon$\end{tabular}}}}%
    \put(0,0){\includegraphics[width=\unitlength,page=2]{all_reduced.pdf}}%
    \put(0.8154503,0.36742403){\color[rgb]{0,0,1}\makebox(0,0)[lt]{\lineheight{1.25}\smash{\begin{tabular}[t]{l}$\lambda$\end{tabular}}}}%
    \put(0,0){\includegraphics[width=\unitlength,page=3]{all_reduced.pdf}}%
    \put(0.82183692,0.60743109){\color[rgb]{1,0,0}\makebox(0,0)[lt]{\lineheight{1.25}\smash{\begin{tabular}[t]{l}$\gamma_1$\end{tabular}}}}%
    \put(0.39877571,0.57201224){\color[rgb]{1,0,0}\makebox(0,0)[lt]{\lineheight{1.25}\smash{\begin{tabular}[t]{l}$\gamma_2$\end{tabular}}}}%
  \end{picture}%
\endgroup%
}%
    \caption{Schematic representation of the mixed-dimensional objects: the
    three-dimensional domain $\Upsilon$, the two-dimensional liner $\lambda$,
    and the one-dimensional electrodes $\gamma_i$.}
    \label{fig:all}
\end{figure}

For simplicity we assume that the electrodes and the liner are not intersecting
nor in contact, so
the global problem considers both contributions presented in Problem \ref{pb:electrode} and Problem
\ref{eq:liner} in a straightforward case. The mixed-dimensional model consider a set of $N$ electrodes,
indexed with $i=1, \ldots, N$, and a single liner, see Figure \ref{fig:all}.
The equations are presented
as following.
\begin{problem}[The mixed-dimensional DC problem]\label{eq:model}
    Find $(\bm{J}, \varphi)$ in $\Upsilon$,
    $(j_{\gamma_i}, \bm{J}_{\gamma_i}, \varphi_{\gamma_i})$
    in $\gamma_i$ for $i=1, \ldots, N$
    and
    $(j_\lambda, \bm{J}_\lambda,
    \varphi_\lambda)$
    in $\lambda$ such that
    \begin{gather*}
        \begin{aligned}
            &\bm{J} + \sigma \nabla \varphi = \bm{0}\\
            &\nabla \cdot \bm{J} + {\textstyle\sum_{i}} j_{\gamma_i} \delta_{\gamma_i} =  0
        \end{aligned}
        \quad \text{in } \Upsilon
        \qquad
        \begin{aligned}
            &\bm{J}_{\gamma_i} + \pi r_i^2 \sigma_{\Gamma_i} \nabla
            \varphi_{\gamma_i} = \bm{0}\\
            &\nabla \cdot \bm{J}_{\gamma_i} - j_{\gamma_i} = 0
        \end{aligned}
        \quad \text{in } \gamma_i
        \qquad
        \begin{aligned}
            &\bm{J}_\lambda + \varepsilon \sigma_\Lambda \nabla \varphi_\lambda = \bm{0}\\
            &\nabla \cdot \bm{J}_\lambda - j_\lambda = 0
        \end{aligned}
        \quad \text{in } \lambda
    \end{gather*}
    coupled with the following interface condition
    \begin{gather*}
        j_{\gamma_i} + \sigma_{\gamma_i} (\varphi_{\gamma_i} - \varphi) = 0
        \quad \text{on } \gamma_i
        \qquad
        \varepsilon j_\lambda + \sigma_\lambda ( \varphi_\lambda - \varphi) = 0
        \quad \text{on } \lambda
    \end{gather*}
    with boundary conditions given by
    \begin{gather*}
        \begin{aligned}
            &\bm{J} \cdot \bm{n}_\partial = 0 && \text{on } \partial \Upsilon\\
            &\bm{J}_\lambda \cdot \bm{n}_\partial = 0 &&  \text{on } \partial
            \lambda
        \end{aligned}
        \qquad
        \begin{aligned}
            &\bm{J}_{\gamma_i} \cdot \bm{n}_\partial = \pi r_i^2
            \overline{J}_{\Gamma_i}
            && \text{on } \partial_t \gamma_i\\
            &\bm{J}_{\gamma_i} \cdot \bm{v}_{\Gamma_i} = 0 &&  \text{on } \partial_b
            \gamma_i\\
        \end{aligned}
    \end{gather*}
\end{problem}
By balancing the current density imposed on the boundary,
Problem \ref{eq:model} admits a unique solution for $(\bm{J}, \bm{J}_{\gamma_i},
j_{\gamma_i}, \bm{J}_\lambda,  j_\lambda)$, however the electric potentials
$(\varphi, \varphi_{\gamma_i}, \varphi_\lambda)$
are defined up to a constant $c \in \mathbb{R}$. To uniquely define this
constant
it is possible to consider different strategies like imposing the following
null average condition
\begin{gather*}
    \int_\Upsilon \varphi + \sum_{i=1}^{N} \int_{\gamma_i} \varphi_{\gamma_i} +
    \int_\lambda \varphi_\lambda = 0.
\end{gather*}
However, in our application, we are interested in potential differences, thus the
constant $c$ does not have any practical influence.

\section{Numerical approximation}\label{sec:approximation}

In this section we discuss the numerical approximation of Problem
\ref{eq:model},
detailing the strategies adopted to represent the mixed-dimensional objects in
an efficient and effective way.

We start by considering each object separately and, in Subsection
\ref{subsec:coupling}, we discuss their coupling. By substituting the
constitutive equations of Problem \ref{eq:model} into their respective
conservation equation, we obtain a problem in primal form in terms only of the electric
potentials and $j_\lambda$ and $j_{\gamma_i}$. Being $\Upsilon$, $\lambda$ and $\gamma_i$ objects of dimension 3,
2 and 1, respectively, we construct a simplicial grid composed of tetrahedra for
$\Upsilon$, triangles for $\lambda$ and segments for each electrode $\gamma_i$.
For their construction we rely on Gmsh \cite{Geuzaine2009}; more details
are given in Subsection \ref{subsec:grid}.

For the numerical approximation of the
problem, we consider two finite-volume cell-centered schemes normally used in
the context of flows in porous media and other elliptic problems: the
two-point flux approximation (TPFA) as presented in \cite{Aavatsmark2007a}
and the multi-point flux approximation (MPFA) as introduced in
\cite{Aavatsmark2007,Aavatsmark2002}. Both schemes have one degree of freedom
for each cell but different stencils. TPFA is consistent only for the so called
$k$-orthogonal grids,
which are grids
where the faces normals and the connecting lines between the centers of neighboring
cells are orthogonal with respect to the inverse of the diffusion coefficient.
Moreover, TPFA
is very efficient in the construction and solution of the associated linear
problem. MPFA is more expensive but it is convergent for simplicial grids, with
only weak regularity requirements on the mesh \cite{Aavatsmark2007,Aavatsmark2002}.
We will test in Section \ref{sec:experiments} the performance of both
methods with respect to laboratory experiments.

\subsection{Coupling between dimensions}\label{subsec:coupling}

A key aspect in the numerical approximation of Problem \ref{eq:model} is the
actual implementation of the coupling conditions between the electrodes, or the
liner, and $\Upsilon$. We consider the approach  presented in
\cite{Nordbotten2018,Keilegavlen2020}, where new grids, called mortar grids, are constructed to handle
the coupling between dimensions. In particular, since $j_\lambda$ is defined on
both sides of the liner and couples it with the surrounding domain, we construct
two mortar grids that put in communication $\lambda$ with $\Upsilon$. In
principle, these grids may be non-conforming and transfer operators
between them need to be constructed that are not, in general, identity maps. To keep our presentation simple, we
consider only conforming approximations, see \cite{Nordbotten2018,Keilegavlen2020} for more
details for the non-conforming case.

Since the electrodes are mono-dimensional, and smaller than any other objects present
in the problem, to increase the flexibility of the proposed approach, we have
considered for them a non-matching discretization meaning that each segment
composing $\gamma_i$ is immersed in some of the tetrahedra of the $\Upsilon$ grid.
A single mortar grid, for each
$\gamma_i$, is
constructed for the coupling. In this case, transfer operators between
$\Upsilon$, $\gamma_i$ and the mortar grid are accordingly constructed to map
variables and then discretize the coupling condition \eqref{eq:varsigma_electrde}. In this case, $j_{\gamma_i}$ is
the variable defined on the mortar grid that connects $\gamma_i$ with $\Upsilon$.

\begin{figure}[tb]
    \centering
    \includegraphics[width=0.75\textwidth]{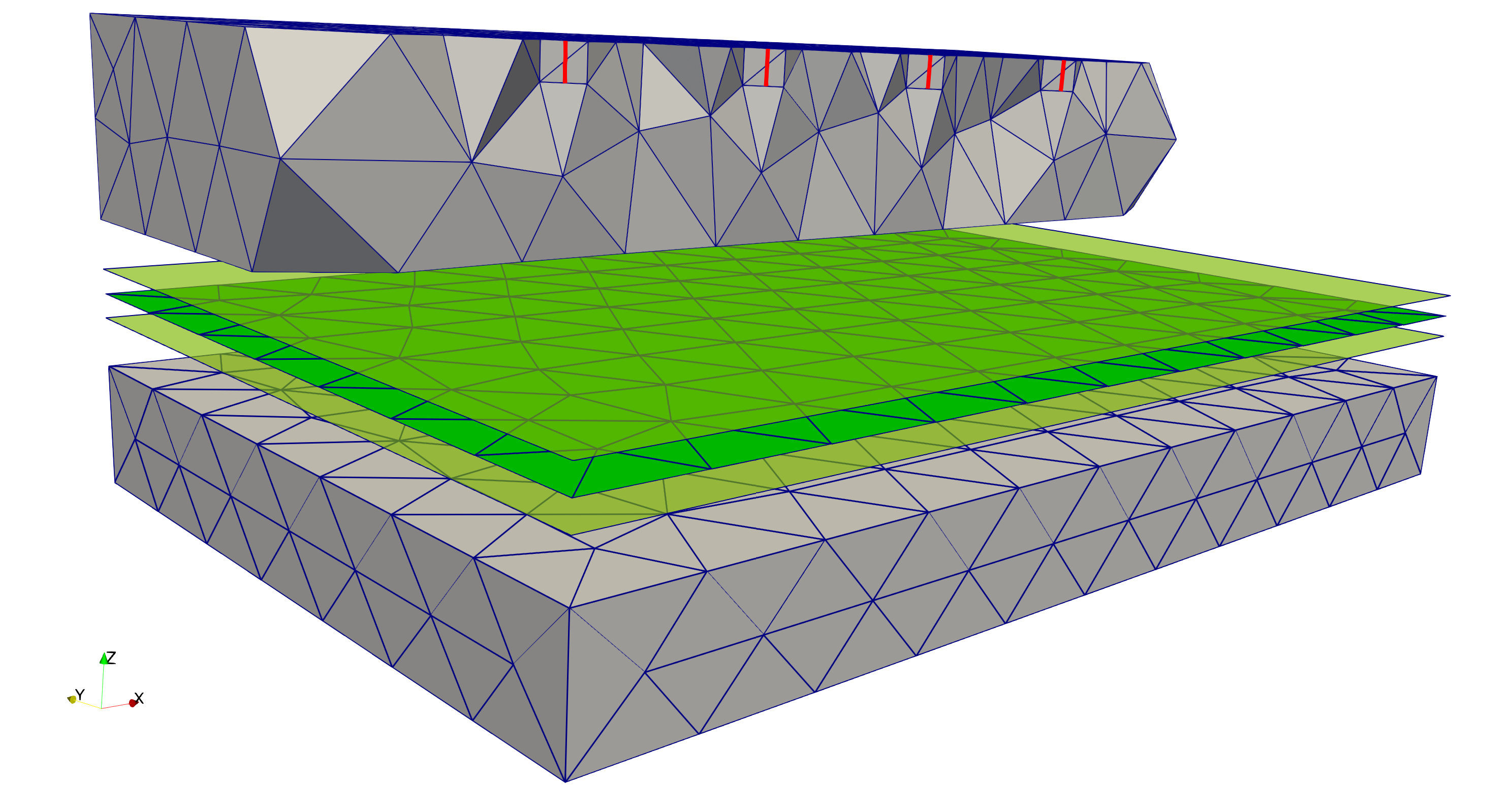}%
    \caption{Example of computational grid constructed for the approximation of
    the problem. The different objects are separated to make the visualization clearer.}
    \label{fig:grid_coupling}
\end{figure}

An example of discretized problem is reported in Figure \ref{fig:grid_coupling}, where the liner is
approximated with the green grid and the two mortar grids surrounding it are in
transparent light green. At the  top of the domain, in red, a set of four electrodes
is represented as non-matching with the grid of $\Upsilon$.

\subsection{Construction of the computational grid}\label{subsec:grid}

\begin{figure}[tb]
    \centering
    \resizebox{0.55\textwidth}{!}{\fontsize{1cm}{2cm}\selectfont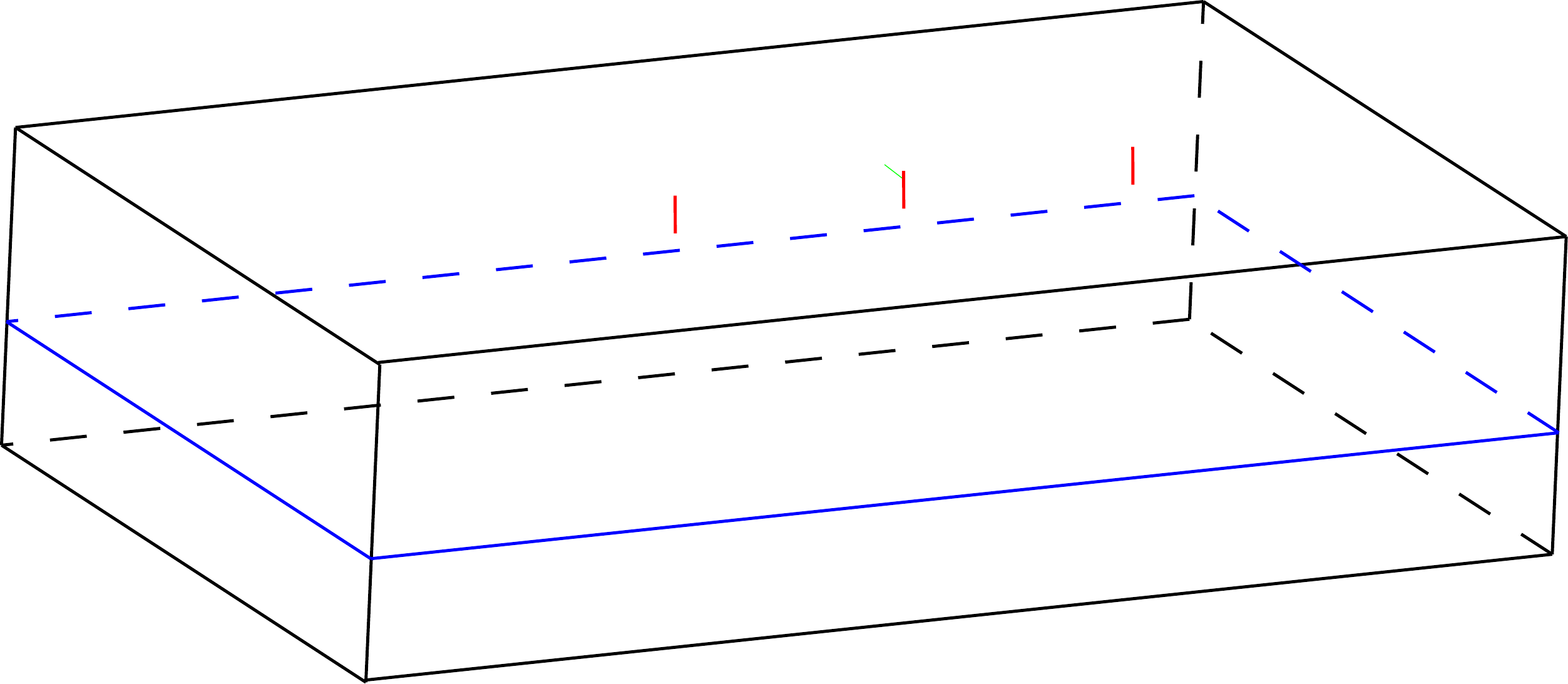}%
    \hspace*{0.125\textwidth}%
    \raisebox{0.035\textwidth}{\includegraphics[width=0.15\textwidth]{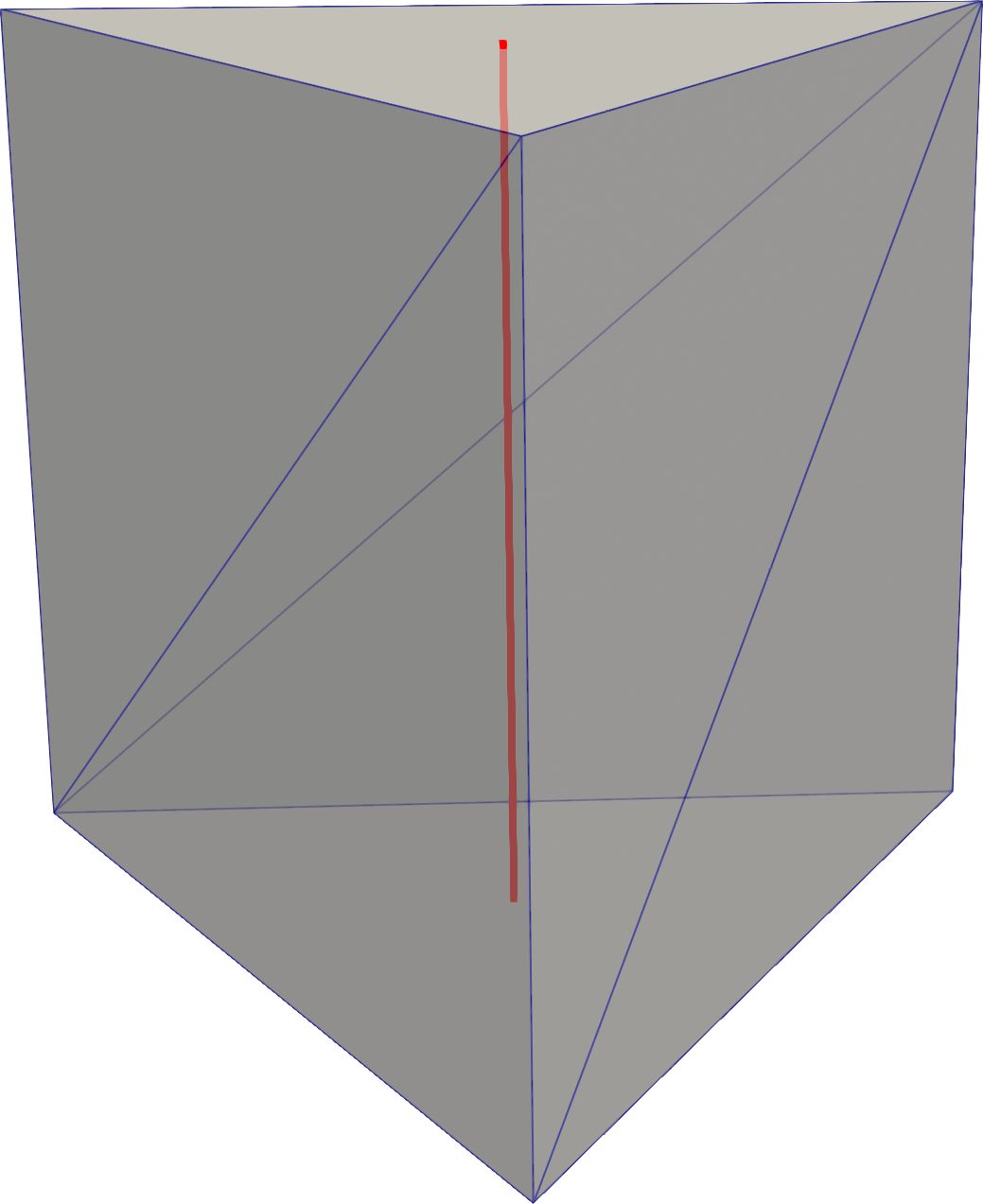}}%
    \hspace*{0.025\textwidth}%
    \raisebox{0.035\textwidth}{\includegraphics[width=0.15\textwidth]{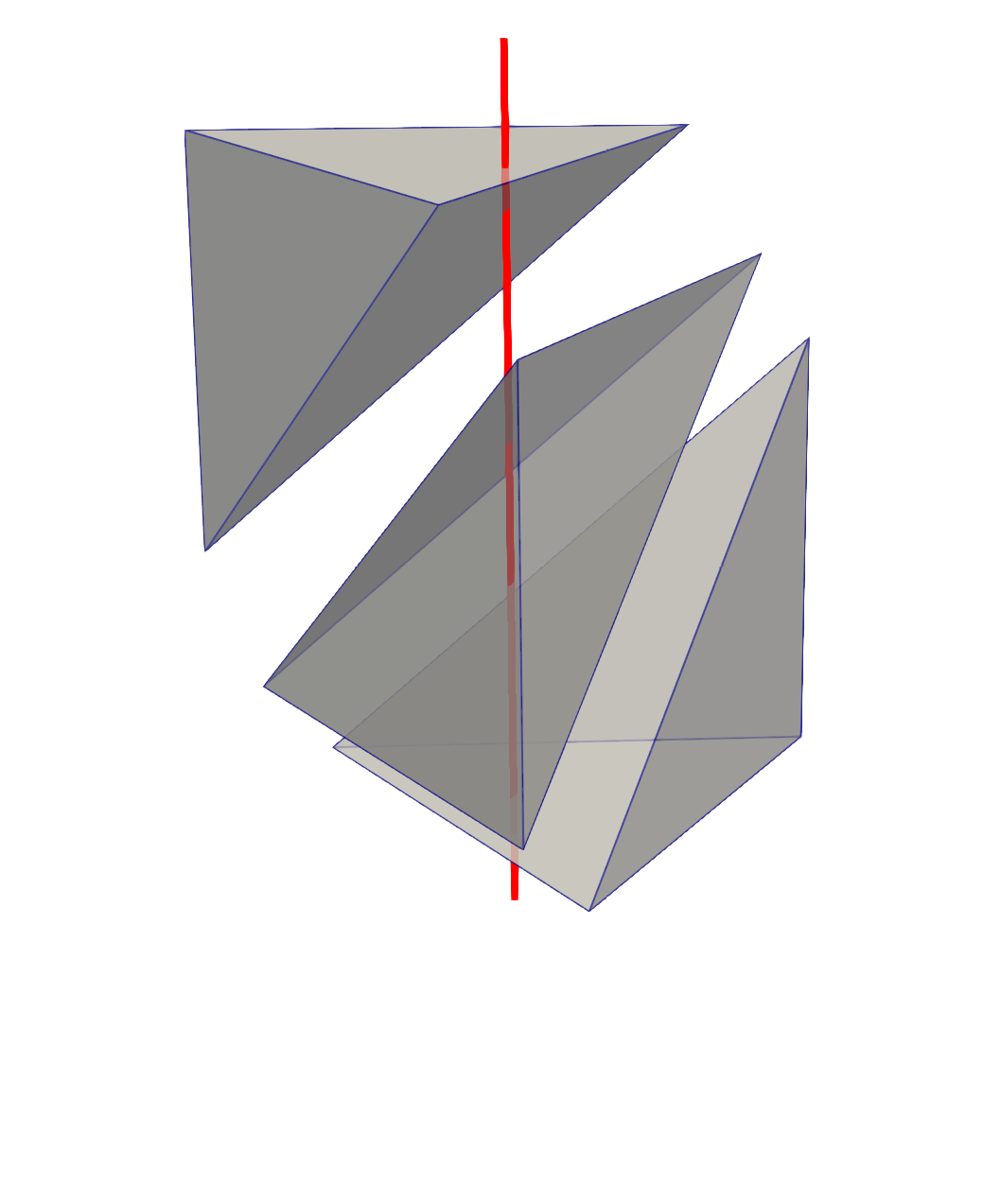}}%
    \caption{Graphical representation of the domain $\Upsilon$, the liner
    $\lambda$ in blue and the electrodes $\gamma_i$ in red. The different groups
    of green tetrahedra
    are associated to each electrodes. On the right, a detail on the construction
    of cells around an electrode.}
    \label{fig:geometry}
\end{figure}

Another important task in the discretization is to guide the gridding tool to create
appropriate grids, \cite{Aavatsmark2007,Aavatsmark2002}, and thus limit the
numerical error introduced. As discussed before,
the grid of the liner is conforming with the surrounding domain, however the
electrodes are assumed to be immersed in the cells of $\Upsilon$. The presented approach is
more effective when the electrodes are inserted in the centre of the cells, far
from their edges.
To match this condition, our approach is to construct a-priori a set of cells surrounding the electrodes
and force the gridding tool to include them in the construction of the grid. The
error introduced by the non-conforming approximation is thus minimised and the
procedure is fully automatized. Figure \ref{fig:geometry} shows the geometry
given to the gridding  tool to construct the full three-dimensional grid. On the
top of the domain, each of the four electrodes
is immersed in a wedge discretized with three tetrahedra.

Finally, to detect which cell of $\Upsilon$ intersects each electrode segment,
we have considered a fast algorithm based on an alternating digital tree
(ADT), see \cite{Bonet1991,Thompson1998} for more details. Once these cells are identified, the maps
between them and the intersected electrodes segments are constructed. In general,
these are not identity maps.

\section{Numerical experiments}\label{sec:experiments}

The purpose of this section is to validate the proposed approach against laboratory experiments and to study the distribution of current density and electric potential in presence of a highly-resistive liner.
We fill a $52 \sib{\centi\metre} \times34\sib{\centi\metre}\times40\sib{\centi\metre}$
plastic box with different volumes of tap water, whose
conductivity is estimated with a Crison MM40+ multimeter to be equal to $\sigma = 1/29$
\sib{\siemens\per\meter}.

\begin{figure}
    \centering
    \raisebox{0.05\textwidth}{\resizebox{0.45\textwidth}{!}{\fontsize{0.75cm}{2cm}\selectfont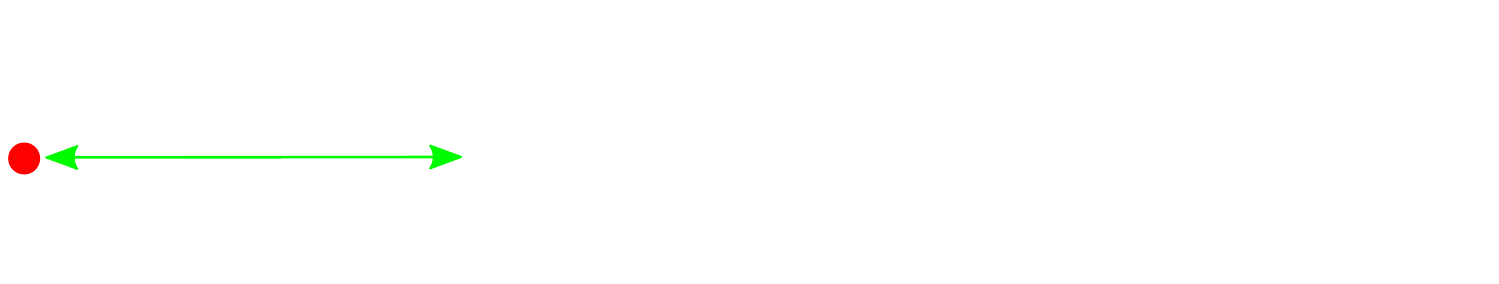}}%
    \resizebox{0.55\textwidth}{!}{\fontsize{2cm}{2cm}\selectfont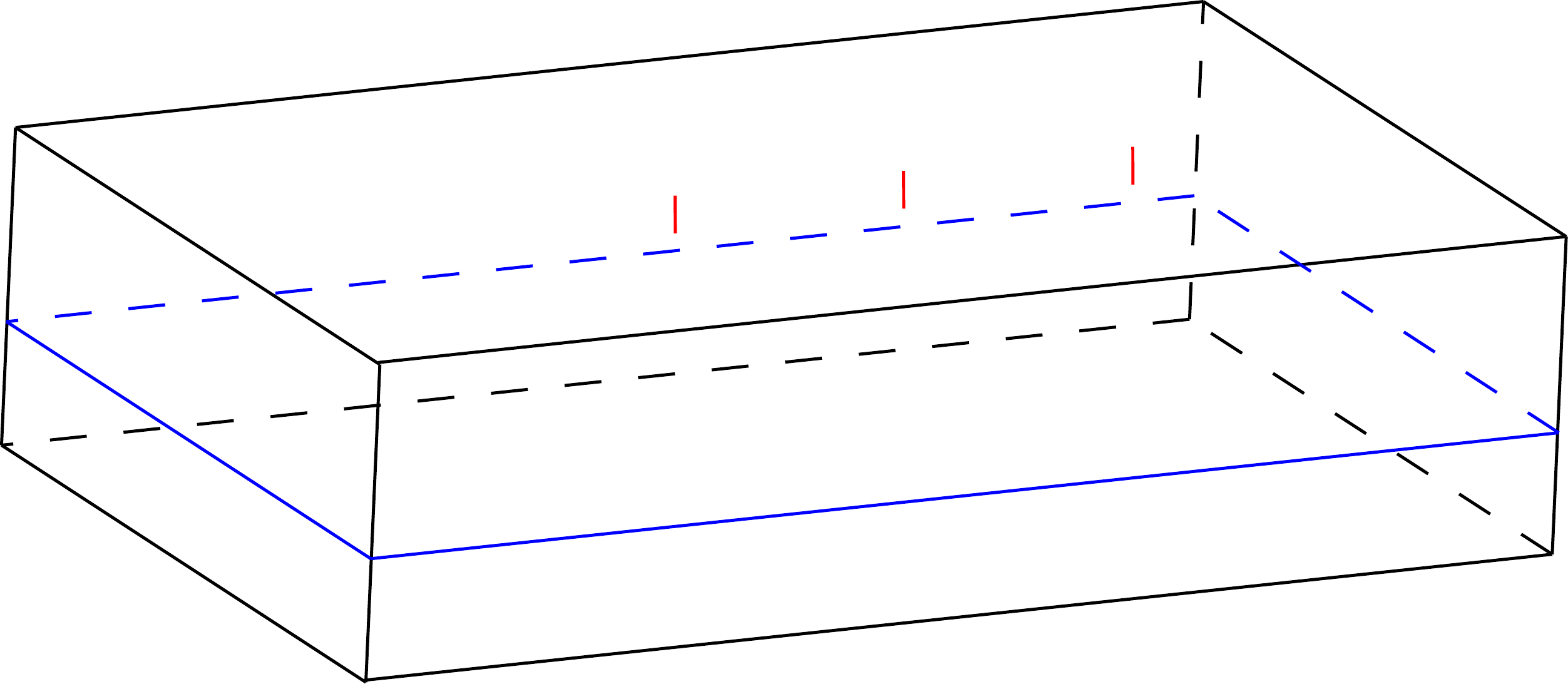}%
    \caption{Wenner-$\alpha$ configuration with four electrodes with spacing $a$.}
    \label{fig:wenner_config}
\end{figure}

We deploy a spread of 4 electrodes according to a
Wenner-$\alpha$ configuration at the centre of the box along the $x$-axis (Figure \ref{fig:wenner_config}).
The two outer electrodes $C_1$ and $C_2$ are used to inject the current
density
with value $\pm \overline{J}_\Gamma$, while the two
inner ones $P_1$ and $P_2$ are used to compute the potential difference
$\varphi_{\gamma_2} - \varphi_{\gamma_3}$, the latter being evaluated at the top of the
electrodes.
The penetration depth of the Wenner-$\alpha$ array for a homogeneous medium can be analytically estimated to be about $0.11L$ or $0.17L$, being $L$ the distance between the current electrodes, either one considers the peak or the median value of the one-dimensional sensitivity function, respectively \cite{roy1971depth, BARKER1989}. However, the penetration depth of a geoelectrical survey only estimates the depth at which the used array is maximally sensitive and may vary considerably if the investigated medium is heterogeneous.

We set $\Upsilon$ to be completely isolated with null current flow through all the boundaries (i.e., homogeneous Neumann boundary condition) as presented in Problem \ref{eq:model}. The domain is
discretized with an unstructured mesh with a characteristic length of 5
\sib{\centi\meter}, which is successively refined to 1 \sib{\centi\meter} near the electrodes. No significant improvements are observed in the numerical results considering smaller characteristic lengths.

The
electrodes have radius equal to $r=1\sib{\milli\meter}$, length of $5
\sib{\milli\meter}$ and conductivity equal to $\sigma_\gamma = 1.45 \cdot 10^6
\sib{\siemens\per\meter}$. The ratio between the radius and the length of the electrodes is larger than what stated previously for convenience in the lab experiments. However, numerical tests with smaller ratios down to 1/100 do not show significant differences. The liner, which will have different shapes for each
example, has thickness equal to $\varepsilon = 1
\sib{\milli\meter}$ and conductivity $\sigma_\lambda =
10^{-9} \sib{\siemens\per\meter}$. As mentioned before, the liner is approximated as a
bi-dimensional object. Contrary to the electrodes, the computational grid is
conforming to the liner.

We consider two sets of experiments with increasing level of complexity: in Subsection \ref{subsec:itact_liner} the
liner is a horizontal plane, while in Subsection \ref{subsec:broken_liner} the
liner has a box shape. The latter mimics the shape of the liner that is placed underneath real landfills to completely isolate the waste body from the surrounding media. In such a case, it is interesting to analyse current and electrical potential distribution according to the position of the electrodes and to the presence of damages in the liner. We will compare the numerical and the experimental values
of the estimated
apparent resistivity $\rho_a$ in \sib{\ohm\metre}, defined for the
Wenner-$\alpha$ configuration as
\begin{gather*}
    \rho_{a} = 2\pi a \frac{\varphi_{\gamma_2} -
    \varphi_{\gamma_3}}{i_\gamma},
\end{gather*}
where the value of $i_\gamma$ is imposed and the difference
$\varphi_{\gamma_2} - \varphi_{\gamma_3}$ comes from measurements or modelling.

\subsection{A horizontal liner}\label{subsec:itact_liner}

\begin{figure}
    \centering
    \includegraphics[scale=0.45]{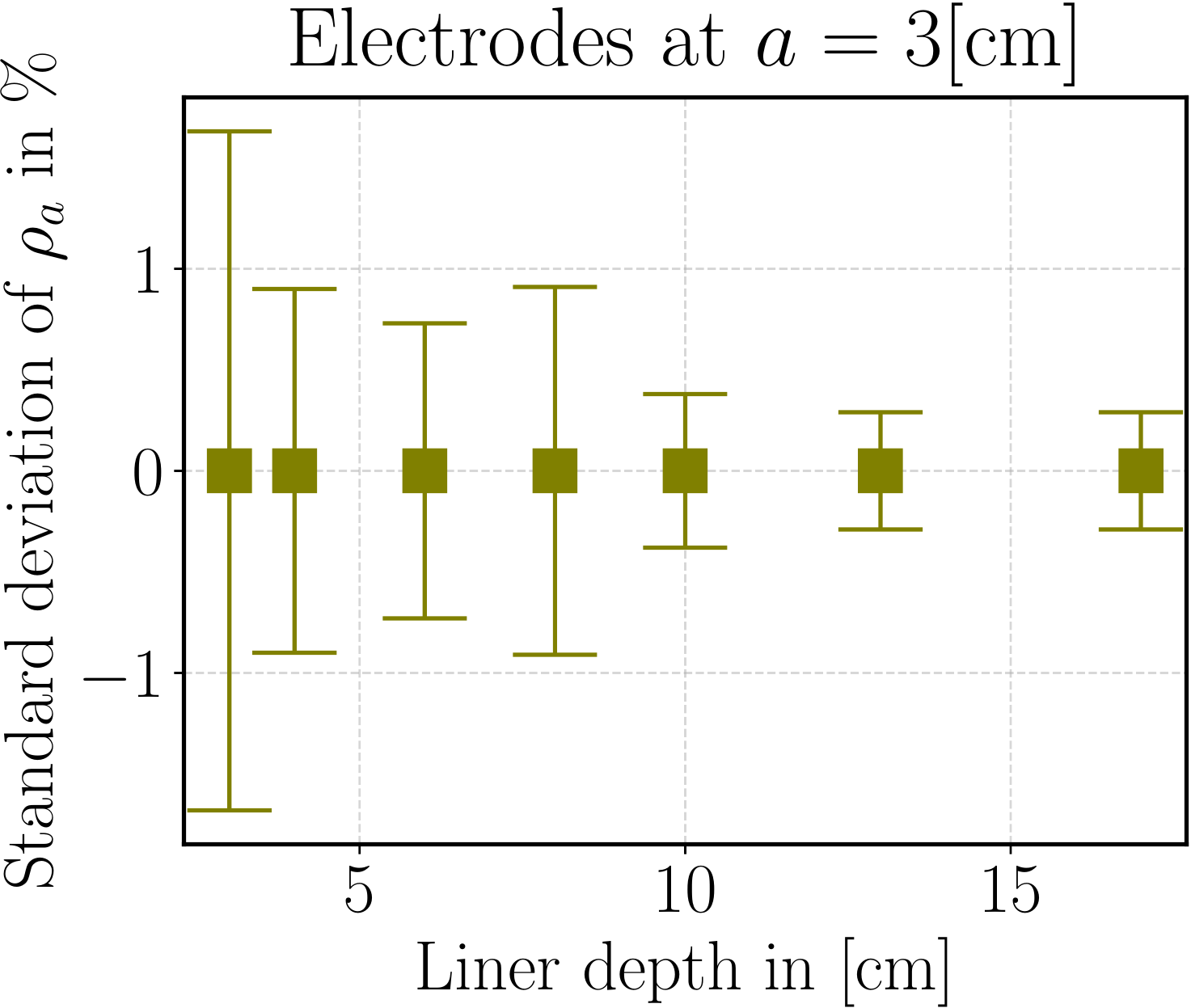}%
    \hspace{0.05\textwidth}%
    \includegraphics[scale=0.45]{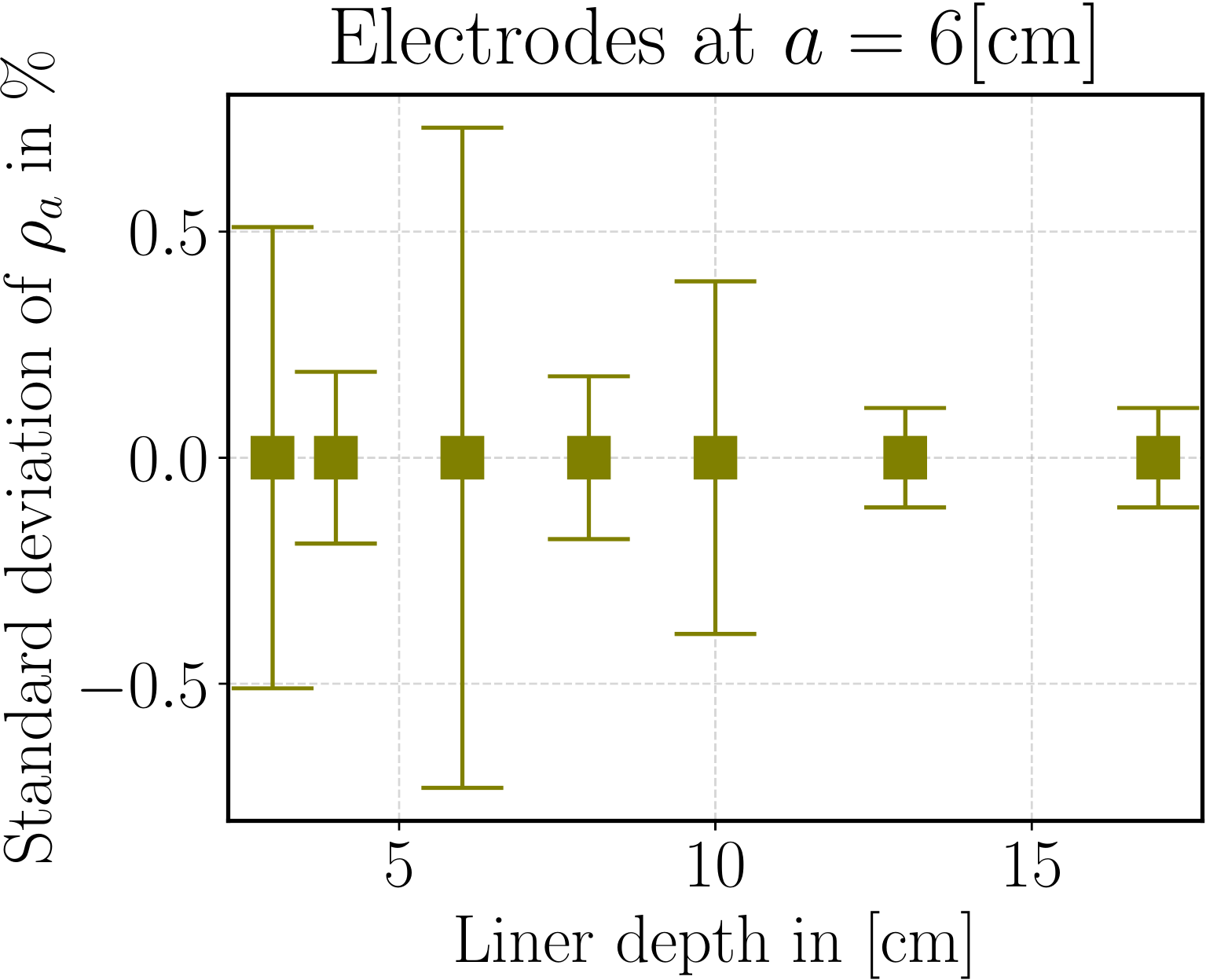}
    \caption{Percentages of standard deviations, centred at their mean value, for the measurements of $\rho_a$
    for the two experiments, for the example in Subsection \ref{subsec:itact_liner}.}
    \label{fig:experiment1result_b_plot}
\end{figure}

In this first test, we want to evaluate the effect of the liner on $\rho_a$.
The liner is perfectly horizontal and touches all the lateral boundaries of $\Upsilon$. We further assume that $\lambda$ is intact, meaning that neither holes nor
defects are present. We consider two settings with different spacing
between the electrodes: $a=3\sib{\centi\meter}$ and
$a=6\sib{\centi\meter}$. Since the liner is a highly resistive
membrane, it can be approximated with a perfect insulator. In the laboratory test, instead of introducing the liner, we
simply decrease the height of the water layer from $h=17
\sib{\centi\metre}$ to $h=3
\sib{\centi\metre}$ to ease operations. Following \cite{Sheriff1990}, as reference value we
consider the apparent resistivity
$\widetilde{\rho}_a$ in \sib{\ohm\meter} analytically expressed as
\begin{gather}\label{eq:rho_teo}
    \widetilde{\rho}_a = \frac{\rho i_\gamma}{2 \pi}
    \left(
    \frac{1}{a}
    + 4 \sum_{n=1}^\infty
    \frac{1}{\sqrt{a^2 + 4n^2 h^2}}
    -
    \frac{1}{\sqrt{4a^2 + 4n^2 h^2}}
    \right),
\end{gather}
with $\rho = 29 \sib{\ohm\meter}$ being the resistivity of the water.
This formula was derived under the condition that $\Upsilon$ is indefinitely
laterally extended (i.e., 1D model) and the electrodes are represented as zero-dimensional
sources, and, in this work, further simplified since the resistivity of the liner is several orders of magnitude higher than the water resistivity.
The instrument used in the laboratory performs several measurements of $\rho_a$ for each depth of the liner, so that we can estimate the mean value and the standard deviation, which can be
used as proxy for the reliability of the measurements, as reported in Figure
\ref{fig:experiment1result_b_plot}. In addition, measurement errors were also checked with reciprocal measurements, i.e. switching current and potential electrodes, but values are rather small confirming that the measured $\rho_a$ values are reliable.

\begin{figure}
    \centering
    \resizebox{0.35\textwidth}{!}{\fontsize{0.5cm}{2cm}\selectfont
\begingroup%
  \makeatletter%
  \providecommand\color[2][]{%
    \errmessage{(Inkscape) Color is used for the text in Inkscape, but the package 'color.sty' is not loaded}%
    \renewcommand\color[2][]{}%
  }%
  \providecommand\transparent[1]{%
    \errmessage{(Inkscape) Transparency is used (non-zero) for the text in Inkscape, but the package 'transparent.sty' is not loaded}%
    \renewcommand\transparent[1]{}%
  }%
  \providecommand\rotatebox[2]{#2}%
  \newcommand*\fsize{\dimexpr\f@size pt\relax}%
  \newcommand*\lineheight[1]{\fontsize{\fsize}{#1\fsize}\selectfont}%
  \ifx\svgwidth\undefined%
    \setlength{\unitlength}{167.42576899bp}%
    \ifx\svgscale\undefined%
      \relax%
    \else%
      \setlength{\unitlength}{\unitlength * \real{\svgscale}}%
    \fi%
  \else%
    \setlength{\unitlength}{\svgwidth}%
  \fi%
  \global\let\svgwidth\undefined%
  \global\let\svgscale\undefined%
  \makeatother%
  \begin{picture}(1,0.46398296)%
    \lineheight{1}%
    \setlength\tabcolsep{0pt}%
    \put(0,0){\includegraphics[width=\unitlength,page=1]{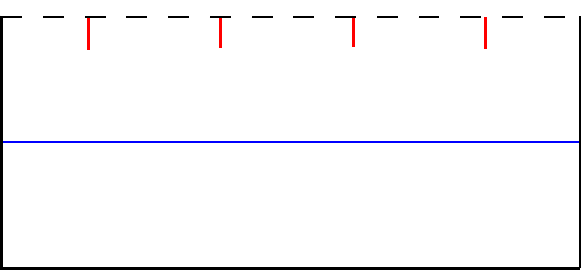}}%
    \put(0.4730842,0.25862307){\color[rgb]{0,0,1}\makebox(0,0)[lt]{\lineheight{1.25}\smash{\begin{tabular}[t]{l}$\lambda$\end{tabular}}}}%
    \put(0.16656363,0.37285249){\color[rgb]{1,0,0}\makebox(0,0)[lt]{\lineheight{1.25}\smash{\begin{tabular}[t]{l}$\gamma_1$\end{tabular}}}}%
    \put(0.39146793,0.37094644){\color[rgb]{1,0,0}\makebox(0,0)[lt]{\lineheight{1.25}\smash{\begin{tabular}[t]{l}$\gamma_2$\end{tabular}}}}%
    \put(0.62083284,0.37094644){\color[rgb]{1,0,0}\makebox(0,0)[lt]{\lineheight{1.25}\smash{\begin{tabular}[t]{l}$\gamma_3$\end{tabular}}}}%
    \put(0.84878351,0.37094644){\color[rgb]{1,0,0}\makebox(0,0)[lt]{\lineheight{1.25}\smash{\begin{tabular}[t]{l}$\gamma_4$\end{tabular}}}}%
  \end{picture}%
\endgroup%
}%
    \hspace*{0.1\textwidth}%
    \resizebox{0.35\textwidth}{!}{\fontsize{0.5cm}{2cm}\selectfont
\begingroup%
  \makeatletter%
  \providecommand\color[2][]{%
    \errmessage{(Inkscape) Color is used for the text in Inkscape, but the package 'color.sty' is not loaded}%
    \renewcommand\color[2][]{}%
  }%
  \providecommand\transparent[1]{%
    \errmessage{(Inkscape) Transparency is used (non-zero) for the text in Inkscape, but the package 'transparent.sty' is not loaded}%
    \renewcommand\transparent[1]{}%
  }%
  \providecommand\rotatebox[2]{#2}%
  \newcommand*\fsize{\dimexpr\f@size pt\relax}%
  \newcommand*\lineheight[1]{\fontsize{\fsize}{#1\fsize}\selectfont}%
  \ifx\svgwidth\undefined%
    \setlength{\unitlength}{167.42576899bp}%
    \ifx\svgscale\undefined%
      \relax%
    \else%
      \setlength{\unitlength}{\unitlength * \real{\svgscale}}%
    \fi%
  \else%
    \setlength{\unitlength}{\svgwidth}%
  \fi%
  \global\let\svgwidth\undefined%
  \global\let\svgscale\undefined%
  \makeatother%
  \begin{picture}(1,0.46398296)%
    \lineheight{1}%
    \setlength\tabcolsep{0pt}%
    \put(0,0){\includegraphics[width=\unitlength,page=1]{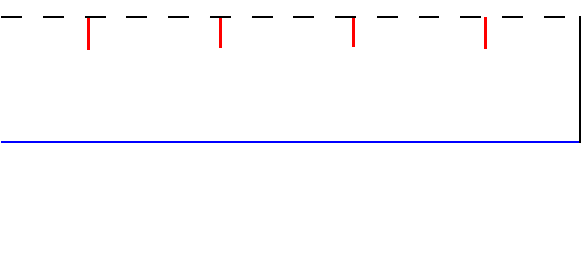}}%
    \put(0.4730842,0.25862307){\color[rgb]{0,0,1}\makebox(0,0)[lt]{\lineheight{1.25}\smash{\begin{tabular}[t]{l}$\lambda$\end{tabular}}}}%
    \put(0,0){\includegraphics[width=\unitlength,page=2]{example1_without.pdf}}%
    \put(0.16660004,0.37284024){\color[rgb]{1,0,0}\makebox(0,0)[lt]{\lineheight{1.25}\smash{\begin{tabular}[t]{l}$\gamma_1$\end{tabular}}}}%
    \put(0.39150435,0.37093419){\color[rgb]{1,0,0}\makebox(0,0)[lt]{\lineheight{1.25}\smash{\begin{tabular}[t]{l}$\gamma_2$\end{tabular}}}}%
    \put(0.62086929,0.37093419){\color[rgb]{1,0,0}\makebox(0,0)[lt]{\lineheight{1.25}\smash{\begin{tabular}[t]{l}$\gamma_3$\end{tabular}}}}%
    \put(0.8488199,0.37093419){\color[rgb]{1,0,0}\makebox(0,0)[lt]{\lineheight{1.25}\smash{\begin{tabular}[t]{l}$\gamma_4$\end{tabular}}}}%
  \end{picture}%
\endgroup%
}%
    \caption{Schematic representation of the domain of example in Subsection \ref{subsec:itact_liner}. The resistive interface of the liner is represented by an interface (left) or by the bottom of the domain (right).}
    \label{fig:domain_example1}
\end{figure}

In the numerical simulations, we
adopt two strategies to account for the liner $\lambda$: i) we consider the entire domain $\Upsilon$ and place a planar surface with very low electric
conductivity at a certain depth; ii) since $\lambda$ is highly
resistive, we simply remove the part of $\Upsilon$ below the liner and impose a
null current density flux condition at the boundary, as schematically shown in Figure \ref{fig:domain_example1}.
Both TPFA and MPFA are considered as discretization schemes.
We can thus evaluate the effect of the liner and compare
the obtained results with the laboratory experiments.

\begin{figure}
    \centering
    \includegraphics[scale=0.45]{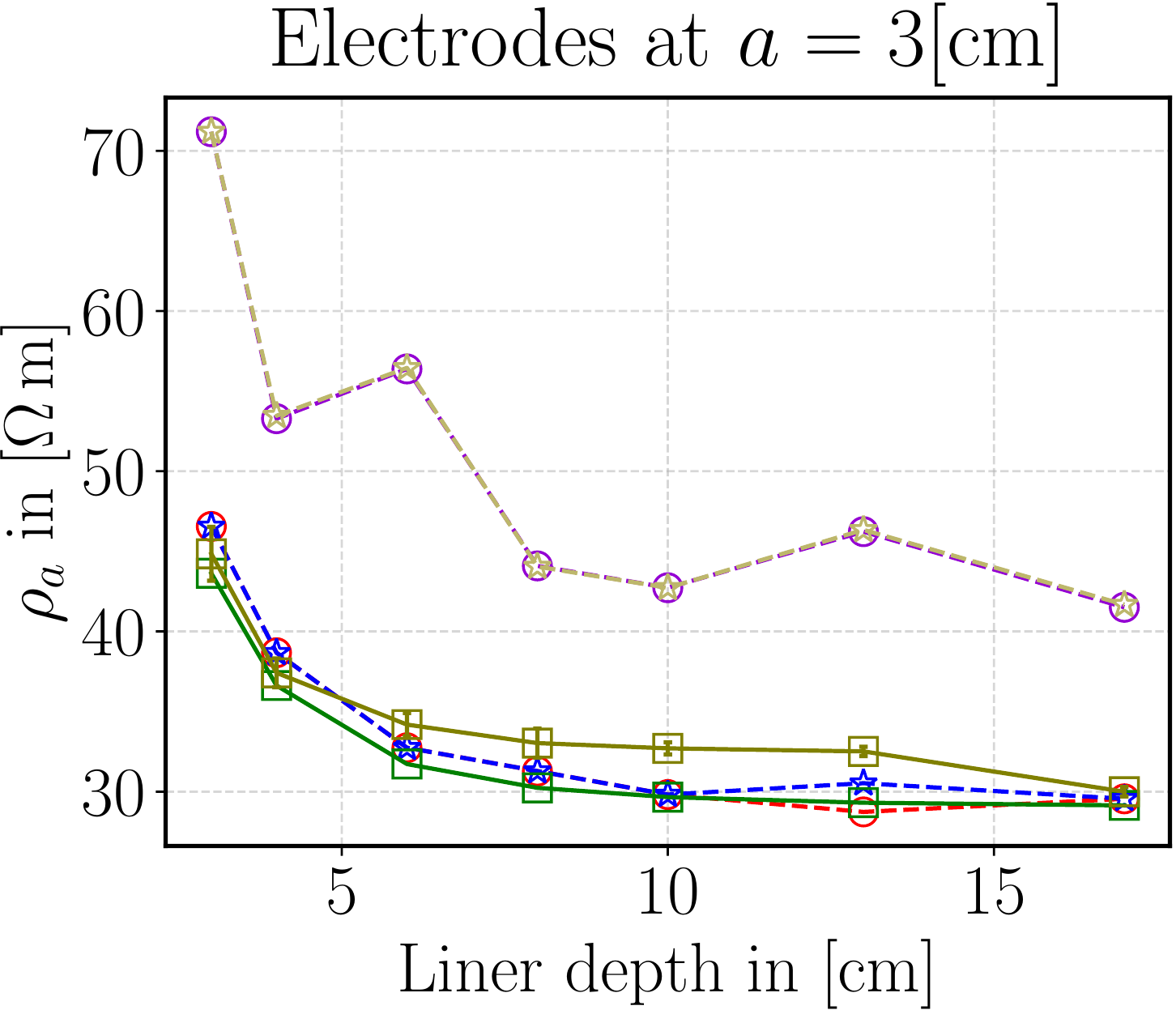}%
    \hspace{0.05\textwidth}%
    \includegraphics[scale=0.45]{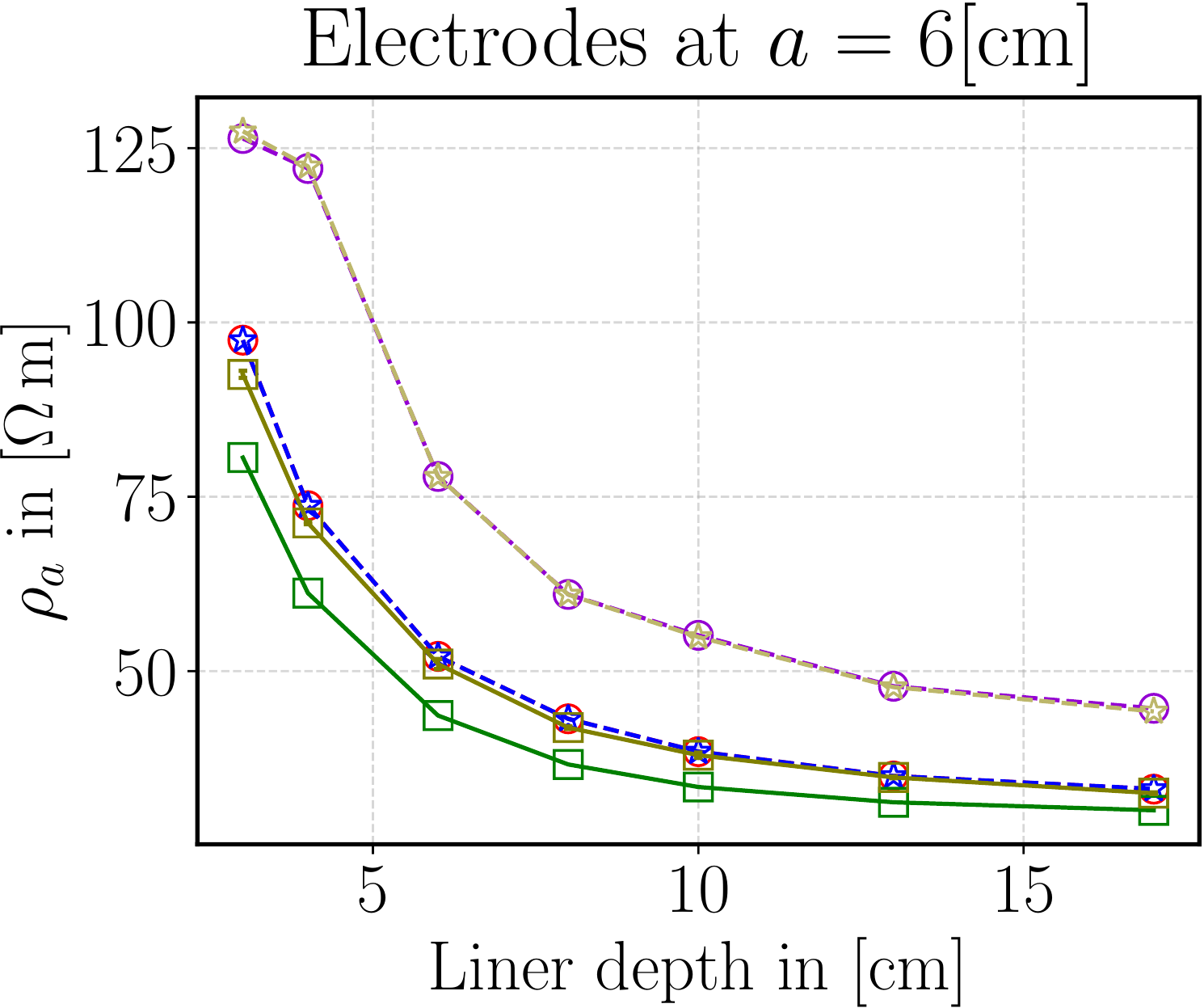}\\[0.25cm]
    \includegraphics[width=0.85\textwidth]{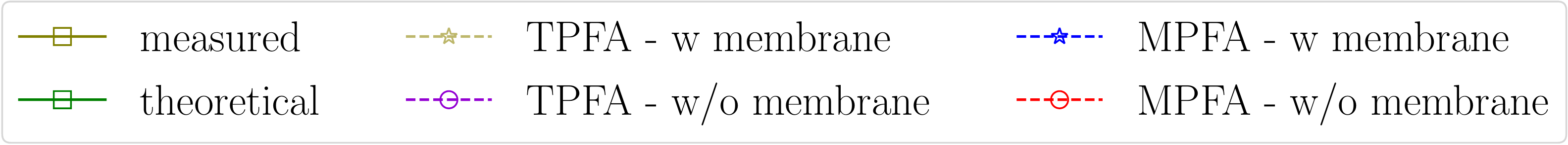}
    \caption{Results for the example in Subsection \ref{subsec:itact_liner}. Apparent resistivities computed and measured for several liner depths and for electrode spacing $a=3\sib{\centi\meter}$ (left) and $a=6\sib{\centi\meter}$ (right).}
    \label{fig:experiment1result}
\end{figure}

In Figure \ref{fig:experiment1result} we show the apparent resistivities
obtained in the simulations, in the laboratory experiments, and with the
analytical relationship \eqref{eq:rho_teo}. The measured standard deviations (Figure \ref{fig:experiment1result_b_plot}) are also represented and now centred
in the actual values of the apparent resistivity. Since these standard
deviations are rather small they are hardly visible in the graph.
All the results show the same trend:
by decreasing the liner depth, the apparent resistivity increases because the liner is closer to the electrodes.

We note that the results of the simulations does not change if either $\lambda$ is included in the modelled domain or $\Upsilon$ is reduced accordingly. This means that the
laboratory setting is appropriate and can be compared with the simulations when
$\lambda$ is considered. Moreover, there is a clear difference between the
results obtained from the TPFA and MPFA schemes, having for the former, for
$a=3\sib{\centi\meter}$, a non-monotone behaviour with respect to the liner
depth.

By comparing the synthetic data with the measured ones, we observe that the numerical solutions computed with MPFA are in good agreement,
while the TPFA ones are less accurate. A possible explanation is that the latter
is not consistent when using simplicial grids and, therefore, the obtained results might not be reliable.

\begin{figure}
    \centering
    \includegraphics[scale=0.45]{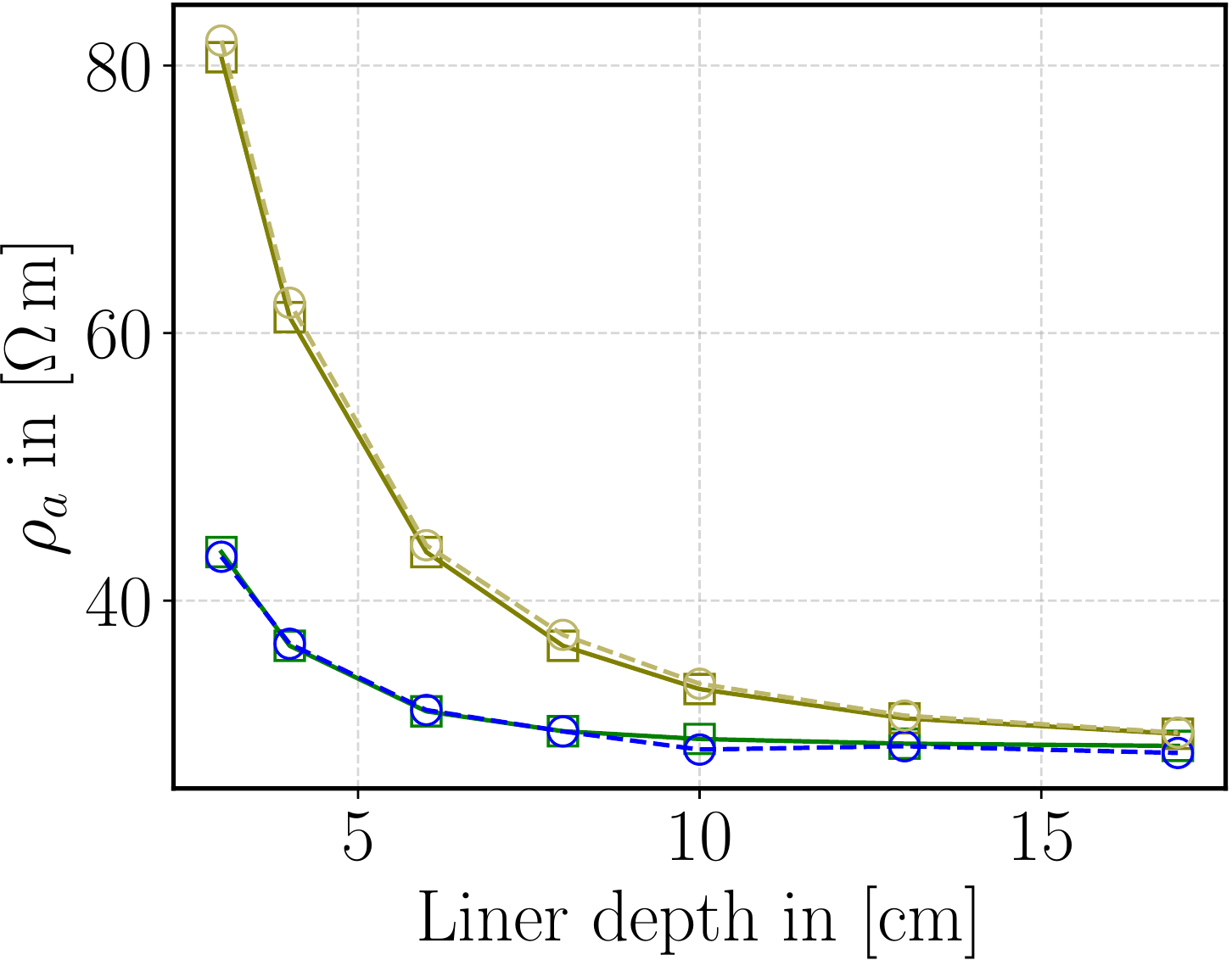}\\[0.25cm]
    \includegraphics[width=0.6\textwidth]{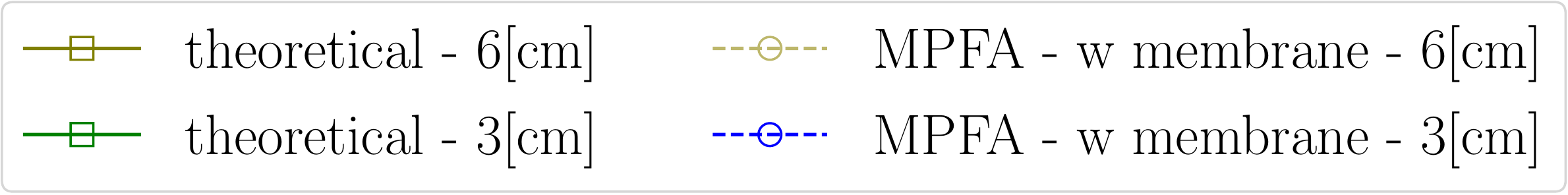}
    \caption{Apparent resistivities computed with MPFA and estimated with relationship \eqref{eq:rho_teo} for different values of $a$ when the computational domain $\Upsilon$ is laterally extended to $100\times100\times40$ \sib{\centi\metre\cubed}.}
    \label{fig:experiment1resulta}
\end{figure}

Finally, we observe that, for a spacing of $a=6\sib{\centi\meter}$, the measured and computed values of $\rho_a$ are higher than the actual one. This can be due to the fact that the considered model is not 1D, thus we could not apply relationship \eqref{eq:rho_teo}.
Figure \ref{fig:experiment1resulta} shows a comparison between the numerical
results and the values obtained with  \eqref{eq:rho_teo} when the computational domain
$\Upsilon$ is laterally extended to $100\times100\times40$
\sib{\centi\metre\cubed}. The results are now in good agreement, confirming that
the discrepancy in Figure \ref{fig:experiment1result} for
$a=6\sib{\centi\meter}$ is mainly due to boundary effects. With respect to the values of penetration depth stated before, our tests show that the Wenner-$\alpha$ array becomes insensitive to the resistive liner at depths approximately greater than 3 times the electrode spacing.

\begin{figure}
    \centering
    \includegraphics[width=0.475\textwidth]{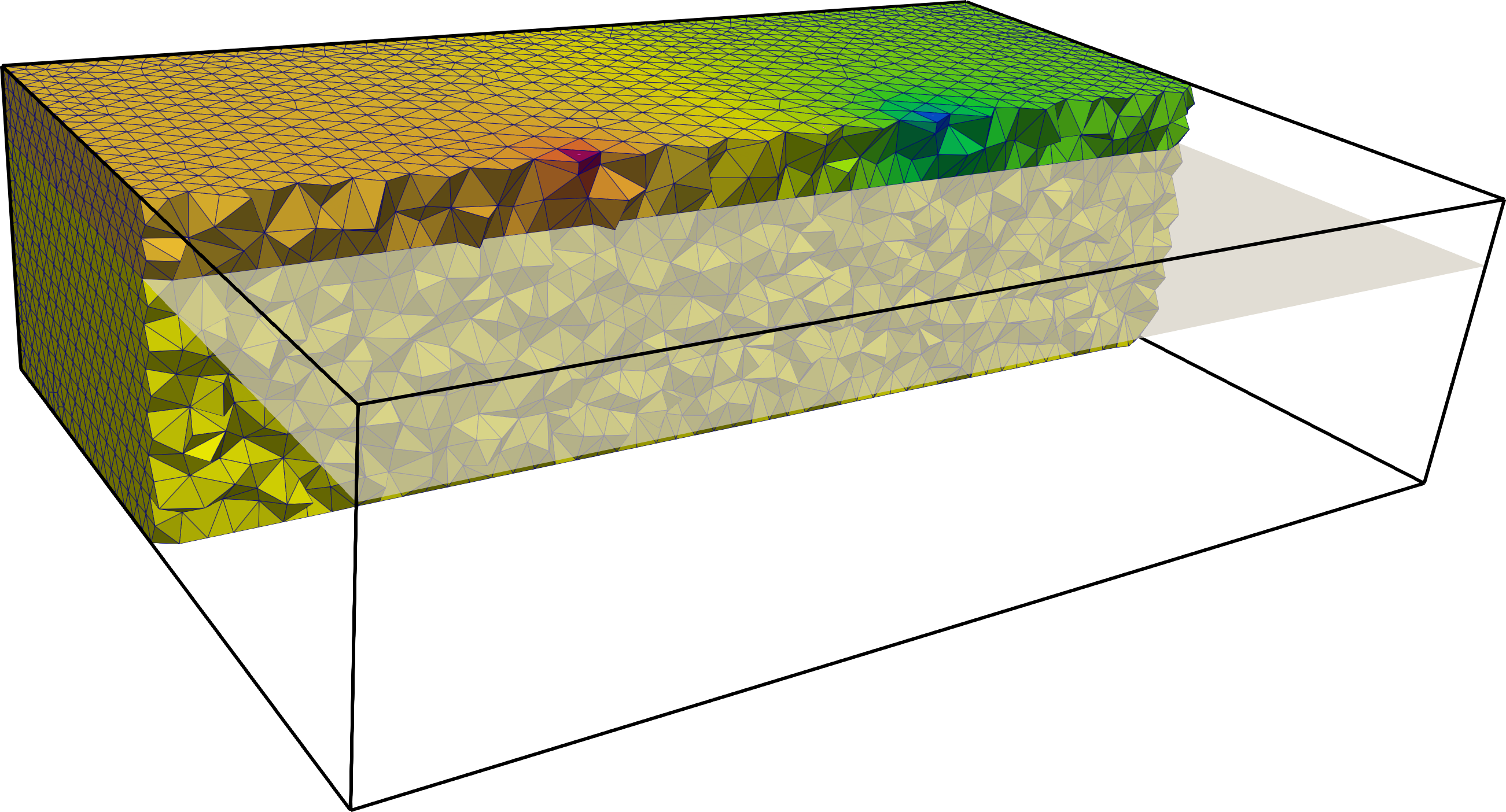}%
    \hspace*{0.05\textwidth}%
    \includegraphics[width=0.475\textwidth]{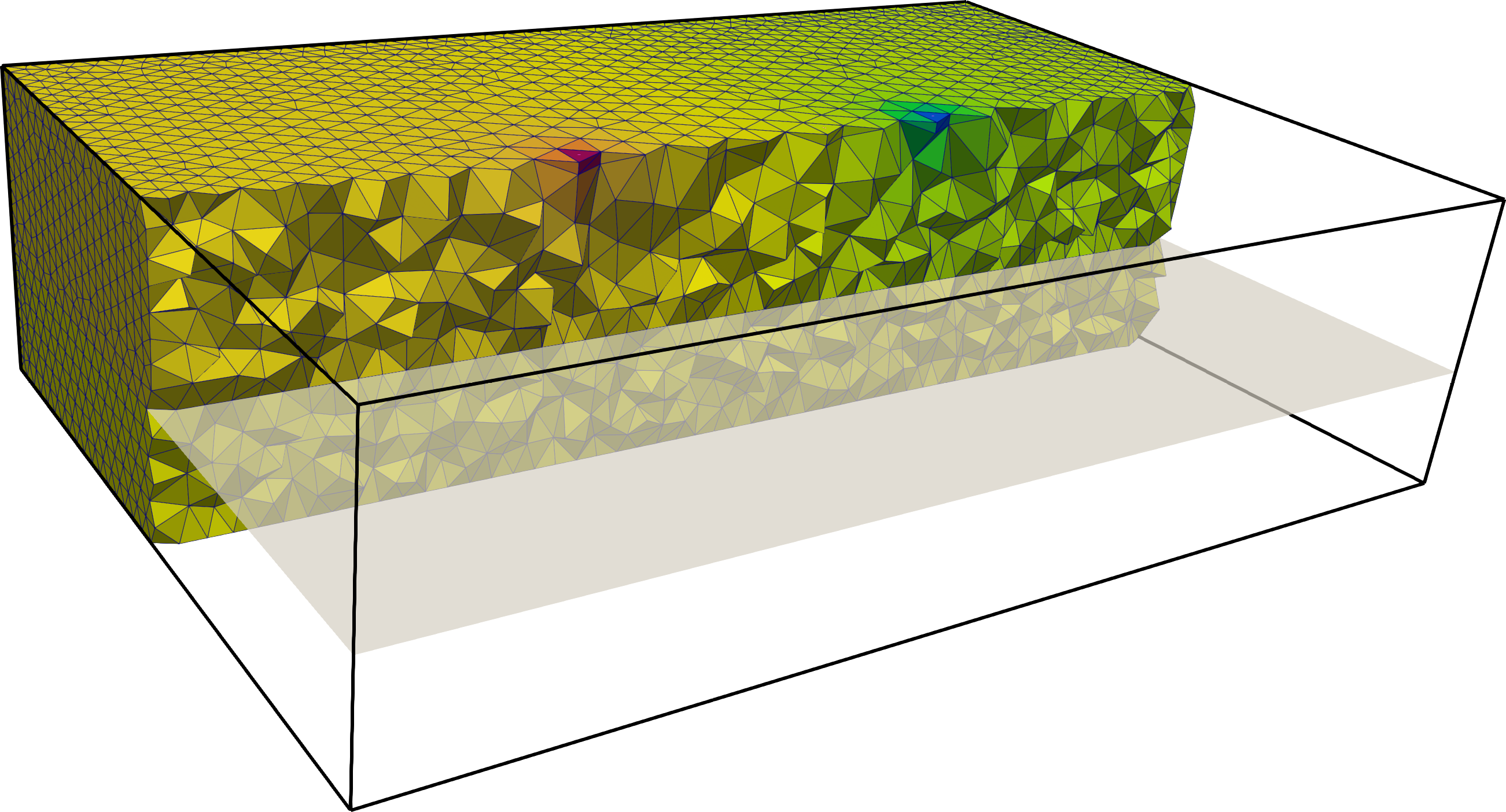}\\
    \includegraphics[width=0.475\textwidth]{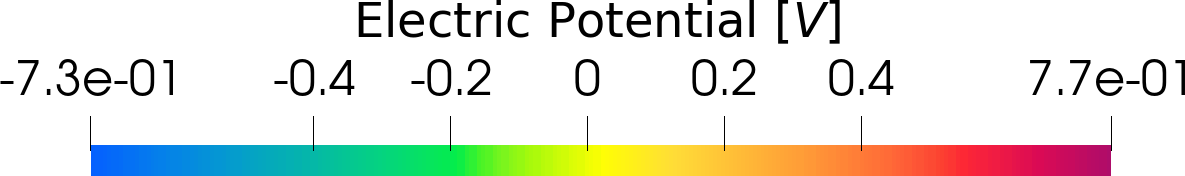}%
    \hspace*{0.05\textwidth}%
    \includegraphics[width=0.475\textwidth]{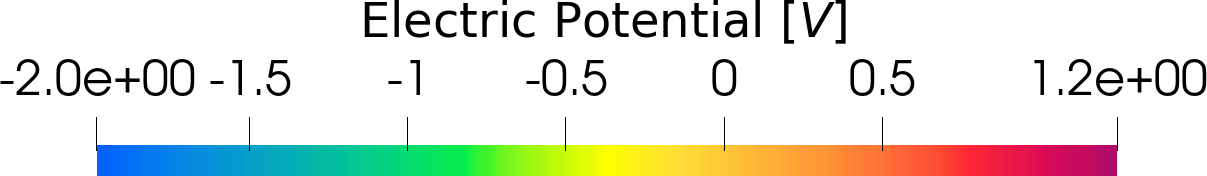}%
    \caption{Results for the example in Subsection \ref{subsec:itact_liner}. Electric potential in the modelled domain with the liner at depth $h=3\sib{\centi\meter}$ (left) and $h=8\sib{\centi\meter}$ (right).}
    \label{fig:experiment1potential}
\end{figure}

\begin{figure}
    \centering
    \includegraphics[width=0.475\textwidth]{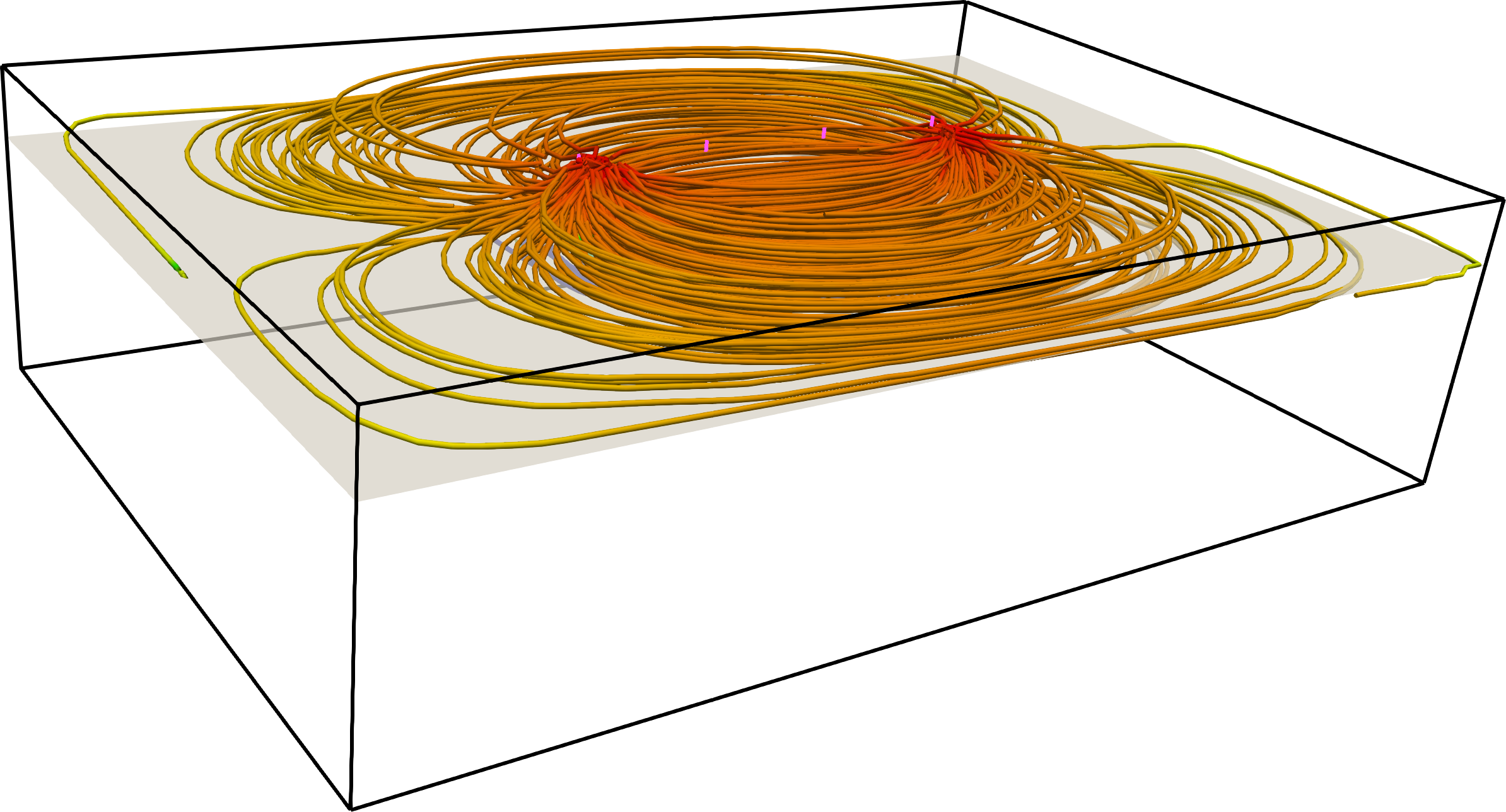}%
    \hspace*{0.05\textwidth}%
    \includegraphics[width=0.475\textwidth]{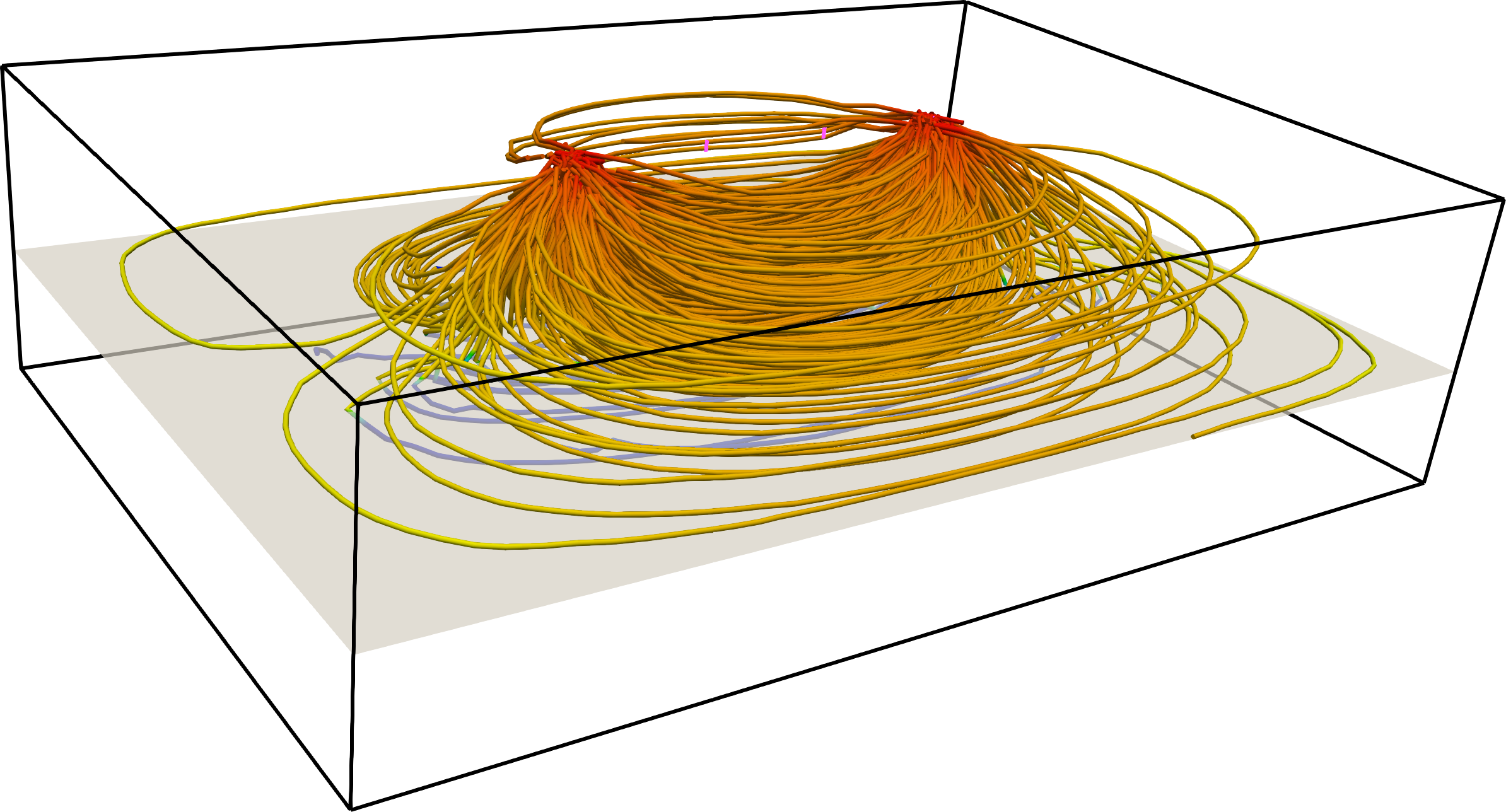}\\
    \includegraphics[width=0.475\textwidth]{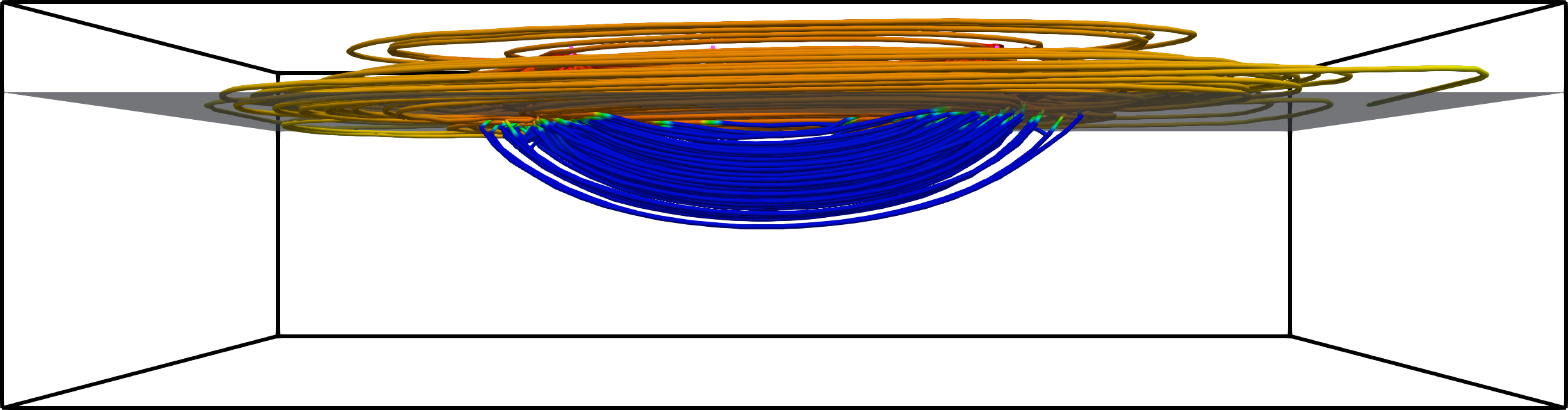}%
    \hspace*{0.05\textwidth}%
    \includegraphics[width=0.475\textwidth]{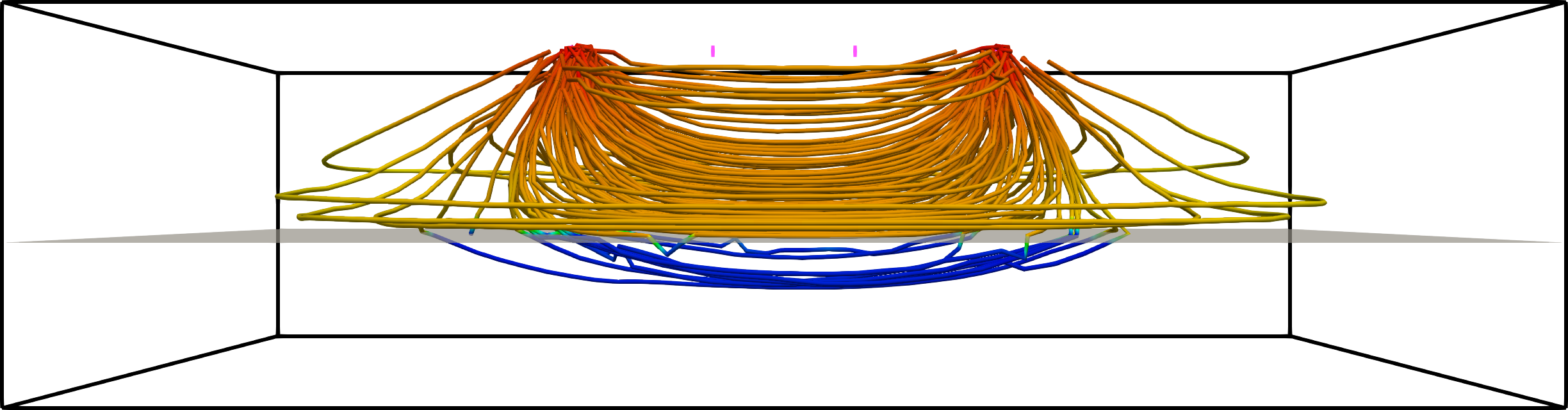}\\
    \includegraphics[width=0.475\textwidth]{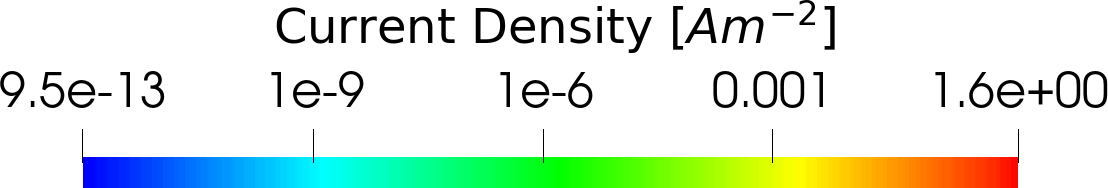}%
    \hspace*{0.05\textwidth}%
    \includegraphics[width=0.475\textwidth]{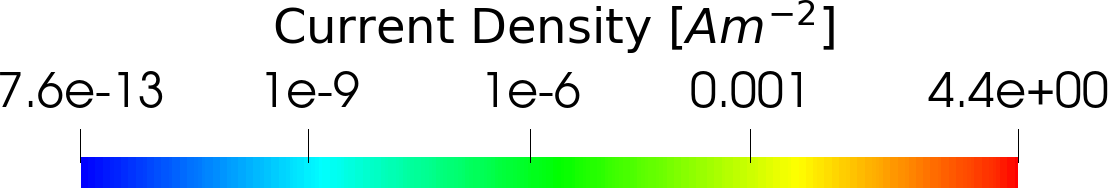}
    \caption{Results for the example in Subsection \ref{subsec:itact_liner}. Current lines and the associated current density in the modelled domain with the liner at depth $h=3\sib{\centi\meter}$ (left) and $h=8\sib{\centi\meter}$ (right).}
    \label{fig:experiment1current}
\end{figure}

Some of the numerical solutions computed with MPFA are illustrated in Figure
\ref{fig:experiment1potential} with the liner at different depths and for an
electrode spacing $a=3\sib{\centi\metre}$. The influence of the liner on the
electric potential is clearly visible, but becomes less evident when the
liner is farther from the electrodes. Below the liner the potential is essentially homogeneous. Figure \ref{fig:experiment1current} shows, for the same setting, the current lines with the associate current density, which is injected in $\Upsilon$ from $\gamma_1$ and is drawn from $\gamma_4$. Again, it is clearly visible how the liner practically confines current circulation in the upper part of the domain only.

\subsection{A box-shaped liner}\label{subsec:broken_liner}

\begin{figure}
    \centering
    \raisebox{0.\textwidth}{\resizebox{0.475\textwidth}{!}{\fontsize{0.75cm}{2cm}\selectfont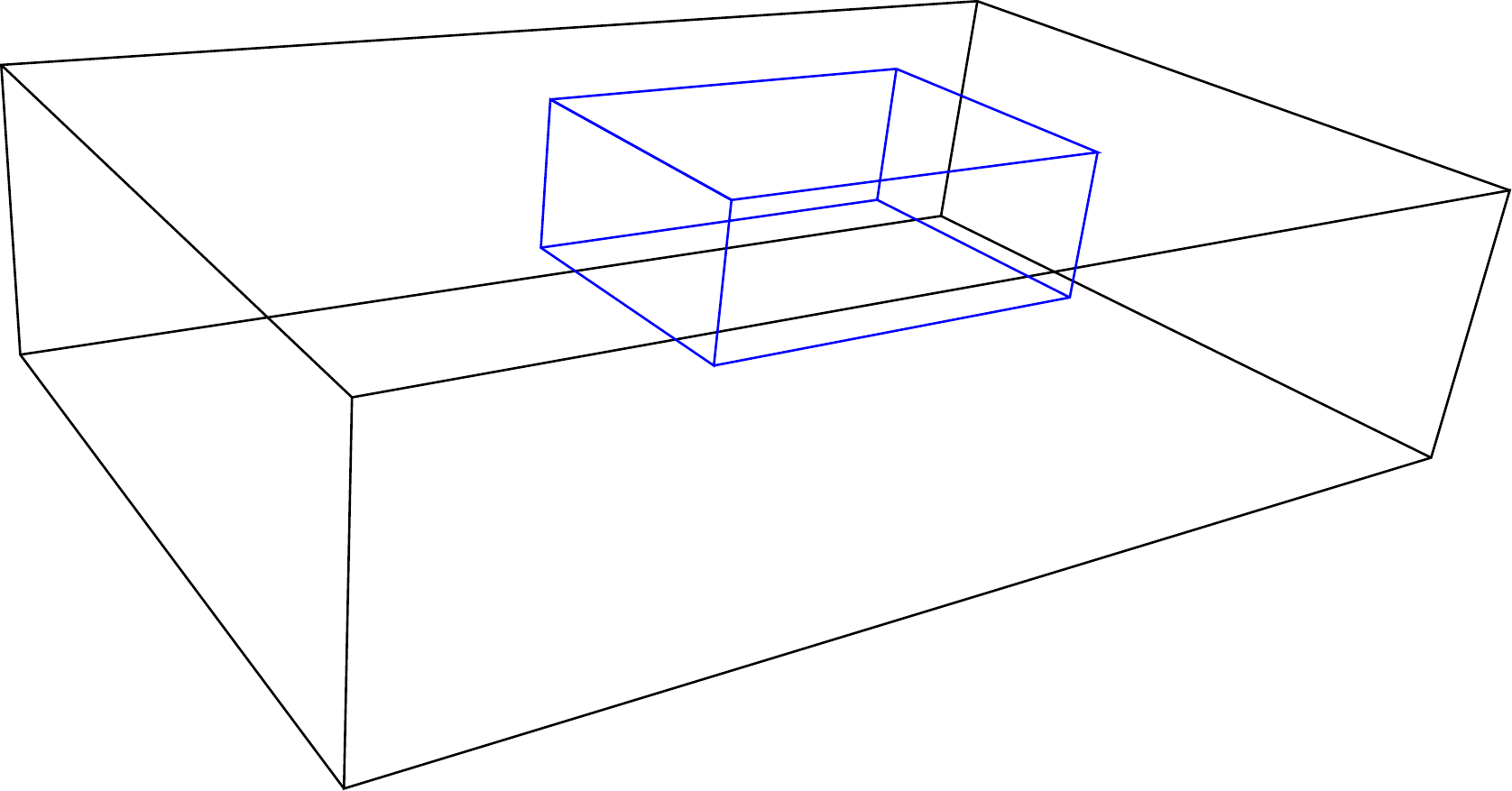}}%
    \hspace*{0.05\textwidth}%
    \raisebox{0.\textwidth}{\resizebox{0.475\textwidth}{!}{\fontsize{0.75cm}{2cm}\selectfont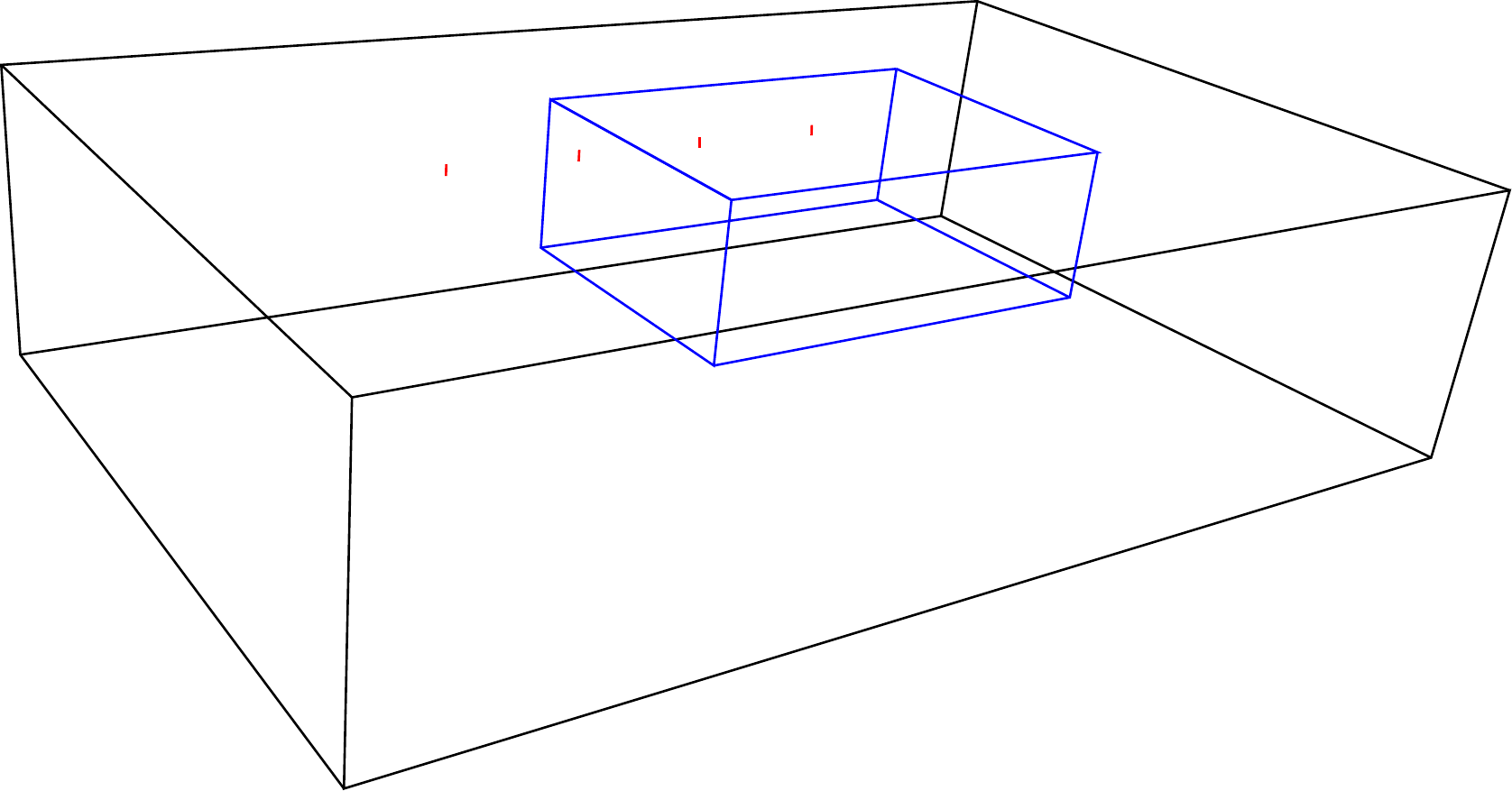}}%
    \caption{
    Domains for the example in Subsection \ref{subsec:broken_liner}:
    electrode
    configurations for \textit{case ii} (left) and \textit{case iii} (right).
    Note the presence of the hole $\eta$ in the bottom surface of the liner. See text for details.}
    \label{fig:domain_box}
\end{figure}

\begin{figure}
    \centering
    \includegraphics[width=0.24\textwidth,height=0.225\textwidth]{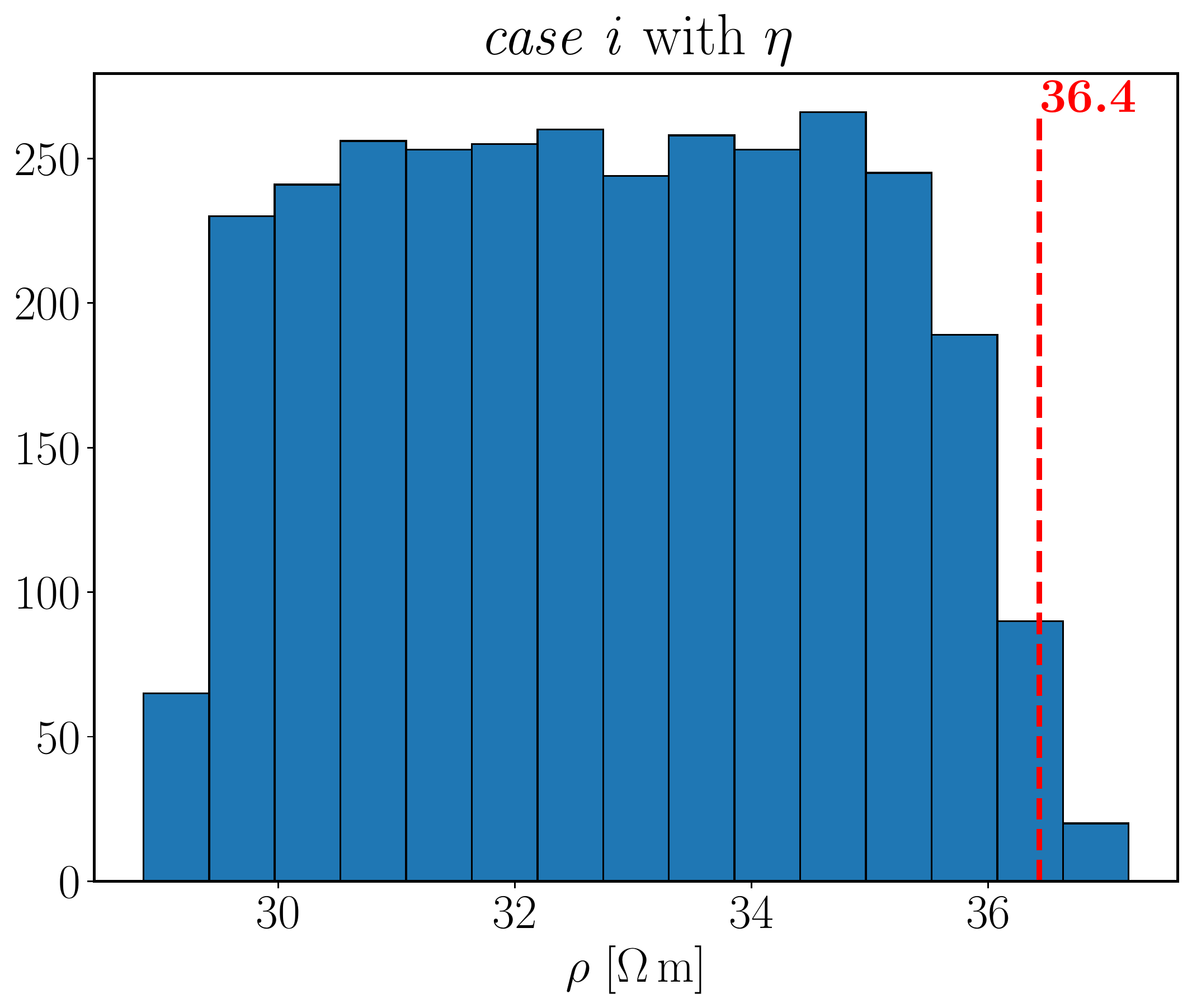}%
    \hspace*{0.01333\textwidth}%
    \includegraphics[width=0.24\textwidth,height=0.225\textwidth]{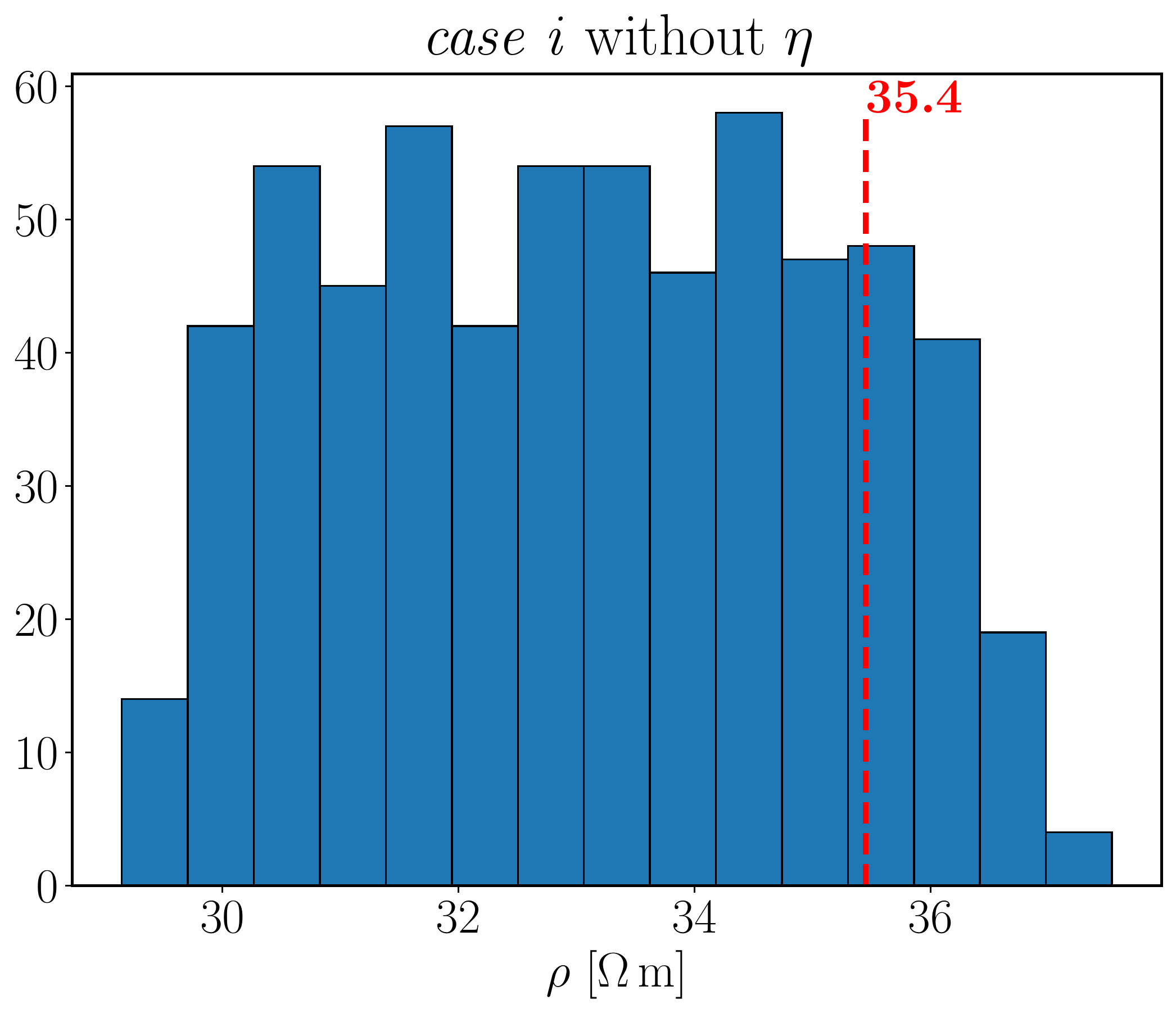}%
    \hspace*{0.01333\textwidth}%
    \includegraphics[width=0.24\textwidth,height=0.225\textwidth]{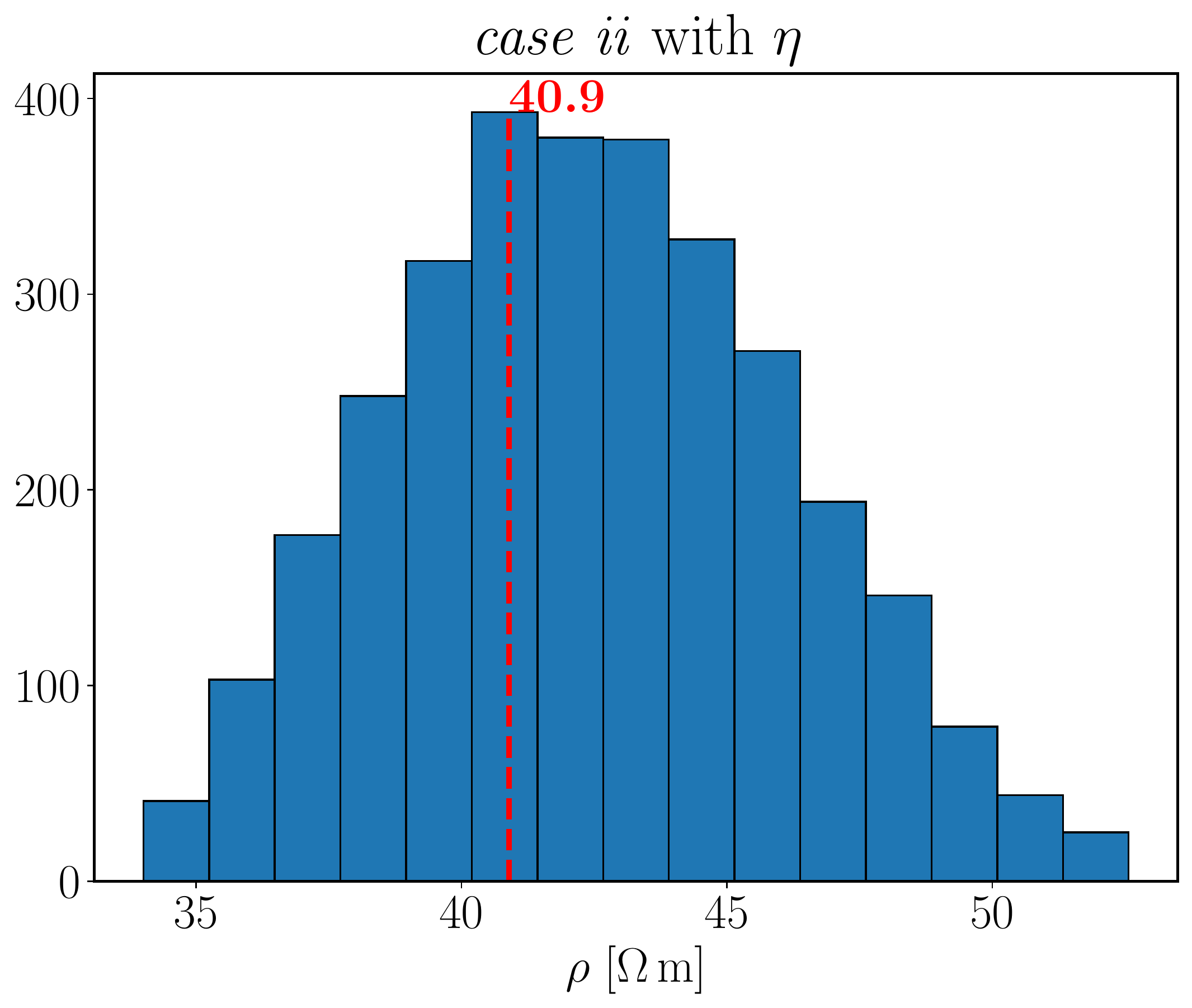}%
    \hspace*{0.01333\textwidth}%
    \includegraphics[width=0.24\textwidth,height=0.225\textwidth]{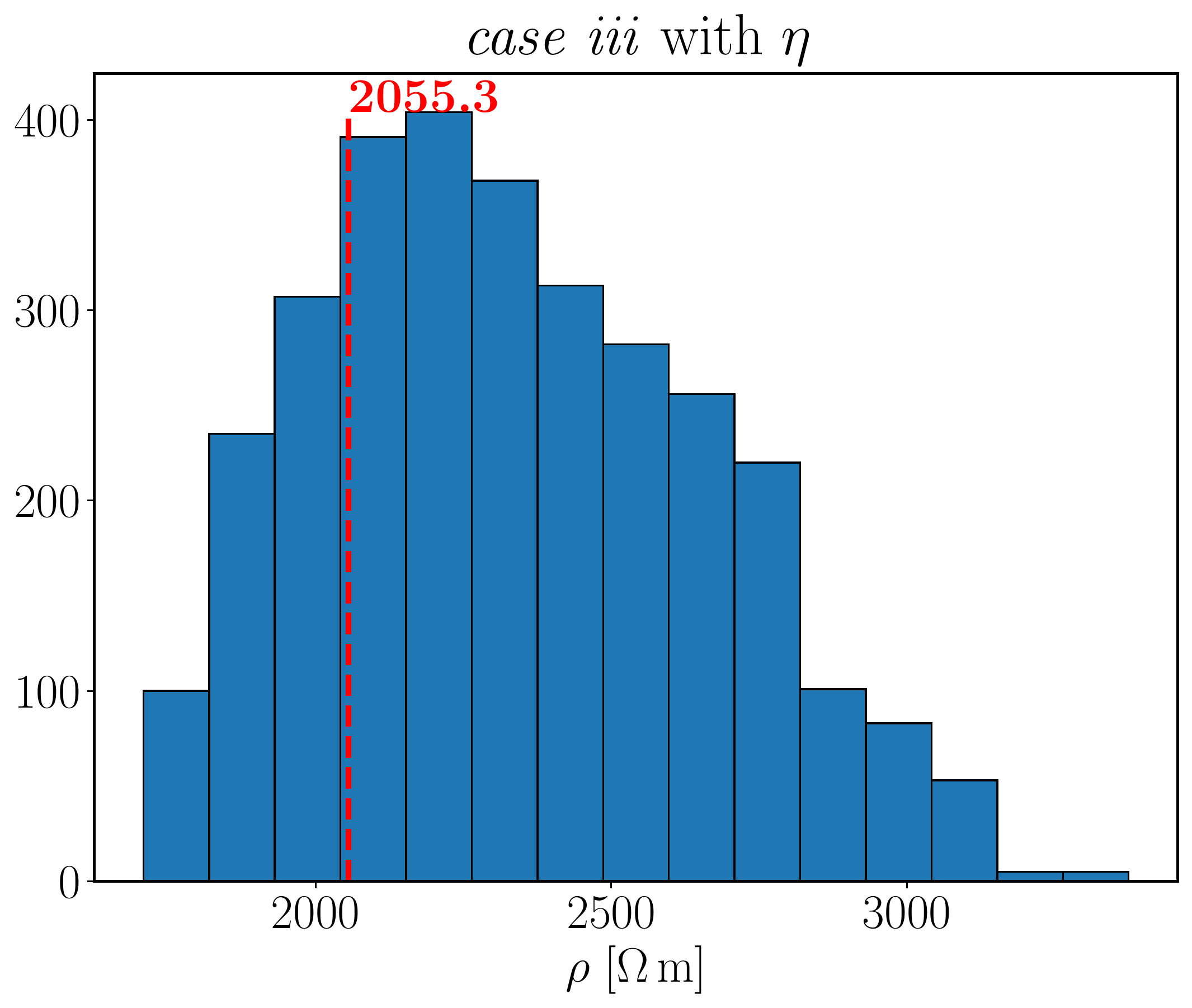}%
    \caption{Modelling results for different setups in Subsection 4.2. Histograms of the apparent resistivities modelled by taking into account uncertainties possibly affecting lab setup parameters. Red dashed lines are mean measured resistivities. See text for details.}
    \label{fig:hist}
\end{figure}

In this section we mainly want to evaluate the effect of a hole in the liner and its
impact on the resulting apparent resistivity according to different electrode deployments. We consider the liner has a box shape, both with and without a hole $\eta$ in its bottom surface. The box has dimension $20 \times 14
 \times 7.2 \sib{\cubic\centi\meter}$
and minimum coordinates in $(20.5, 10.5, 7.2) \sib{\centi\meter}$, the hole has
diameter of
$0.5\sib{\centi\meter}$ and is centred in $(25.5,
18.5, 7.2) \sib{\centi\meter}$ (Figure \ref{fig:domain_box}).
We consider three deployments for the electrodes in a Wenner-$\alpha$ configuration:
\textit{case i}) all electrodes are within the box with a spacing of
$3\sib{\centi\meter}$; \textit{case ii}) just the current electrodes $C_1$ is placed outside the box; \textit{case iii}) both the current electrode $C_1$ and the potential electrode $P_1$ are placed outside the box. For these last two cases the electrode spacing is set to $6\sib{\centi\meter}$. Please note that, for this set of experiments, tap water conductivity is estimated in the lab to be equal to $\sigma = 1/24$ \sib{\siemens\per\meter} and the numerical simulations are performed only with the MPFA scheme. The other parameters of the problem are the same of the previous example. In addition, by taking into account possible uncertainties in the laboratory setup, we perform several simulations with the variation of the water level of $\pm 3 \sib{\milli\meter}$, the variation of the water resistivity of $\pm 2 \sib{\ohm\meter}$, the variation of the radius of $\eta$ of $\pm 0.15 \sib{\milli\meter}$, and the shift of the electrodes along $x$ and $y$, with respect to the box, of $\pm 6 \sib{\milli\meter}$.
 Figure \ref{fig:hist} shows that, for all the cases, the measured apparent resistivities fall into the range of the modelled values. We point out that the measured standard deviations are always lower than 1\% and that, for \textit{case ii}) and  \textit{case iii}), it was not possible to obtain reliable measurements when there is no hole in the liner.

\begin{figure}
    \centering
    \includegraphics[width=0.475\textwidth]{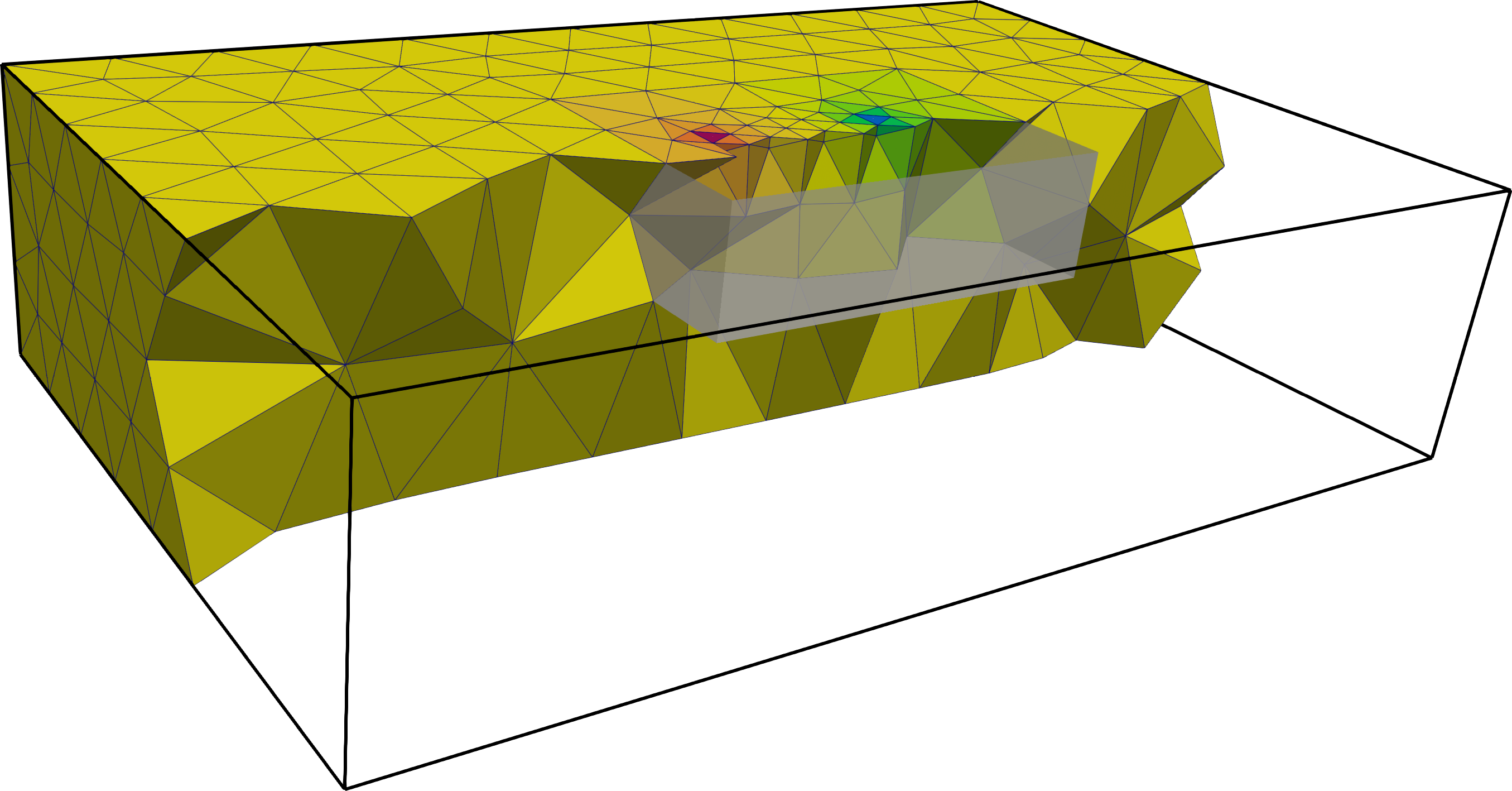}%
    \hspace*{0.05\textwidth}%
    \includegraphics[width=0.475\textwidth]{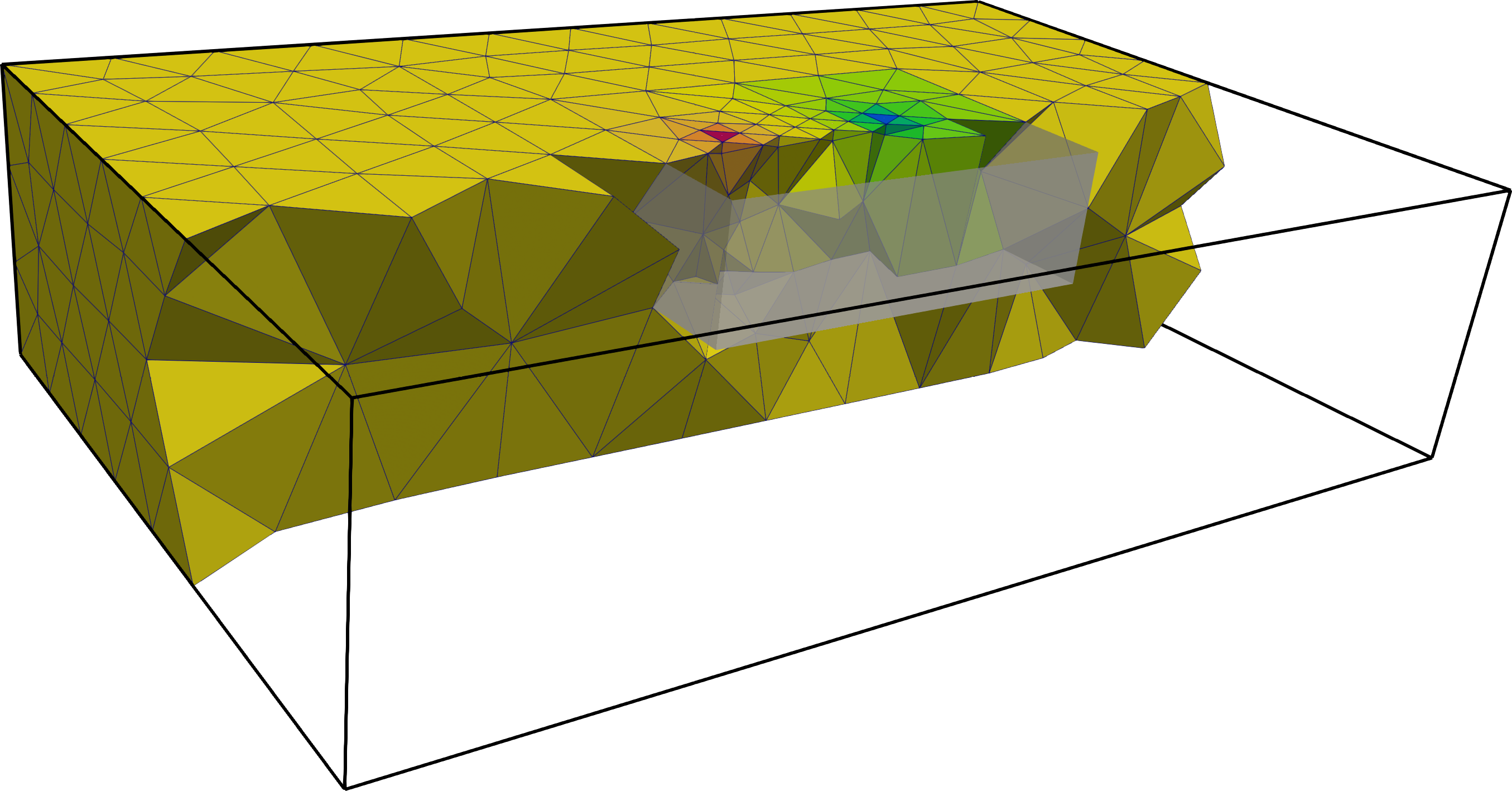}
    \\
    \includegraphics[width=0.475\textwidth]{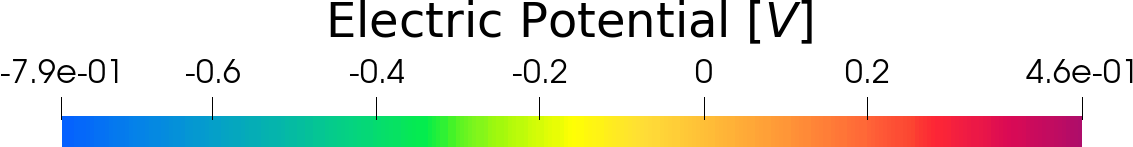}%
    \\[0.75cm]
    \includegraphics[width=0.475\textwidth]{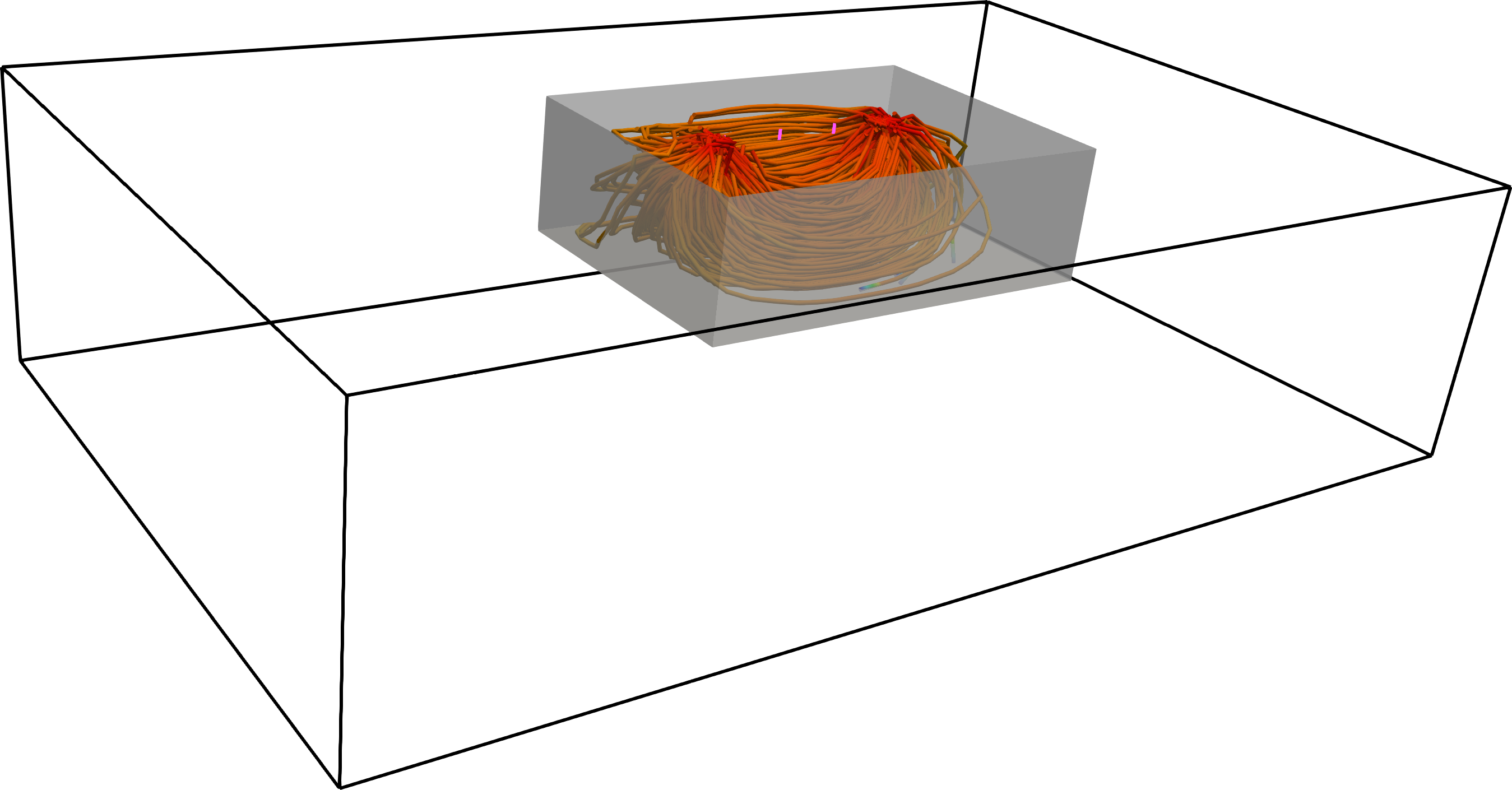}%
    \hspace*{0.05\textwidth}%
    \includegraphics[width=0.475\textwidth]{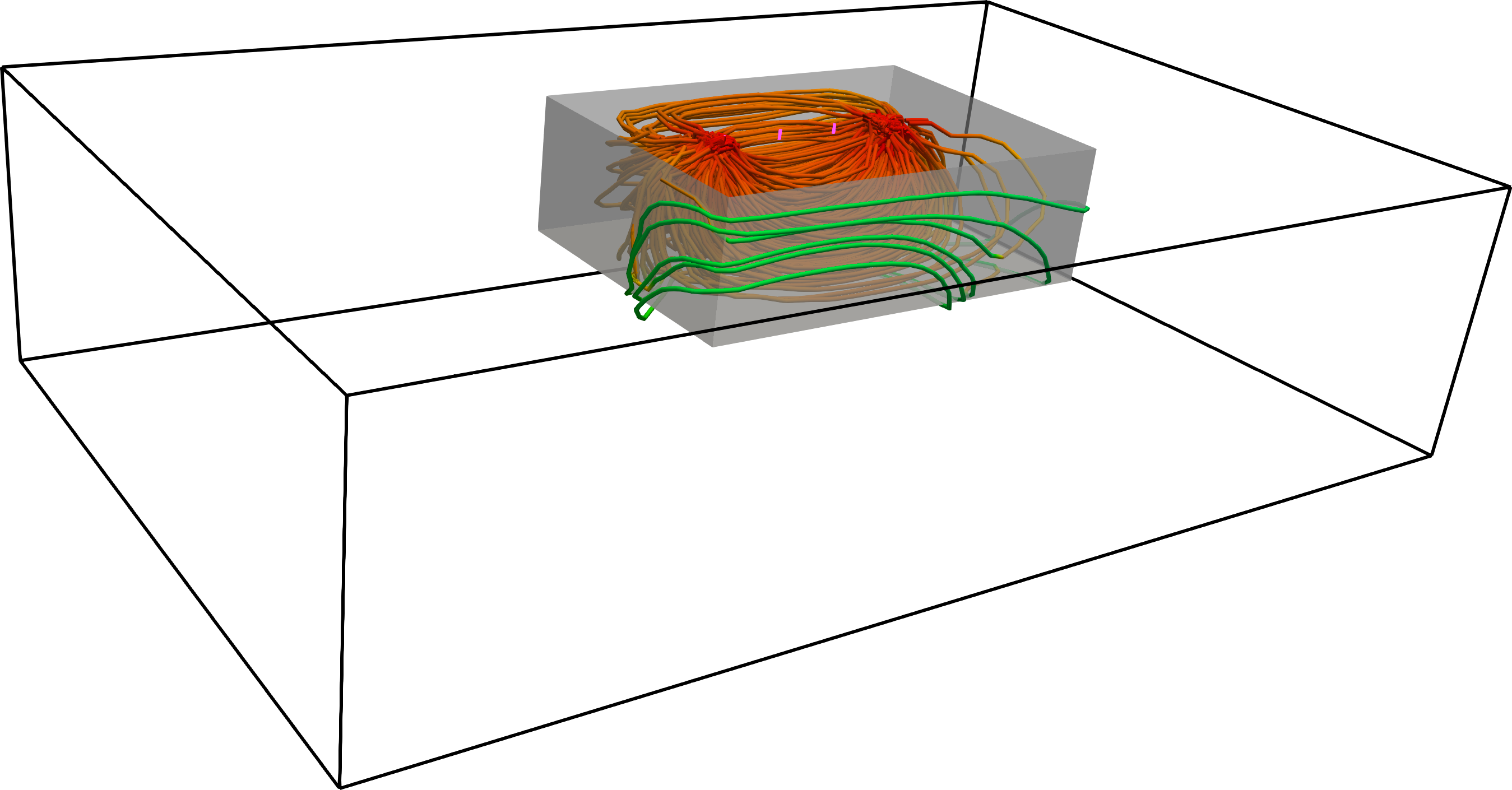}
    \\
    \includegraphics[width=0.475\textwidth]{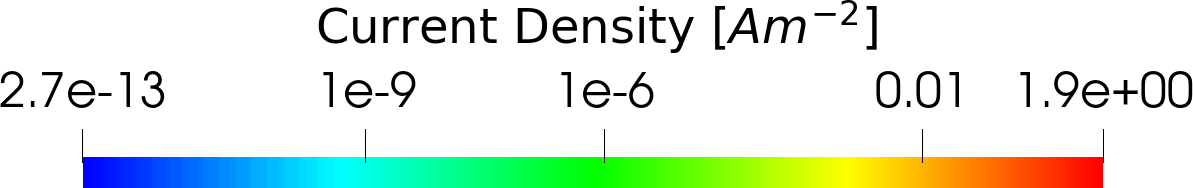}%
    \caption{Results for \textit{case i} described in Subsection \ref{subsec:broken_liner} without (left) and with (right) the hole in the liner. The top and the bottom images show the modelled electric potential and current lines with the associated current density, respectively.}
    \label{fig:experiment2c}
\end{figure}

Regarding \textit{case i}, it is interesting to note that the presence of the hole $\eta$ does not significantly affect the value of the apparent resistivity $\rho_a$, which implies that this electrode configuration may not be helpful for detecting defects of the liner. Accordingly, the numerical solutions for this configuration do not show relevant differences either the hole is present or not (Figure \ref{fig:experiment2c}). The current circulates mostly inside the box-shaped liner, where we observe the variations of the electric potential.

\begin{figure}
    \centering
    \includegraphics[width=0.475\textwidth]{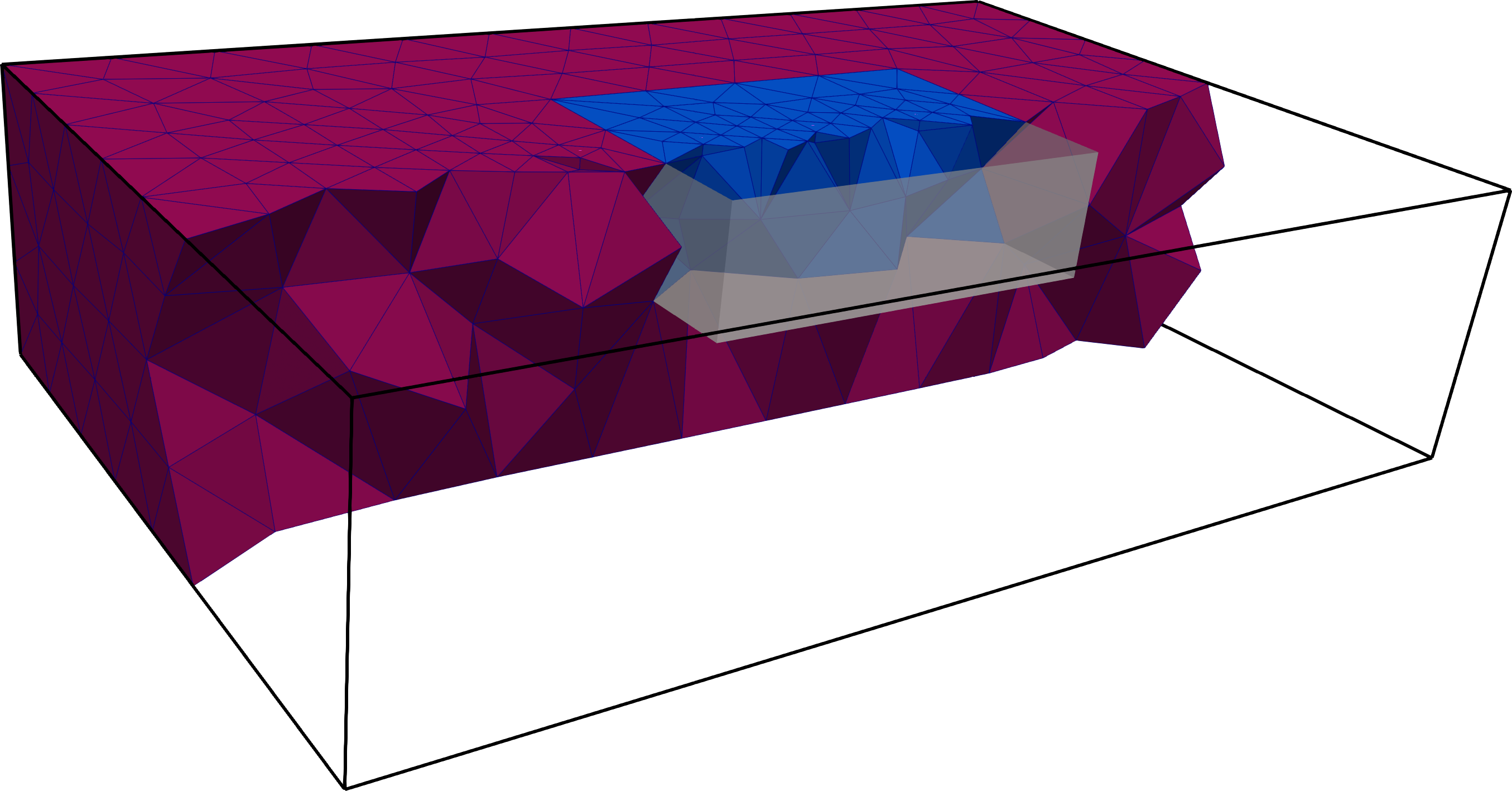}%
    \hspace*{0.05\textwidth}%
    \includegraphics[width=0.475\textwidth]{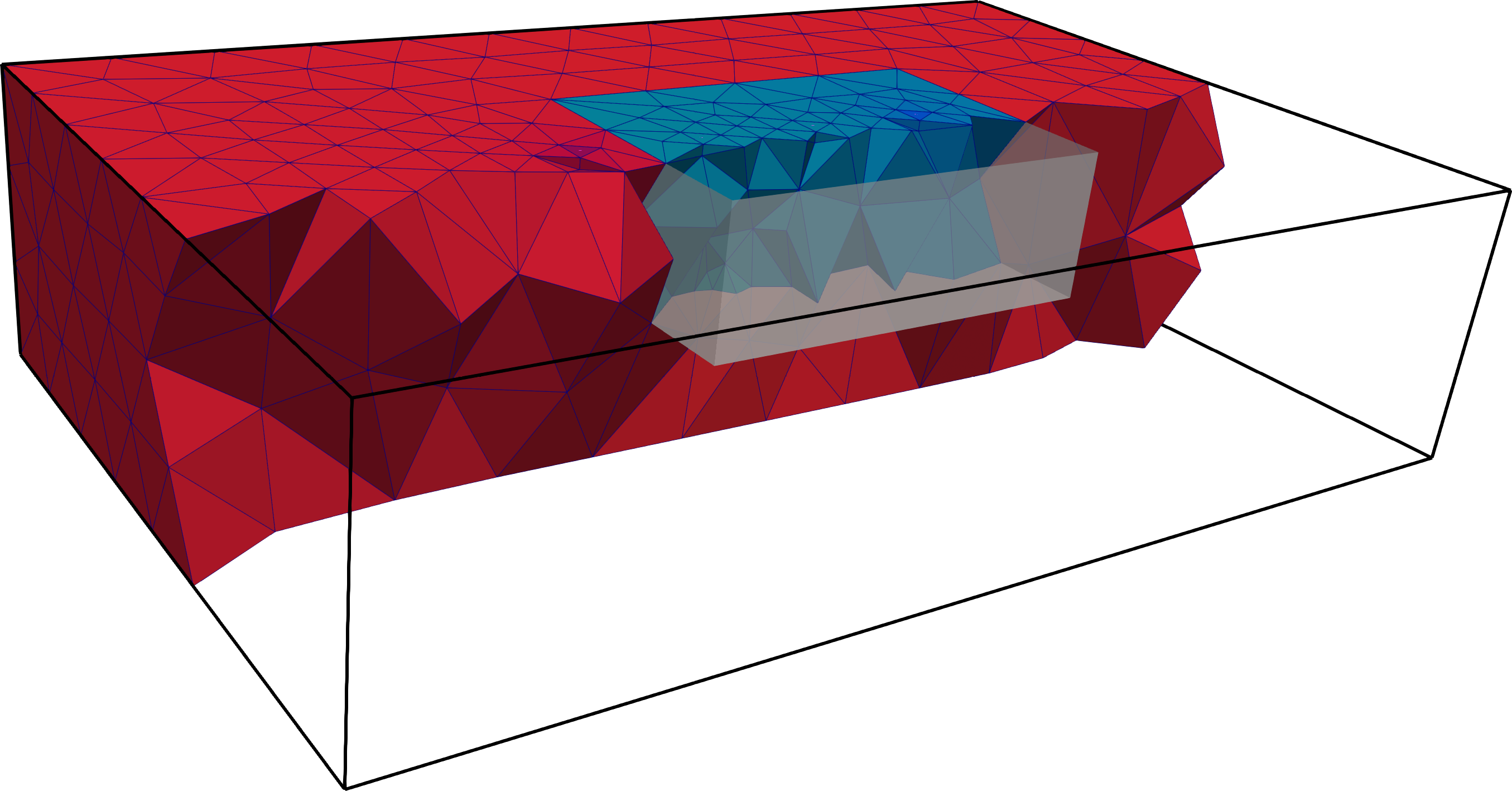}\\
    \includegraphics[width=0.475\textwidth]{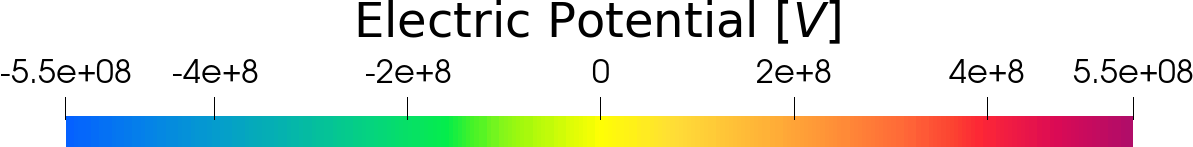}%
    \hspace*{0.05\textwidth}%
    \includegraphics[width=0.475\textwidth]{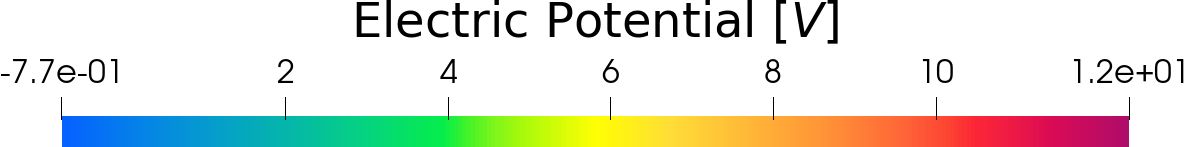}\\[0.75cm]
    \includegraphics[width=0.475\textwidth]{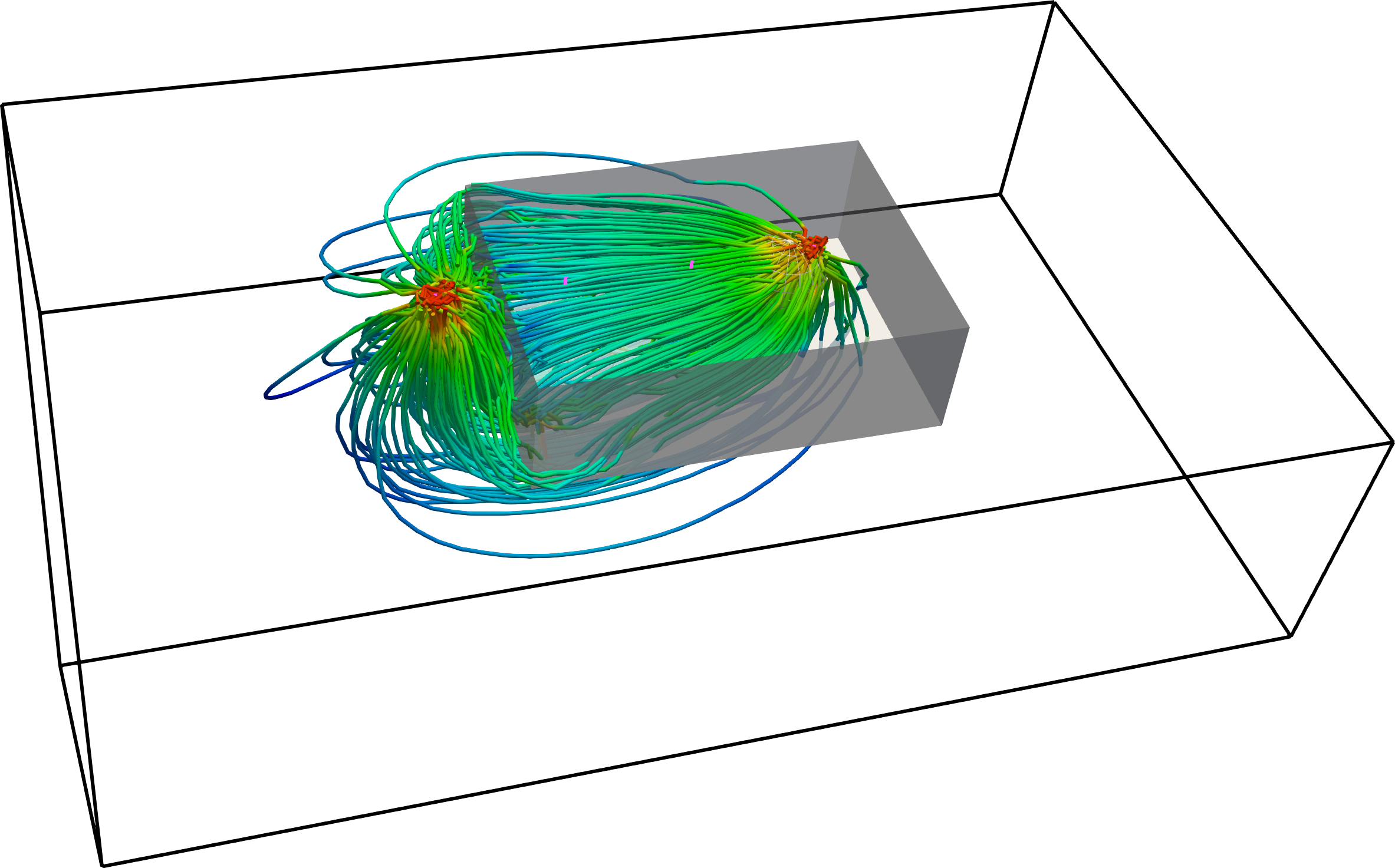}%
    \hspace*{0.05\textwidth}%
    \includegraphics[width=0.475\textwidth]{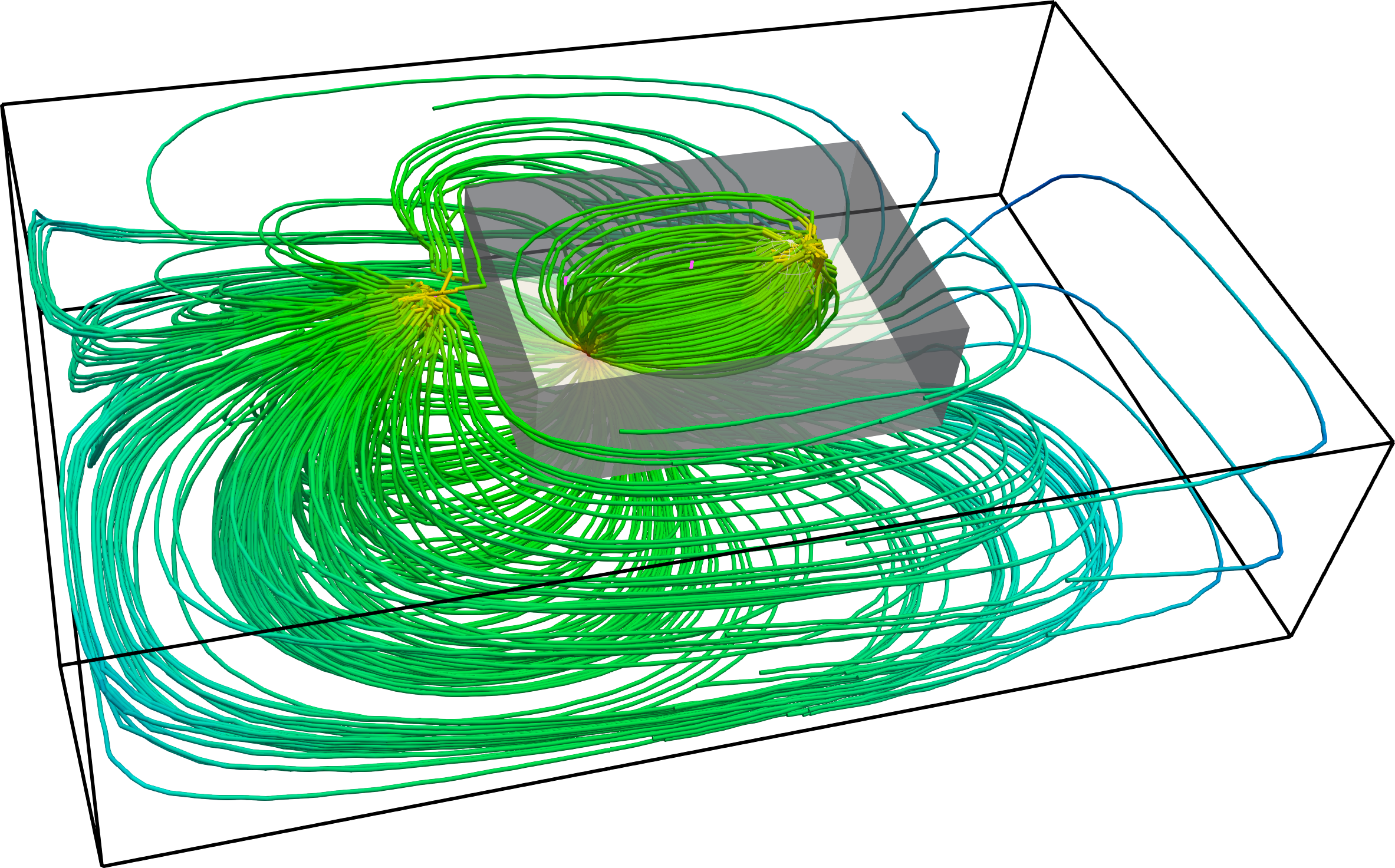}\\
    \includegraphics[width=0.475\textwidth]{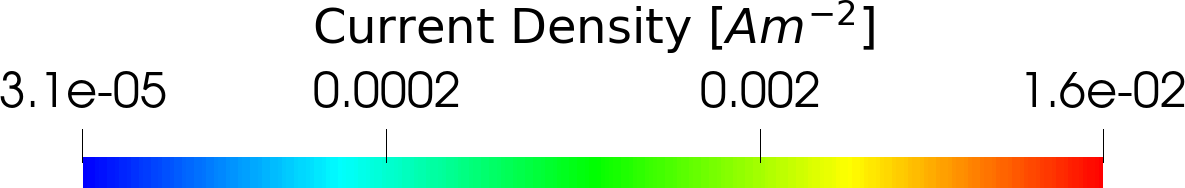}%
    \hspace*{0.05\textwidth}%
    \includegraphics[width=0.475\textwidth]{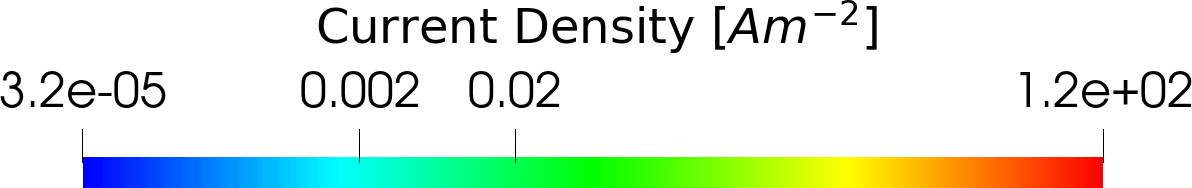}
    \caption{Results for \textit{case ii} described in Subsection \ref{subsec:broken_liner} without (left) and with (right) the hole in the liner. The top and the bottom images show the modelled electric potential and current lines with the associated current density, respectively. Please note the different colormap ranges for each image.}
    \label{fig:experiment2b}
\end{figure}

\begin{figure}
    \centering
    \includegraphics[width=0.475\textwidth]{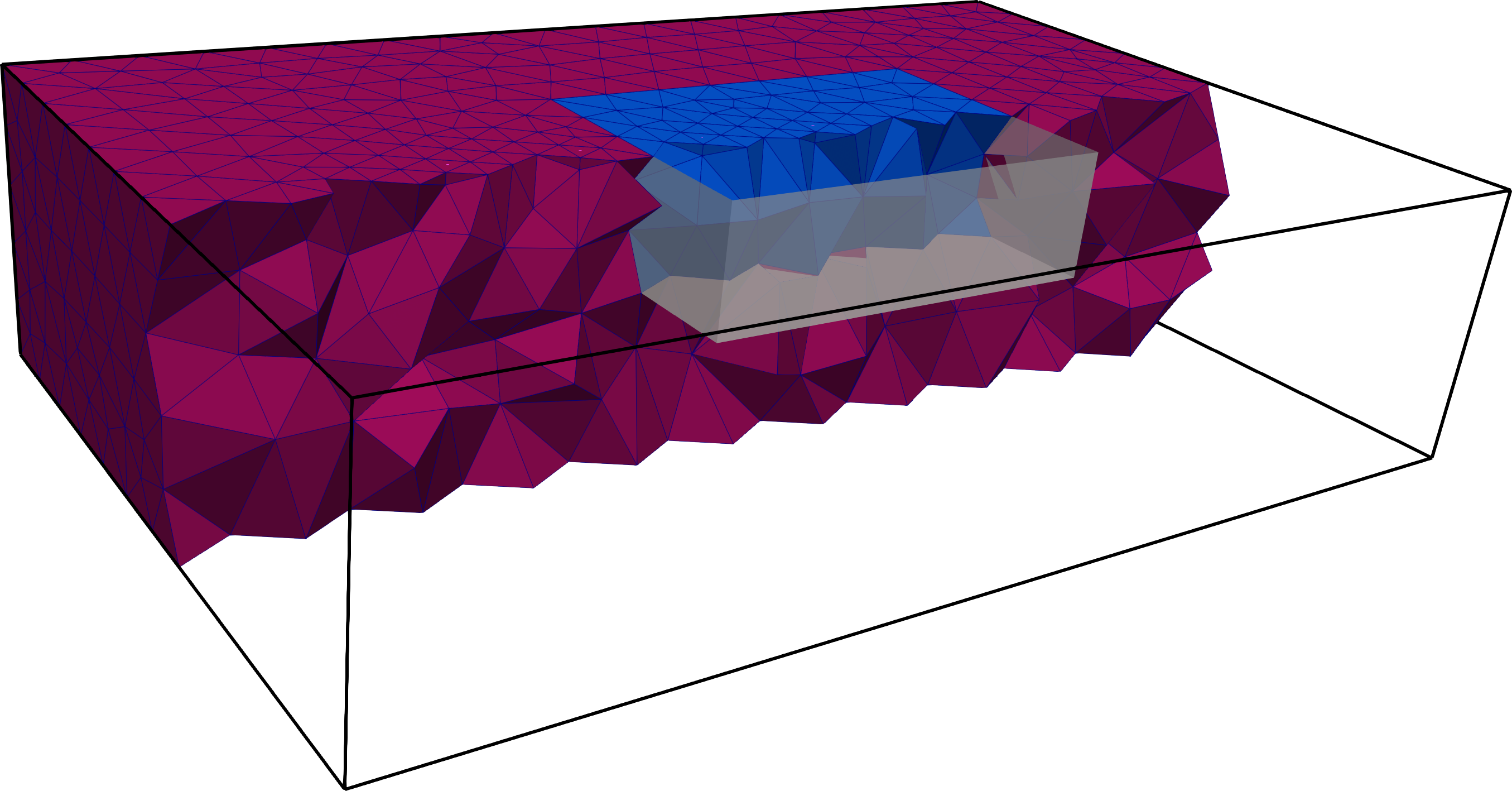}%
    \hspace*{0.05\textwidth}%
    \includegraphics[width=0.475\textwidth]{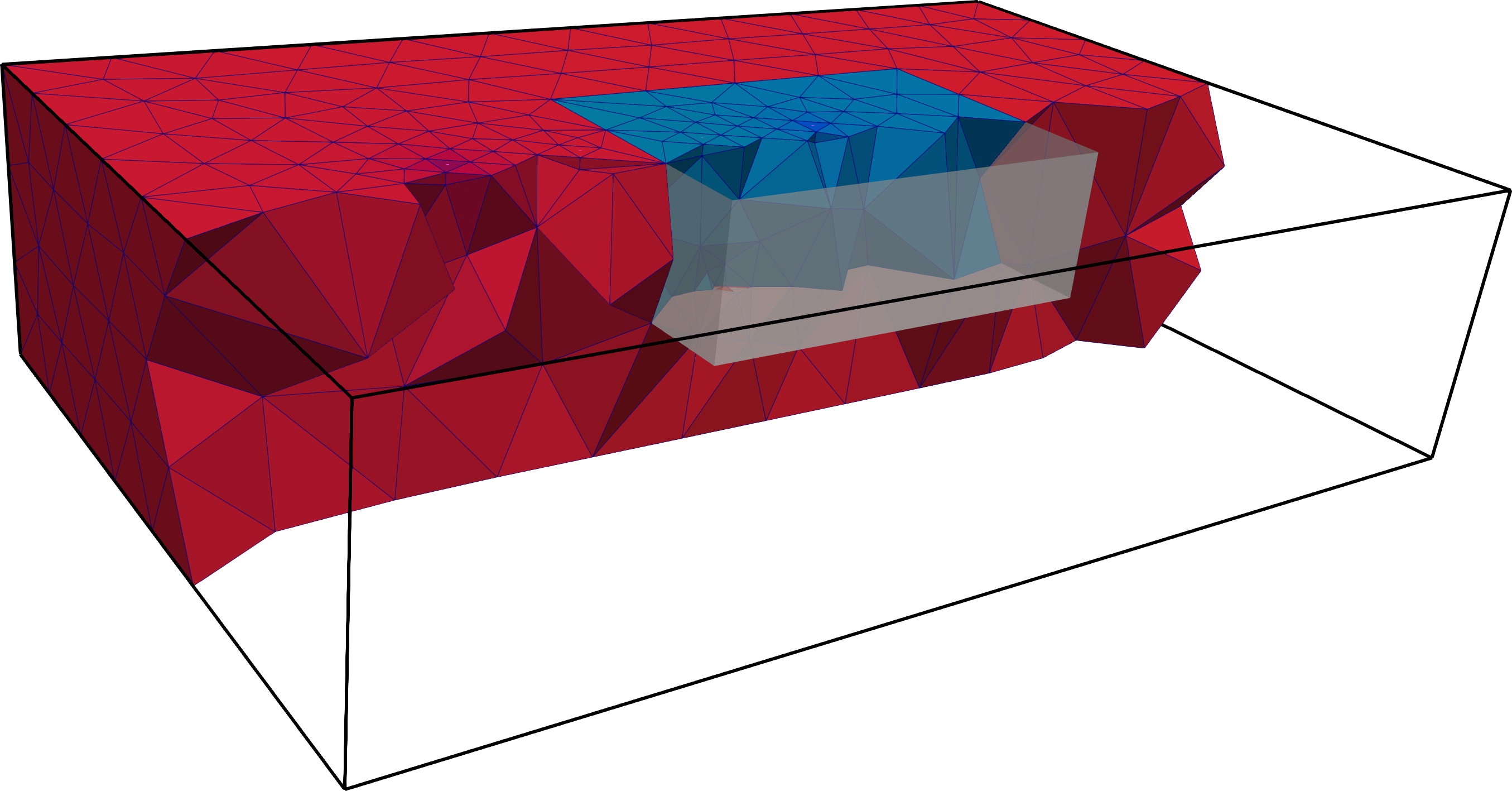}\\
    \includegraphics[width=0.475\textwidth]{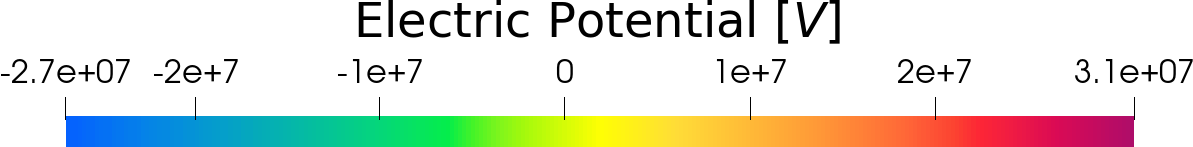}%
    \hspace*{0.05\textwidth}%
    \includegraphics[width=0.475\textwidth]{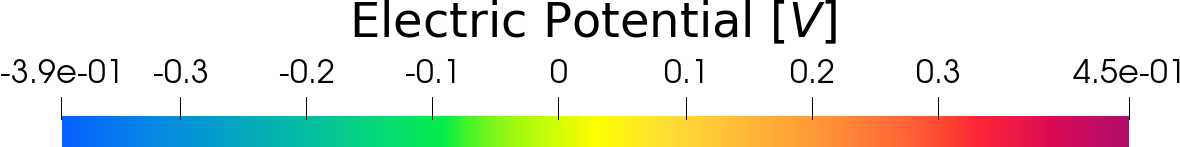}\\[0.75cm]
    \includegraphics[width=0.475\textwidth]{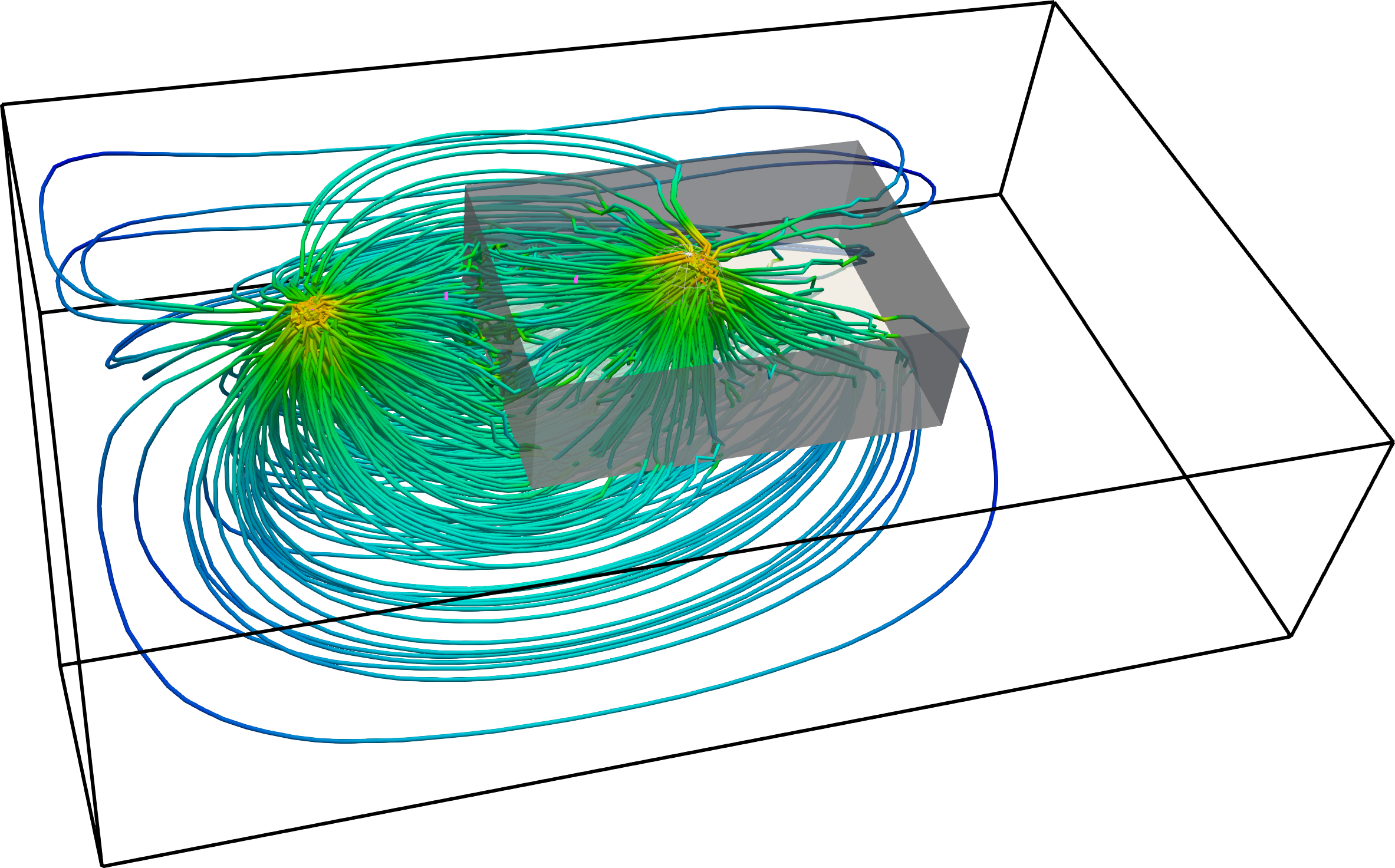}%
    \hspace*{0.05\textwidth}%
    \includegraphics[width=0.475\textwidth]{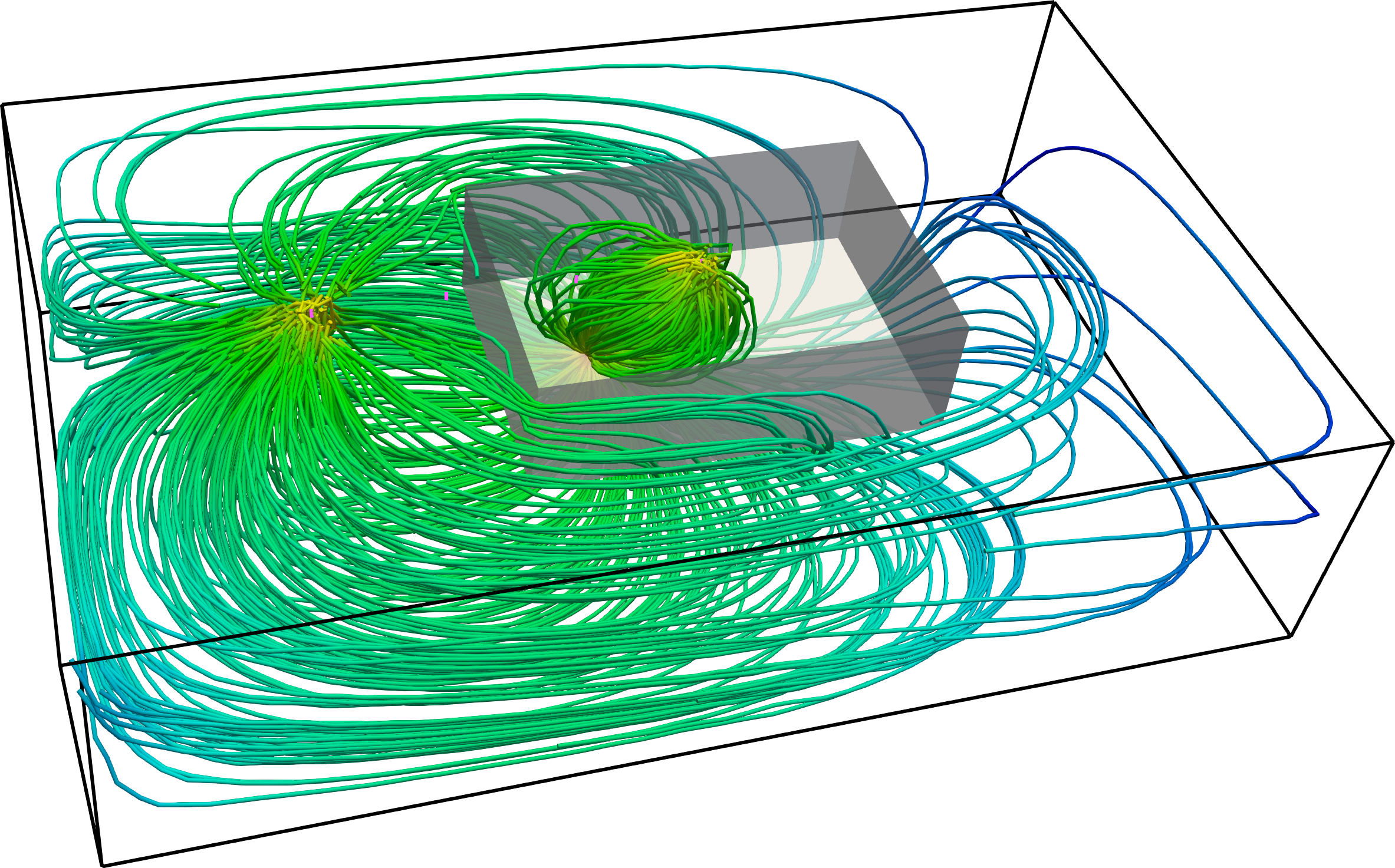}\\
    \includegraphics[width=0.475\textwidth]{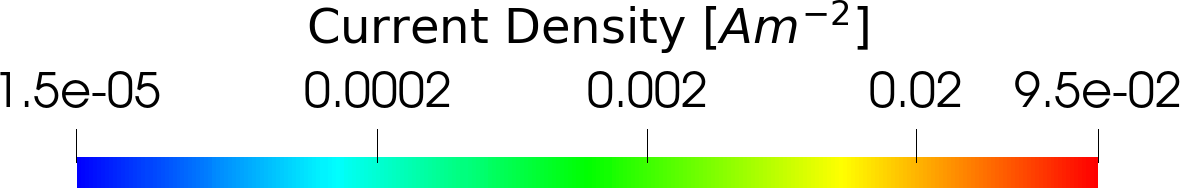}%
    \hspace*{0.05\textwidth}%
    \includegraphics[width=0.475\textwidth]{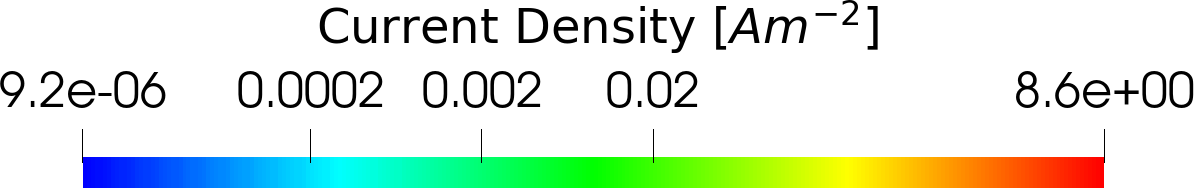}
    \caption{Results for \textit{case iii} described in Subsection \ref{subsec:broken_liner} without (left) and with (right) the hole in the liner. The top and the bottom images show the modelled electric potential and current lines with the associated current density, respectively. Please note the different colormap ranges for each image.}
    \label{fig:experiment2a}
\end{figure}

The modelling results for \textit{case ii} and \textit{case iii} are similar (Figure \ref{fig:experiment2b} and \ref{fig:experiment2a}), except for the fact that the estimated apparent resitivity is approximately two orders of magnitude higher for the latter (Figure Figure \ref{fig:hist}). This is due to the presence of the high-resistivity liner between the voltage electrodes in \textit{case iii}. We generally observe two equi-potential regions, one inside and the other outside the box-shaped liner. When the liner is perforated, the variation of $\varphi$ across $\lambda$ is obviously lower and the main path of the current is from $C_1$ to $C_2$ through $\eta$. On the contrary, when there is no hole, there in no clear preferential direction for current circulation.
For sake of completeness, we report that, for \textit{case ii} and \textit{case iii} with no hole, the estimated apparent resitivities are $26.9 \sib{\ohm\meter}$ and $5.48e10 \sib{\ohm\meter}$, respectively. Such a large gap between those values is obviously related to the position of both current and voltage electrodes with respect to the liner.

\section{Conclusions}\label{sec:conclusions}

In this work we have presented a mixed-dimensional mathematical model to obtain the electric potential and current density in direct current simulations. The model can handle multiple electrodes and a thin high-resistivity liner included in the domain. These objects have one or two of their dimensions that are very small and thus difficult to be represented exactly by equi-dimensional grids in the simulations. The mixed-dimensional approach used in this work approximates them as object with lower dimension: electrodes as one-dimensional and the liner as two-dimensional objects. New equations have been derived along with new interface conditions to couple the electrodes and the liner with the surrounding domain. To improve the efficiency of the simulation, the electrodes are placed in the background computational grid avoiding an excessive refinement around them. The numerical approximation relies on two cell-centred
finite-volume schemes: the two-point and multi-point flux approximations.
The mathematical model has been tested against laboratory experiments
giving reliable solutions, in particular when the multi-point flux approximation was considered. In the first set of laboratory experiments, we validated the code with an analytical relationship  and showed that the liner influences the measured apparent resistivity at depth much higher then the penetration depth of a Wenner-$\alpha$ array. In the second set of experiments, we analyzed the presence of a hole in the liner for different deployments of the electrodes and observed that some configurations are more favourable to detect possible defects.
We can conclude that the proposed mixed-dimensional model is a reliable tool for direct current simulations. The model can handle particular features that are very important when dealing with geoelectrical investigations of MSWLFs, where large variations of electrical resistivity can occur in very limited space.
The approach we developed could also be exported to other application fields, presenting similar governing equations and peculiarities.

\section*{Acknowledgements}

The authors wish to thank Carlo De Falco, Laura Longoni, Monica Papini,
and Luigi Zanzi at Politecnico di Milano for fruitful discussions. We also thank Francesco Ronchetti and Marco Sabattini at Università degli Studi di Modena e Reggio Emilia for their help with laboratory experiments.

\bibliographystyle{alpha}
\bibliography{biblio}

\newcommand{\etalchar}[1]{$^{#1}$}
\begin{thebibliography}{GKNW19}

\bibitem[Aav02]{Aavatsmark2002}
Ivar Aavatsmark.
\newblock An introduction to multipoint flux approximations for quadrilateral
  grids.
\newblock {\em Computational Geosciences}, 6:405--432, 2002.

\bibitem[Aav07a]{Aavatsmark2007a}
Ivar Aavatsmark.
\newblock Interpretation of a two-point flux stencil for skew parallelogram
  grids.
\newblock {\em Computational Geosciences}, 11(3):199--206, 2007.

\bibitem[Aav07b]{Aavatsmark2007}
Ivar Aavatsmark.
\newblock Multipoint flux approximation methods for quadrilateral grids.
\newblock In {\em The 9\textsuperscript{th} International Forum on Reservoir
  Simulation, Abu Dhabi}, pages 9--13, 2007.

\bibitem[AFS{\etalchar{+}}21]{aguzzoli2021inversion}
Alessandro Aguzzoli, Alessio Fumagalli, Anna Scotti, Luigi Zanzi, and Diego
  Arosio.
\newblock Inversion of synthetic and measured 3d geoelectrical data to study
  the geomembrane below a landfill.
\newblock In {\em 4th Asia Pacific Meeting on Near Surface Geoscience \&
  Engineering}, volume~1, pages 1--5. European Association of Geoscientists \&
  Engineers, 2021.

\bibitem[AHZA20]{Aguzzoli2020}
A.~Aguzzoli, A.~Hojat, L.~Zanzi, and D.~Arosio.
\newblock Two dimensional ert simulations to check the integrity of
  geomembranes at the base of landfill bodies.
\newblock 2020(1):1--5, 2020.

\bibitem[BBF{\etalchar{+}}20]{Berre2020a}
Inga Berre, Wietse~M. Boon, Bernd Flemisch, Alessio Fumagalli, Dennis
  Gl\"{a}ser, Eirik Keilegavlen, Anna Scotti, Ivar Stefansson, Alexandru
  Tatomir, Konstantin Brenner, Samuel Burbulla, Philippe Devloo, Omar Duran,
  Marco Favino, Julian Hennicker, I-Hsien Lee, Konstantin Lipnikov, Roland
  Masson, Klaus Mosthaf, Maria Giuseppina~Chiara Nestola, Chuen-Fa Ni, Kirill
  Nikitin, Philipp Sch\"{a}dle, Daniil Svyatskiy, Ruslan Yanbarisov, and
  Patrick Zulian.
\newblock Verification benchmarks for single-phase flow in three-dimensional
  fractured porous media.
\newblock {\em Advances in Water Resources}, 147, 2020.

\bibitem[BD03]{Binley2003}
Andrew Binley and William Daily.
\newblock The performance of electrical methods for assessing the integrity of
  geomembrane liners in landfill caps and waste storage ponds.
\newblock {\em Journal of Environmental and Engineering Geophysics},
  8(4):227--237, 2003.

\bibitem[BDOH00]{Bernstone2000}
C.~Bernstone, T.~Dahlin, T.~Ohlsson, and H.~Hogland.
\newblock Dc-resistivity mapping of internal landfill structures: two
  pre-excavation surveys.
\newblock {\em Environmental Geology}, 39(3):360--371, Jan 2000.

\bibitem[BDR97]{Binley1997}
Andrew Binley, William Daily, and Abelardo Ramirez.
\newblock Detecting leaks from environmental barriers using electrical current
  imaging.
\newblock {\em Journal of Environmental and Engineering Geophysics},
  2(1):11--19, Mar 1997.

\bibitem[BGS22]{Berrone2022}
Stefano Berrone, Denise Grappein, and Stefano Scial{\'o}.
\newblock 3d-1d coupling on non conforming meshes via a three-field
  optimization based domain decomposition.
\newblock {\em Journal of Computational Physics}, 448:110738, 2022.

\bibitem[BK05]{Binley2005}
Andrew Binley and Andreas Kemna.
\newblock {\em DC Resistivity and Induced Polarization Methods}, pages
  129--156.
\newblock Springer Netherlands, Dordrecht, 2005.

\bibitem[BP91]{Bonet1991}
Javier Bonet and Jaime Peraire.
\newblock An alternating digital tree (adt) algorithm for 3d geometric
  searching and intersection problems.
\newblock {\em International Journal for Numerical Methods in Engineering},
  31(1):1--17, 1991.

\bibitem[BSB{\etalchar{+}}20]{BLANCHY2020104423}
Guillaume Blanchy, Sina Saneiyan, Jimmy Boyd, Paul McLachlan, and Andrew
  Binley.
\newblock Resipy, an intuitive open source software for complex geoelectrical
  inversion/modeling.
\newblock {\em Computers \& Geosciences}, 137:104423, 2020.

\bibitem[CLZ19]{Cerroni2019}
Daniele Cerroni, Federica Laurino, and Paolo Zunino.
\newblock Mathematical analysis, finite element approximation and numerical
  solvers for the interaction of 3d reservoirs with 1d wells.
\newblock {\em GEM - International Journal on Geomathematics}, 10(1):4, Jan
  2019.

\bibitem[DC17]{DEDONNO2017302}
Giorgio {De Donno} and Ettore Cardarelli.
\newblock Tomographic inversion of time-domain resistivity and chargeability
  data for the investigation of landfills using a priori information.
\newblock {\em Waste Management}, 59:302--315, 2017.

\bibitem[DFL{\etalchar{+}}18]{DIMAIO2018629}
R.~{Di Maio}, S.~Fais, P.~Ligas, E.~Piegari, R.~Raga, and R.~Cossu.
\newblock 3d geophysical imaging for site-specific characterization plan of an
  old landfill.
\newblock {\em Waste Management}, 76:629--642, 2018.

\bibitem[DPC{\etalchar{+}}13]{DeCarlo2013}
Lorenzo {De Carlo}, Maria~Teresa Perri, Maria~Clementina Caputo, Rita Deiana,
  Michele Vurro, and Giorgio Cassiani.
\newblock Characterization of a dismissed landfill via electrical resistivity
  tomography and mise-{\'a}-la-masse method.
\newblock {\em Journal of Applied Geophysics}, 98:1--10, 2013.

\bibitem[DQ08]{DAngelo2008}
Carlo D'Angelo and Alfio Quarteroni.
\newblock On the coupling of 1d and 3d diffusion-reaction equations:
  application to tissue perfusion problems.
\newblock {\em Mathematical Models and Methods in Applied Sciences},
  18(08):1481--1504, 2008.

\bibitem[DRL10]{Dahlin2010}
T.~Dahlin, H.~Rosqvist, and V.~Leroux.
\newblock Resistivity-ip mapping for landfill applications.
\newblock {\em First Break}, 28(8), 2010.

\bibitem[DS12]{DAngelo2011}
Carlo D'Angelo and Anna Scotti.
\newblock A mixed finite element method for {D}arcy flow in fractured porous
  media with non-matching grids.
\newblock {\em Mathematical {M}odelling and {N}umerical {A}nalysis},
  46(02):465--489, 2012.

\bibitem[FK18]{Fumagalli2016a}
Alessio Fumagalli and Eirik Keilegavlen.
\newblock Dual virtual element method for discrete fractures networks.
\newblock {\em SIAM Journal on Scientific Computing}, 40(1):B228--B258, 2018.

\bibitem[FKS19]{Scialo2017}
Alessio Fumagalli, Eirik Keilegavlen, and Stefano Scial\`o.
\newblock Conforming, non-conforming and non-matching discretization couplings
  in discrete fracture network simulations.
\newblock {\em Journal of Computational Physics}, 376:694--712, 2019.

\bibitem[Fra97]{Frangos1997}
William Frangos.
\newblock Electrical detection of leaks in lined waste disposal ponds.
\newblock {\em Geophysics}, 62(6):1737--1744, 1997.

\bibitem[FS20]{Fumagalli2019a}
Alessio Fumagalli and Anna Scotti.
\newblock A multi-layer reduced model for flow in porous media with a fault and
  surrounding damage zones.
\newblock {\em Computational Geosciences}, 24(3):1347--1360, 2020.

\bibitem[GKN20]{Gjerde2020}
Ingeborg~G{\aa}seby Gjerde, Kundan Kumar, and Jan~Martin Nordbotten.
\newblock A singularity removal method for coupled 1d--3d flow models.
\newblock {\em Computational Geosciences}, 24(2):443--457, Apr 2020.

\bibitem[GKNW19]{Gjerde2019}
Ingeborg~G{\aa}seby Gjerde, Kundan Kumar, Jan~Martin Nordbotten, and Barbara
  Wohlmuth.
\newblock Splitting method for elliptic equations with line sources.
\newblock {\em ESAIM: M2AN}, 53(5):1715--1739, 2019.

\bibitem[GR09]{Geuzaine2009}
Christophe Geuzaine and Jean-Fran{\c{c}}ois Remacle.
\newblock Gmsh: A 3-d finite element mesh generator with built-in pre- and
  post-processing facilities.
\newblock {\em International Journal for Numerical Methods in Engineering},
  79(11):1309--1331, 2009.

\bibitem[JFT98]{Thompson1998}
Nigel P.~Weatherill Joe F.~Thompson, Bharat K.~Soni, editor.
\newblock {\em Handbook of Grid Generation}.
\newblock CRC Press, 1998.

\bibitem[KBF{\etalchar{+}}20]{Keilegavlen2020}
Eirik Keilegavlen, Runar Berge, Alessio Fumagalli, Michele Starnoni, Ivar
  Stefansson, Jhabriel Varela, and Inga Berre.
\newblock Porepy: An open-source software for simulation of multiphysics
  processes in fractured porous media.
\newblock {\em Computational Geosciences}, 2020.

\bibitem[KLMZ21]{Kuchta2021}
Miroslav Kuchta, Federica Laurino, Kent-Andre Mardal, and Paolo Zunino.
\newblock Analysis and approximation of mixed-dimensional pdes on 3d-1d domains
  coupled with lagrange multipliers.
\newblock {\em SIAM Journal on Numerical Analysis}, 59(1):558--582, 2021.

\bibitem[LDR10]{Leroux2010}
V.~Leroux, T.~Dahlin, and H.~Rosqvist.
\newblock Time-domain ip and resistivity sections measured at four landfills
  with different contents.
\newblock {\em Near Surface 2010}, P09, 2010.

\bibitem[Lok22]{Loke2022}
Meng~Heng Loke.
\newblock Tutorial : 2-d and 3-d electrical imaging surveys.
\newblock In {\em Geotomo Software}, 2022.

\bibitem[LRDC19]{Loke2019}
M.H. Loke, D.~Rucker, T.~Dahlin, and J.E. Chambers.
\newblock Recent advances in the geoelectrical method and new challenges: A
  software perspective.
\newblock {\em FastTimes}, 24(4):56--62, 2019.

\bibitem[LRQ{\etalchar{+}}19]{Ling2019}
C.~Ling, A.~Revil, Y.~Qi, F.~Abdulsamad, P.~Shi, S.~Nicaise, and L.~Peyras.
\newblock Application of the mise-{\`a}-la-masse method to detect the bottom
  leakage of water reservoirs.
\newblock {\em Engineering Geology}, 261:105272, Nov 2019.

\bibitem[Mej00]{MEJU2000115}
Maxwell~A. Meju.
\newblock Geoelectrical investigation of old/abandoned, covered landfill sites
  in urban areas: model development with a genetic diagnosis approach.
\newblock {\em Journal of Applied Geophysics}, 44(2):115--150, 2000.

\bibitem[Men18]{Menke2018}
William Menke.
\newblock {\em Geophysical Data Analysis: Discrete Inverse Theory}.
\newblock Academic Press, 4th edition, 2018.

\bibitem[MJR05]{Martin2005}
Vincent Martin, J{\'e}r{\^o}me Jaffr{\'e}, and Jean~Elisabeth Roberts.
\newblock Modeling {F}ractures and {B}arriers as {I}nterfaces for {F}low in
  {P}orous {M}edia.
\newblock {\em SIAM J. Sci. Comput.}, 26(5):1667--1691, 2005.

\bibitem[NBFK19]{Nordbotten2018}
Jan~Martin Nordbotten, Wietse Boon, Alessio Fumagalli, and Eirik Keilegavlen.
\newblock Unified approach to discretization of flow in fractured porous media.
\newblock {\em Computational Geosciences}, 23(2):225--237, 2019.

\bibitem[Pea78]{Peaceman1978}
Donald~W. Peaceman.
\newblock {Interpretation of Well-Block Pressures in Numerical Reservoir
  Simulation(includes associated paper 6988 )}.
\newblock {\em Society of Petroleum Engineers Journal}, 18(03):183--194, 06
  1978.

\bibitem[Pea83]{Peaceman1983}
Donald~W. Peaceman.
\newblock {Interpretation of Well-Block Pressures in Numerical Reservoir
  Simulation With Nonsquare Grid Blocks and Anisotropic Permeability}.
\newblock {\em Society of Petroleum Engineers Journal}, 23(03):531--543, 06
  1983.

\bibitem[RA71]{roy1971depth}
A~Roy and A~Apparao.
\newblock Depth of investigation in direct current methods.
\newblock {\em Geophysics}, 36(5):943--959, 1971.

\bibitem[R.D89]{BARKER1989}
Barker R.D.
\newblock Depth of investigation of collinear symmetrical four-electrode
  arrays.
\newblock {\em Geophysics}, 54:1031--1037, 1989.

\bibitem[RG11]{Ruecker2011}
Carsten R{\"u}cker and Thomas G{\"u}nther.
\newblock The simulation of finite ert electrodes using the complete electrode
  model.
\newblock {\em GEOPHYSICS}, 76(4):F227--F238, 2011.

\bibitem[RGW17]{RUCKER2017106}
Carsten R{\"u}cker, Thomas G{\"u}nther, and Florian~M. Wagner.
\newblock pygimli: An open-source library for modelling and inversion in
  geophysics.
\newblock {\em Computers \& Geosciences}, 109:106--123, 2017.

\bibitem[TVFT14]{Tsourlos2014}
P.~Tsourlos, G.N. Vargemezis, I.~Fikos, and G.N. Tsokas.
\newblock Dc geoelectrical methods applied to landfill investigation: case
  studies from greece.
\newblock {\em First Break}, 32(8), 2014.

\bibitem[VSP{\etalchar{+}}11]{VAUDELET2011738}
P.~Vaudelet, M.~Schmutz, M.~Pessel, M.~Franceschi, R.~Gu{\'e}rin, O.~Atteia,
  A.~Blondel, C.~Ngomseu, S.~Galaup, F.~Rejiba, and P.~B{\'e}gassat.
\newblock Mapping of contaminant plumes with geoelectrical methods. a case
  study in urban context.
\newblock {\em Journal of Applied Geophysics}, 75(4):738--751, 2011.

\bibitem[WT90]{Sheriff1990}
R.E.~Sheriff W.M.~Telford, L.P.~Geldart.
\newblock {\em Applied Geophysics}.
\newblock Cambridge University Press, 2nd edition, 1990.

\end{thebibliography}

\end{document}